\documentclass[11pt]{amsart}
\usepackage{verbatim}
\usepackage{hyperref}
\usepackage{mathdots}
\usepackage{amsmath}
\usepackage{enumerate}
\usepackage{amssymb}
\usepackage{color}
\usepackage{url}
\usepackage{graphicx}
\usepackage{fullpage}
\usepackage{lscape,ulem,esvect}
\usepackage[ruled,vlined]{algorithm2e}
\usepackage{longtable}
\usepackage{booktabs}
 \usepackage{caption}

\newcommand\+{\;\lower\plusheight\hbox{$+$}\;}
\newcommand\lldots{\;\lower\plusheight\hbox{$\cdots$}\;}

\newtheorem{Theorem}{Theorem}[section]
\newtheorem{Lemma}[Theorem]{Lemma}
\newtheorem{Corollary}[Theorem]{Corollary}
\newtheorem{Definition}[Theorem]{Definition}

\newtheorem{Remark}[Theorem]{Remark}

\newtheorem{Conjecture}[Theorem]{Conjecture}

\DeclareMathOperator{\spn}{Span}

\setbox0=\hbox{$+$}
\newdimen\plusheight
\plusheight=\ht 0 \setbox0=\hbox{$-$}
\newdimen\minusheight
\minusheight=\ht0

\setbox0=\hbox{$\cdots$}
\newdimen\cdotsheight
\cdotsheight=\plusheight

\title{Ramanujan Type Congruences for Quotients of Klein forms}
\author{Timothy Huber, Nathaniel Mayes, Jeffery Opoku and Dongxi Ye}

\begin{document}
\begin{abstract}
In this work, Ramanujan type congruences modulo powers of primes $p \ge 5$ are derived for a general class of products that are modular forms of level $p$. These products are constructed in terms of Klein forms and subsume generating functions for $t$-core
partitions known to satisfy Ramanujan type congruences for $p=5,7,11$. The vectors of exponents corresponding to products that are modular forms for $\Gamma_{1}(p)$ are subsets of bounded polytopes with explicit parameterizations. This allows for the derivation of a complete list of products that are modular forms for $\Gamma_{1}(p)$ of weights $1\le k \le 5$ for primes $5\le p \le 19$ and whose Fourier coefficients satisfy Ramanujan type congruences for all powers of the primes. For each product satisfying a congruence, cyclic permutations of the exponents determine additional products satisfying congruences.
Common forms among the exponent sets lead to products satisfying Ramanujan type congruences for a broad class of primes, including $p> 19$. Canonical bases for modular forms of level $5\le p \le 19$ are constructed by summing weight one Hecke Eisensten series of levels $5\le p \le 19$ and expressing the result as a quotient of Klein forms. Generating sets for the graded algebras of modular forms for $\Gamma_{1}(p)$ and $\Gamma(p)$ are formulated in terms of permutations of the exponent sets. A sieving process is described by decomposing the space of modular forms of weight $1$ for $\Gamma_{1}(p)$ as a direct sum of subspaces of modular forms for $\Gamma(p)$ of the form $q^{r/p}\mathbb Z[[q]]$. Since the relevant bases generate the graded algebra of modular forms for these groups, the weight one decompositions determine series dissections for modular forms of higher weight that lead to additional classes of congruences. 
\end{abstract}
\address{
School of Mathematical and Statistical Sciences, University of Texas Rio Grande
Valley, Edinburg, Texas 78539, USA}
\email{timothy.huber@utrgv.edu}
\email{npmayes@gmail.com}
\email{opokujeffery5@gmail.com}

\address{
School of Mathematics (Zhuhai), Sun Yat-sen University, Zhuhai 519082, Guangdong,
People's Republic of China}
\email{yedx3@mail.sysu.edu.cn}

\thanks{\textit{Statements and Declarations.} Dongxi Ye was supported by the Guangdong Basic and Applied Basic Research
Foundation (Grant No. 2023A1515010298).}

\maketitle
\numberwithin{equation}{section}
\allowdisplaybreaks

\section{Introduction}

In combinatorial theory, one of the most famous partition counting functions is the so-called $t$-core partition function $\mathcal{P}_{t}(n)$ first introduced by James and Kerber \cite{JK}, whose generating function is given by the simple form
$$
\sum_{n=0}^{\infty}\mathcal{P}_{t}(n)q^{n}=\frac{(q^{t};q^{t})^{t}_{\infty}}{(q;q)_{\infty}},
$$
where $(a;q)_{\infty} = \prod_{n=1}^{\infty} (1 -
aq^{n-1})$, and $(a,b;q)_{\infty}=(a;q)_{\infty}(b;q)_{\infty}$. The generating function for $\mathcal{P}_{t}(n)$ differs from a holomorphic modular form by an integral power of $q=e^{2\pi i\tau}$. In their work \cite{GKS}, Garvan et al prove that for $t=5,7$ and~11, $\mathcal{P}_{t}(n)$ satisfy the following Ramanujan type congruence 
\begin{align} \label{tcp}
\mathcal{P}_{t}(t^{j}n-\delta_{t})\equiv0\pmod{t^{j}}    
\end{align}
for any nonnegative integer $j$, where $\delta_{t}=\frac{t^{2}-1}{24}$. When \eqref{tcp} is specialized to $j=1$, Ramanujan's congruences \cite{R1} for the partition function may be deduced. Namely, if $\sum_{n=0}^{\infty} p(n)q^{n} = (q;q)_{\infty}^{-1}$,
$$
p(tn-\delta_{t})\equiv0\pmod{t}.
$$
Congruences such as these have been termed Ramanujan type congruences. A summary of many
arithmetic functions that satisfy
similar congruences can be found in \cite{Wa}. This work focuses on a class of infinite products satisfying Ramanujan type congruences defined for primes $p\ge 5$ by


\begin{align}
  \label{eq:p}
  (q^{p}; q^{p})_{\infty}^{a_{0}}\prod_{i=1}^{(p-1)/2}(q^{i}, q^{p-i};
      q^{p})_{\infty}^{a_{i}} := \sum_{n=0}^{\infty} P_{a_{0},\ldots,a_{(p-1)/2}}(n)q^{n}, \qquad a_{i} \in \Bbb Z.
\end{align}
For certain choices of exponents $a_{i}$, the coefficients of these infinite products satisfy analogous classes of Ramanujan type congruences defined by Definition \ref{defcong}.  
Like the subsumed generating functions for enumerating $t$-core partitions $\mathcal{P}_{t}(n) = P_{t-1,-1,\ldots,-1}(n)$, the products \eqref{eq:p} differ from modular forms by a power of $q$ for restricted values of the exponents. This instantiation naturally motivates one to seek Ramanujan type congruences satisfied by $P_{a_{0},\ldots,a_{(p-1)/2}}(n)$ in cases where the generating function is modular. We show that if $\ell=\ell(a_{0},\ldots,a_{(p-1)/2})\in\mathbb{Z}$ is defined by
  \begin{align} \label{lp}
     \ell\left(a_{0},\ldots,a_{(p-1)/2}\right)
  =\frac{p}{24 }\left (a_{0}+2\sum_{i=1}^{(p-1)/2}a_{i} \right )
  -\frac{p}{2}\sum_{i=1}^{(p-1)/2}\left \langle\frac{i}{p}
  \right \rangle\left(1-\left \langle\frac{i}{p} \right \rangle\right)
  a_{i},
  \end{align}
  where $\langle r\rangle=r-\lfloor r\rfloor$ denotes the fractional
  part of~$r$, then for a fixed positive {\it even} integer $a_{0}$, 
  \begin{align}
\label{faaa}
  &f_{a_{0},\ldots,a_{(p-1)/2}}(\tau)
    :=   \sum_{n=\ell}^{\infty}
    P_{a_{0},\ldots,a_{(p-1)/2}}\left(n-\ell\right)q^{n}=q^{\ell} (q^{p}; q^{p})_{\infty}^{a_{0}}\prod_{i=1}^{(p-1)/2}(q^{i}, q^{p-i};
      q^{p})_{\infty}^{a_{i}}
\end{align}
is a holomorphic modular form of weight $a_{0}/2$ for $\Gamma_{1}(p)$ if $(a_{0},\ldots,a_{(p-1)/2})$ satisfy Lemma~\ref{m1r}. The number of such forms will be shown to be finite. For fixed values of $a_{0}$, we identify all lattice points $(a_{0},\ldots,a_{(p-1)/2})$ of the polytope determined by Lemma~\ref{m1r} such that $f_{a_{0},\ldots,a_{(p-1)/2}}(\tau)$ is a holomorphic modular form for $\Gamma_{1}(p)$ whose Fourier coefficients satisfy a Ramanujan type congruence
$$
P_{a_{0},\ldots,a_{(p-1)/2}}(p^{j}n-\ell)\equiv0\pmod{p^{j}}.
$$

\begin{Remark} \label{defcw}
Although combinatorial interpretations of the coefficients $P_{a_{0},\ldots,a_{(p-1)/2}}(n)$ do not play a role in the present work, the coefficients count colored partitions according to weights
determined by the parity of the number of parts in congruence
classes modulo $p$.

In particular, for $a_{0} \in 2\Bbb Z^{+}$, and $(a_{0},\ldots,a_{(p-1)/2})\in\mathbb{Z}^{\frac{p+1}{2}}$, the coefficient $P_{a_{0},\ldots,a_{(p-1)/2}}(n)$ 
associated to
$(a_{0},\ldots,a_{(p-1)/2})$ provides a weighted enumeration of
  partitions $\lambda$ of $n$ such that:
 \begin{itemize}
 \item If $a_{i}\le 0 $, parts congruent to $i,
   p-i \pmod{p}$ appear with $|a_{i}|$ colors.
\item If $a_{i} >0$, parts congruent to
  $i,p-i \pmod{p}$ are $a_{i}$-colored,
    distinct, and $\lambda$ is weighted by $\prod_{\{i\mid a_{i} >0\}} (-1)^{\kappa_{i}(\lambda)}$, where $\kappa_{i}(\lambda)$ is the number of parts congruent to $i,p-i\pmod{p}$.
  \end{itemize}
\end{Remark}

In the following, a precise definition of Ramanujan type congruences is given. The adopted definition of Ramanujan type congruences is augmented by including 
congruences that hold for only finitely many powers of primes termed \textit{chimeral} that may be representative of larger classes.

\begin{Definition}
\label{defcong}
Let $f(\tau)=\sum_{n=0}^{\infty}a(n)q^{n}$ and $a(n)\in\mathbb{Z}$ be a holomorphic modular
form for $\Gamma \le SL_{2}(\Bbb{Z})$. For a prime $p$, we say that $f(\tau)$ satisfies a Ramanujan
type congruence modulo~$p$ if $a(p^{j}n)\equiv 0\pmod{p^{\alpha j}}$
for a positive integer $\alpha$ and for all $n,j\geq1$. If the congruences are satisfied for $1\le
j\le r$ but fail to be satisfied for $j=r+1$, then we
call the congruences chimeral modulo $p^{r}$, or modulo~$p$ of order~$r$.
\end{Definition}
An important feature of our work is the algebra of Ramanujan type congruences for the products under consideration. Each product satisfying a Ramanujan congruence gives rise to a set of products satisfying congruences whose subscripts are related by permutations from $S_{(p-1)/2}$. To this end, here and throughout the rest of the work, for the fixed minimal positive $\alpha$ for which  $(\mathbb{Z}/p\mathbb{Z})^{\times}/\{\pm1\}=\langle\alpha\rangle$, we define the permutation $\sigma_{p}$
 \begin{equation} 
 \label{sigmap}
 \sigma_{p}=\left(1,[\alpha],[\alpha^{2}],\ldots,[\alpha^{\frac{p-3}{2}}]\right)\in S_{(p-1)/2},
 \end{equation}
 where $[m]\in (\mathbb{Z}/p\mathbb{Z})^{\times}/\{\pm1\}$ denotes the least positive representative so that $[m]\in\{1,\ldots,(p-1)/2\}$, and $\sigma_{p}$ viewed as a cycle acts on a vector of length $(p+1)/2$, namely $(a_{0}, a_{1}, \ldots, a_{(p-1)/2})$, by omitting the first component: $$\sigma_{p} (a_{0}, a_{1}, \ldots, a_{(p-1)/2}) = (a_{0}, a_{\sigma_{p}(1)}, a_{\sigma_{p}(2)}, \ldots,a_{\sigma_{p}((p-1)/2)}).$$
 We also write $\sigma_{p}(f_{a_{0},\ldots,a_{(p-1)/2}})$ for $f_{\sigma_{p}(a_{0},\ldots,a_{(p-1)/2})}$, and linearly extend action of the permutation on linear combinations of $f_{a_{0},\ldots,a_{(p-1)/2}}$. For primes $5 \le p \le 19$, we show that a product satisfies a Ramanujan type congruence modulo $p$ if and only if each product in its orbit under $\sigma_{p}$ satisfies the congruences. This applies to products satisfying standard or chimeral Ramanujan type congruences.

 \begin{Theorem} 
 \label{z1}
 Let $5\leq p\leq 19$ be prime, let $\sigma_{p}$ be defined by \eqref{sigmap}, and let $f_{a_{0}, \ldots, a_{(p-1)/2}}(\tau)$ be a modular form for $\Gamma_{1}(p)$, defined as in \eqref{faaa}.
 Then $f_{a_{0},  \ldots, a_{(p-1)/2}}(\tau)$ satisfies a Ramanujan type congruence modulo $p$ if and only if  $f_{\sigma_{p}(a_{0},  \ldots, a_{(p-1)/2})}(\tau)$ satisfies a Ramanujan type congruence. 
\end{Theorem}



Theorem \ref{mainth1} provides a complete list of products $f_{a_{0}, a_{1}, \ldots, a_{(p-1)/2}} \in M_{a_{0}/2}(\Gamma_{1}(p))$ satisfying Ramanujan type congruences for weights $2\leq a_{0}\leq 10$ and primes $5\le p \le 19$. These represent a larger class of quotients of the form \eqref{faaa} satisfying Ramanujan type congruences for larger primes.

\begin{Theorem}\label{mainth1}
Let $a_{0}$ be an even integer such that $2\le a_{0}\le 10$.
\begin{enumerate}[(1)]
\item Then 
$$\sum_{n=\ell}^{\infty} P_{a_{0},\ldots,a_{2}}(n - \ell)q^{n} \in M_{a_{0}/2}(\Gamma_{1}(5))$$ and satisfies $P_{a_{0},\ldots,a_{2}}(5^{k}n-\ell)\equiv0\pmod{5^{k}}$ if and only if
\begin{equation*}
    (a_{0},a_{1},a_{2})\in\left\{\begin{array}{cccc}
    (4,-1,-1),\sigma_{5}^{j}(6,1,-4), \sigma_{5}^{j}(8,3,-7), (8, -2, -2), (8, 4, 4), \\  \sigma_{5}^{j}(10,5,-10), \sigma_{5}^{j}(10,0,-5), \sigma_{5}^{j}(10,1,6)
    \end{array}\right\}_{0\le j \le 1}.
    \end{equation*}

\item Define $\sigma_{7}=(1,2,3)$. Then 
$$\sum_{n=\ell}^{\infty} P_{a_{0},\ldots,a_{3}}(n - \ell)q^{n} \in M_{a_{0}/2}(\Gamma_{1}(7))$$ and satisfies $P_{a_{0},\ldots,a_{3}}(7^{k}n-\ell)\equiv0\pmod{7^{k}}$ if and only if  
\begin{equation*}    (a_0, a_1, a_2, a_3) \in\left\{\begin{array}{c}     \sigma_{7}^{j}(4,1,-1,-2),\sigma_{7}^{j}(6,1,0,-4),\sigma_{7}^{j}(6,2,-3,-2), \sigma_{7}^{j}(6,-2,-5,4),  \\   (6,-1,-1,-1),(6,3,3,3),       \sigma_{7}^{j}(8,7, -3, -8), \sigma_{7}^{j}(8, 4, -1, -7), \\ \sigma_{7}^{j}(8, 1, 1, -6), \sigma_{7}^{j}(8, 3, -5, -2), \sigma_{7}^{j}(8, 5, -4, -5), \sigma_{7}^{j}(8, 2, -2, -4), \\ \sigma_{7}^{j}(8, 0, -3, -1), \sigma_{7}^{j}(8, 0, 6, 2), \sigma_{7}^{j}(8, 1, 3, 4)     \end{array}\right\}_{0\le j \le 2}.     \end{equation*}
Moreover, for $(a_{0},\ldots,a_{3})=\sigma_{7}^{j}(6,-4,1,0)$, $P_{a_{0},\ldots,a_{3}}(7^{k}n-\ell)\equiv0\pmod{7^{2k}}$. 

\item Define $\sigma_{11}=(1,2,4,3,5)$. Then
$$\sum_{n=\ell}^{\infty} P_{a_{0},\ldots,a_{5}}(n - \ell)q^{n} \in M_{a_{0}/2}(\Gamma_{1}(11))$$ and satisfies $P_{a_{0},\ldots,a_{5}}(11^{k}n-\ell)\equiv0\pmod{11^{k}}$ if and only if  
$$
 (a_{0},a_{1},a_{2},a_{3},a_{4},a_{5})\in\left\{\begin{array}{cccc}\sigma_{11}^{j}(4,-2,0,0,1,-1),\sigma_{11}^{j}(4,-2,0,1,-2,1),\sigma_{11}^{j}(6,-4,1,0,0,0),\\ (8,4,4,4,4,4),
      (10,-1,-1,-1,-1,-1)
                                                \end{array}\right\}_{0\leq j\leq 4}.
$$
Moreover, for $(a_{0},\ldots,a_{5})=\sigma_{11}^{j}(6,-4,1,0,0,0)$, $P_{a_{0},\ldots,a_{5}}(11^{k}n-\ell)\equiv0\pmod{11^{2k}}$. 
    \item Define $\sigma_{13}=(1,2,4,5,3,6)$. Then  $$\sum_{n=\ell}^{\infty} P_{a_{0},\ldots,a_{6}}(n - \ell)q^{n} \in M_{a_{0}/2}(\Gamma_{1}(13))$$ and satisfies $P_{a_{0},\ldots,a_{6}}(13^{k}n-\ell)\equiv0\pmod{13^{k}}$ if and only if  
   \begin{align*}
   ( a_{0},  \ldots, a_{6}) \in \left \{ \begin{array}{c} \sigma_{13}^{j}(4,-2, 0, 0, 1, -2, 1), \sigma_{13}^{j}(4,-2, 0, 0, 0, 1, -1), \\ 
       \sigma_{13}^{j} (4,-2, 1, -2, 1, 0, 0), \sigma_{13}^{j}(6,1, 0, 0, 0, 0, -4), 
       \\ \sigma_{13}^{j}(6,-3,-1,2,-3,2,0) 
    \end{array}  \right \}_{ 0 \le j \le 5}.
        \end{align*}
     Moreover,  for  $( a_{0},  \ldots, a_{6})= \sigma_{13}^{j}(6,1,0,0,0,0,-4)$, $P_{a_{0},\ldots,a_{6}}(13^{k}n-\ell)\equiv0\pmod{13^{2k}}$.
       
        \item Define $\sigma_{17}=(1,3,8,7,4,5,2,6)$. Then 
        $$\sum_{n=\ell}^{\infty} P_{a_{0},\ldots,a_{8}}(n - \ell)q^{n} \in M_{a_{0}/2}(\Gamma_{1}(17))$$ and satisfies $P_{a_{0},\ldots,a_{8}}(17^{k}n-\ell)\equiv0\pmod{17^{k}}$ if and only if  
   \begin{align*}
   ( a_{0},  \ldots, a_{8}) \in \left \{ \begin{array}{c}  \sigma_{17}^{j}(4,-2,0,0,0,0,1,-2,1),  \sigma_{17}^{j}(4,-2,0,0,0,0,0,1,-1)),  \\ \sigma_{17}^{j}(4,-2, 1, -2, 1, 0, 0, 0, 0), \sigma_{17}^{j}(4,-2,0,1,-2,1,0,0,0), \\ \sigma_{17}^{j}(6,1, 0, 0, 0, 0, 0,0,-4), 
   \end{array}  \right \}_{0 \le j \le 7}.
        \end{align*} 
       Moreover, for $( a_{0},  \ldots, a_{8})= \sigma_{17}^{j}(6,1, 0, 0, 0, 0, 0,0,-4)$, $P_{a_{0},\ldots,a_{8}}(17^{k}n - \ell)\equiv0\pmod{17^{2k}}$.
        
        \item Define $\sigma_{19}=(1,2,4,8,3,6,7,5,9)$. Then 
        $$\sum_{n=\ell}^{\infty} P_{a_{0},\ldots,a_{9}}(n - \ell)q^{n} \in M_{a_{0}/2}(\Gamma_{1}
        (19))$$ and satisfies $P_{a_{0},\ldots,a_{9}}(19^{k}n - \ell)\equiv0\pmod{19^{k}}$ if and only if  
   \begin{align*}
   ( a_{0}, \ldots, a_{9}) \in \left \{ \begin{array}{c}  \sigma_{19}^{i}(4,-2,0,0,0,0,0,1,-2,1),  \sigma_{19}^{i}(4,-2,0,0,0,0,0,0,1,-1),\\
                      \sigma_{19}^{i}(4,-2, 1, -2, 1, 0, 0, 0, 0, 0),\sigma_{19}^{i}(4,-2, 0, 1, -2, 1, 0, 0, 0, 0), \\ \sigma_{19}^{j}(6,1, 0, 0,0, 0, 0, 0,0,-4)
                    \end{array}  \right \}_{0 \le j \le 8}.
        \end{align*} 
        Moreover, for $( a_{0},  \ldots, a_{9})= \sigma_{19}^{j}(6,1,0, 0, 0, 0, 0, 0,0,-4)$, $P_{a_{0},\ldots,a_{9}}(19^{k}n - \ell)\equiv0\pmod{19^{2k}}$.
\end{enumerate}
\end{Theorem}
In particular if $\vv{\mathbf{a}}_{k}$ is the vector with entries
$a$ of length $k$, then the coefficients of 
\begin{align}
  \label{eq:7}
\sum_{n=0}^{\infty} P_{(p-1),\vv{\mathbf{(-1)}}_{(p-1)/2}}(n)q^{n} =
  \frac{(q^{p}, q^{p})_{\infty}^{p}}{(q;q)_{\infty}}
\end{align}
count the number of $p$-core partitions; 
$P_{6,3,3,3}(n)$ are the Fourier coefficients of the CM
modular form $\eta(\tau)^{3}\eta(7\tau)^{3}$; 
and $P_{8,4,4,4,4,4}(n)$ are the Fourier coefficients of the eta product $\eta(\tau)^{4}\eta(11\tau)^{4}$.

The exploration of products of small weight and primes points to the existence of Ramanujan type congruences for products beyond the range of prime levels considered in here. Based on the common forms of the products appearing in Theorem \ref{mainth1}, certain products are identified that may satisfy Ramanujan type congruences for all sufficiently large primes. These products include 
\begin{align*}
 \frac{(q^{p};q^p)_{\infty}^{6}(q,q^{p-1};q^{p})_{\infty}}{(q^{(p-1)/2},q^{(p+1)/2};q^{p})_{\infty}^{4}} &=  \sum_{n=0}^{\infty} P_{6,1,0,...,0,-4}(n)q^{n}
\end{align*}
and  \begin{align*}
     \frac{\left (q^{(p-5)/2}, q^{(p+5)/2}, q^{(p-1)/2}, q^{(p+1)/2};q^{p} \right )_{\infty} \left (q^{p}; q^{p} \right )_{\infty}^{4}}{ \left (q, q^{p-1}, q^{(p-3)/2}, q^{(p+3)/2};q^{p}\right)_{\infty}^{2}} = \sum_{n=0}^{\infty} P_{4,-2,0,\ldots,0,1,-2,1}(n)q^{n}.
 \end{align*}
These are among a larger class that appear to satisfy Ramanujan type congruences for sufficiently large primes, which we prove for primes $p\le 101$ through Theorem \ref{61004} and conjecture to be true for all larger primes. The congruences in Corollary \ref{square} are proven through Theorem \ref{61004} by showing that the corresponding products $f_{\vec{a}}$ are eigenfunctions for the linear operator $U_{p}:= U_{p,0}$, where
\begin{equation}
\label{defu}
U_{p,r}(f)=\sum_{n\equiv r\pmod{p}}b(n)q^{n/p}  \quad \text{for }
f=f(\tau)=\sum_{n=0}^{\infty}b(n)q^{n} \in
M_{k}(\Gamma_{1}(p)).
\end{equation}
Throughout the rest of this work, modular forms refer to holomorphic modular forms, and $M_{k}(\Gamma)$ denotes the space of holomorphic modular forms of weight~$k$ for a subgroup $\Gamma$ of ${\rm SL}_{2}(\mathbb{R})$.

\begin{Theorem}
\label{61004}
Suppose $\ell_{j}=\ell(\sigma_{p}^{j}(\vec{a}_{p}))$, and define $b_{j}(n)$ by 
 $$f_{\sigma_{p}^{j}(\vec{a}_{p})}:=\sum_{n=0}^{\infty}b_{j}(n)q^{n}.$$ 
\begin{enumerate}
    \item Let $p\geq5$ be a prime and $\vec{a}_{p}=(6,1,0,\ldots,0,-4)$. Then $$U_{p}(f_{\vec{a}_{p}})=p^{2}f_{\vec{a}_{p}}$$  if and only if $b_{j}(ps)=0$ for $s<\ell_{j}$ and $2\leq j\leq\frac{p-1}{2}$.
    \item Let $p\geq11$ be a prime and let
 $$\vec{a}_{p}\in \{(4,-2,0,\ldots 0,1,-2,1),(4,-2,0,\ldots,0,1,-1),(4,1,1,0,\ldots,0,-2,-2)\}.
 $$ Then  $$U_{p}(f_{\vec{a}_{p}})=pf_{\vec{a}_{p}}$$ if and only if
 $b_{j}(ps)=0$ for $s<\ell_{j}$ and $2\leq j\leq\frac{p-1}{2}$.
\end{enumerate}
 
\end{Theorem}
In addition to the products above satisfying Ramanujan type congruences, common exponent vectors are identified in Conjecture \ref{hl1} for products that appear to satisfy chimeral Ramanujan type congruences for sufficiently large primes. Our proofs of congruences use the fact that $U_{p,0}$ is a linear transformation on $M_{k}(\Gamma_{1}(p))$ and $U_{p,r}$ is a linear transformation from $M_{k}(\Gamma_{1}(p))$ to $M_{k}(\Gamma(p))$. It is known that the graded ring of modular forms on $\Gamma(p)$ and $\Gamma_{1}(p)$ for $5\le p \le 19$ is generated by weight one forms. We provide a canonical construction of the generating sets for the graded rings of modular forms in terms of quotients of Klein forms. 
  We identify vectors of weight one forms $\vec{v}_{r} \subset
    q^{r/p}\Bbb Z[[q]]$ for the primes $5 \le p \le 19$ that provide a $p$-dissection of the vector space
    \begin{align} \label{d1}
      {M}_{1}(\Gamma(p)) = \bigoplus_{r=0}^{p-1}
      \spn_{\Bbb C} \vec{v}_{r}.
    \end{align}
    Section \ref{Spdissect} provides $p$-dissections over $M_{1}(\Gamma(p))$ for generators of the graded algebra $M(\Gamma_{1}(p))$, $5 \le p \le 19$. In section 6, we consider the graded algebras generated by these weight one forms and apply a sifting procedure based on the orders at the cusp $\infty$ to formulate congruences satisfied by the coefficients of images of modular forms under dissection operators  $U_{p,r}$. 

The remainder of the paper is organized as follows. In Section \ref{fundamentals}, we introduce Klein forms and other basic building blocks for the products considered here. We prove that the products defined in \eqref{faaa} belong to $M_{a_{0}/2}(\Gamma_{1}(p))$ if and only if the exponents are lattice points in $\Bbb Z^{(p-2)/2}$ satisfying a set of linear congruences in a prescribed polyhedron. The corresponding polyhedron are shown to be bounded. We provide parameterizations for the polytopes that allow us to construct the number of products from \eqref{faaa} that are modular with respect to $\Gamma_{1}(p)$. In Section 3, by considering quotients of Klein forms that arise as simple sums of Hecke eisenstein series, explicit bases are constructed for weight $1$ modular forms on $\Gamma_{1}(p)$ that generate all modular forms of integer weights. We deduce how the products are permuted by $\sigma_{p}$ and trace the image of products $f_{\vec{a}} \in q^{r/p} \Bbb Z[[q]]$ under $\sigma_{p}$ and the operators $U_{p,r}$. This allows us to formulate a sufficient condition on the eigendecomposition of $f_{\vec{a}}$ with respect to $U_{p,0}$ in Lemma \ref{upeee} for $f_{\vec{a}}$ to satisfy a Ramanujan type congruence. The Ramanujan type congruences in Theorem \ref{mainth1} are proven in Section \ref{eigen} by using the eigendecompositions. In Section \ref{Spdissect}, we determine explicit decompositions of weight one forms for $\Gamma_{1}(p)$ into the subspaces $$
M_{1}(\Gamma(p))=\bigoplus_{r=0}^{p-1}M_{1,r}(\Gamma(p)),
$$
where for $0 \le r \le p-1$,
\begin{align} \label{mkr}
    M_{k,r}(\Gamma(p))=\{g\in M_{k}(\Gamma(p)):\,\,g\in q^{\frac{r}{p}}\mathbb{Z}[[q]]\}.
\end{align}
Since weight $1$ forms for $\Gamma(p)$ generate the graded algebra of modular forms for $\Gamma(p)$, this allows us to provide a constructive proof of dissection congruences for any weight resulting from application of $U_{p,r}$. In Section 6, we list applications of these dissection formulas to find additional congruences. In Section 7, we prove Theorem \ref{61004} and consider relevant aspects such as the size of the orbits of products satisfying Ramanujan type congruences and chimeral Ramanujan type congruencess of each level and weight. Section 8 concludes by discussing computational aspects of the paper, including algorithms applied to determine products satisfying the congruences.


{\bf Acknowledgment} The authors thank Prof. François Brunault, Prof. Winnie Li, Prof. Abhishek Saha and Dr. Wei-Lun Tsai for helpful discussions. They would also like to thank the anonymous referee for his/her comments, suggestions and corrections.

\section{Polytope Construction} \label{fundamentals}
 For certain values of $a_{0}, a_{1}, \ldots, a_{(p-1)/2}$, the functions $f_{a_{0}, a_{1}, \ldots, a_{(p-1)/2}}$ defined in \eqref{faaa} are modular forms of weight $a_{0}/2$ for $\Gamma_{1}(p)$ that may be represented in terms of the Dedekind eta function $\eta(\tau)=q^{1/24}(q;q)_{\infty}$ and Klein forms, defined for $z=(Q_{1},Q_{2})\in\mathbb{Q}^{2}-\mathbb{Z}^{2}$ and $q_{z}=e^{2\pi i(Q_{1}\tau+Q_{2})}$ by
$$
K_{(Q_{1},Q_{2})}(\tau)=e^{\pi iQ_{2}(Q_{1}-1)}q^{\frac{1}{2}Q_{1}(Q_{1}-1)}(1-q_{z})\prod_{n=1}^{\infty}(1-q_{z}q^{n})(1-q_{z}^{-1}q^{n})(1-q^{n})^{-2},
$$
whose basic properties are summarized as follows.

\begin{Lemma}[\cite{KL}]\label{klf}
Let $K_{(Q_{1},Q_{2})}(\tau)$ be defined as above. Then the following assertions hold.
\begin{enumerate}[(1)]

\item For $(Q_{1},Q_{2})\in\mathbb{Q}^{2}-\mathbb{Z}^{2}$ and $(s_{1},s_{2})\in\mathbb{Z}^{2}$, one has
  \begin{align*}
    K_{(-Q_{1},-Q_{2})}(\tau)&=-K_{(Q_{1},Q_{2})}(\tau) \\
    K_{(Q_{1},Q_{2})+(s_{1},s_{2})}(\tau) &=(-1)^{s_{1}s_{2}+s_{1}+s_{2}}e^{-\pi i(s_{1}Q_{2}-s_{2}Q_{1})}K_{(Q_{1},Q_{2})}(\tau).
  \end{align*}
\item For any $\begin{pmatrix}a&b\\ c&d\end{pmatrix}\in {\rm SL}_{2}(\mathbb{Z})$, 
  \begin{align*}
K_{(Q_{1},Q_{2})}\left(\frac{a\tau+b}{c\tau+d}\right)=(c\tau+d)^{-1}K_{(Q_{1}a+Q_{2}c,Q_{1}b+Q_{2}d)}(\tau).     
  \end{align*}
\item 
The order of vanishing of $K_{(Q_{1},Q_{2})}(\tau)$ at the cusp $i\infty$ is given by
\begin{align*}
{\rm ord}_{\infty}(K_{(Q_{1},Q_{2})})=\frac{1}{2}\langle Q_{1}\rangle\left(\langle Q_{1}\rangle -1\right),   
\end{align*}
where $\langle r\rangle=r-\lfloor r\rfloor$ denotes the fractional part of~$r$.
\end{enumerate}
\end{Lemma}
We work with the special values 
\begin{align}
   K_{p,i}=K_{p,i}(\tau)= K_{\left(i/p,0\right)}(p\tau): = q^{\frac{i (i-p)}{2 p}} \frac{(q^{i}, q^{p-i}; q^{p})_{\infty}}{(q^{p}; q^{p})_{\infty}^{2}}, \quad 1 \le i \le \frac{p-1}{2},
\end{align}
denoted subsequently in this work by $K_{p,i}=K_{p,i}(\tau)$. One can easily verify that
\begin{align}
\label{fak}
f_{a_{0},\ldots,a_{(p-1)/2}}(\tau)=\eta(p\tau)^{a_{0}+2\sum_{i=1}^{(p-1)/2}a_{i}}\prod_{i=1}^{(p-1)/2}K_{(i/p,0)}(p\tau)^{a_{i}}.
\end{align}
Lemma \ref{klf} allows us to deduce a characterization for $f_{a_{0},\ldots,a_{(p-1)/2}}(\tau)$ to be a modular form of integral weight for
$\Gamma_{1}(p)$. 

\begin{Lemma} \label{m1r}
Given a prime $p\geq5$, 
suppose that $\left(\mathbb{Z}/p\mathbb{Z}\right)^{\times}/\{\pm1\}=\langle\alpha\rangle$. Then 
$f_{a_{0},\ldots,a_{(p-1)/2}}(\tau)$ lies in $M_{a_{0}/2}(\Gamma_{1}(p))$ if and only if
\begin{align} \label{ineq}
\begin{cases} \displaystyle 
    0\leq  a_{0}+2\sum_{i=1}^{(p-1)/2}a_{i} \equiv 0 \pmod{24},\\\\
   \displaystyle   \frac{p}{2}\sum_{i=1}^{(p-1)/2}\left \langle\frac{i}{p}
  \right \rangle\left(1-\left \langle\frac{i}{p} \right \rangle\right)
  a_{i}\in\mathbb{Z},\\\\
    \displaystyle  \frac{p}{24 }\left (a_{0}+2\sum_{i=1}^{(p-1)/2}a_{i} \right )
      -\frac{p}{2}\sum_{i=1}^{(p-1)/2}\left \langle\frac{i\alpha^{m}}{p}
      \right \rangle\left(1-\left \langle\frac{i\alpha^{m}}{p} \right \rangle\right)
  a_{i}\geq0&\mbox{for all $0\leq m\leq \frac{p-3}{2}$.}
\end{cases}
\end{align}
\end{Lemma}

\begin{proof}
For an odd prime $p$,
\begin{align} \label{eq:77}
  -\sum_{i=1}^{(p-1)/2}\frac{p}{2}\left 
  \langle\frac{i}{p} \right \rangle\left(1-\left \langle\frac{i}{p}
  \right \rangle\right)a_{i}\in\mathbb{Z}
\end{align}
if and only if
\begin{align}
  \label{eq:8}
\sum_{i=1}^{(p-1)/2} i^{2}a_{i} \equiv 0
  \pmod{p}.  
\end{align}
Therefore \eqref{eq:77}  is true if and only if the quotient of Klein forms $
  \prod_{i=1}^{(p-1)/2}K_{(i/p,0)}(p\tau)^{a_{i}}$ satisfies the
  criteria for weak modularity on $\Gamma(p)$ derived in \cite[p.~68]{KL}. In
  fact, each factor $K_{(i/p,0)}(p\tau)^{a_{i}}$ is a weakly modular
  form for $\Gamma(p)$ of weight $-a_{i}$. To see this, let
  $\left  ( \begin{smallmatrix} a&b\\c&d\end{smallmatrix} \right ) \in\Gamma(p)$, and note that
  \begin{align*}
  &K_{\left(\frac{i}{p},0\right)}\left(p\frac{a\tau+b}{c\tau+d}\right) = K_{\left(\frac{i}{p},0\right)}\left(\frac{a(p\tau)+bp}{\frac{c}{p}(p\tau)+d}\right)=(c\tau+d)^{-1}K_{\left(\frac{ia}{p},ib\right)}(p\tau).
  \end{align*}
  Suppose that $a=1+Ap$ and $b=Bp$. Then
  $$
  K_{\left(\frac{ia}{p},ib\right)}(p\tau)=(-1)^{i^{2}Ab+iA+ib}(-1)^{i^{2}B}K_{\left(\frac{i}{p},0\right)}(p\tau)=(-1)^{i^{2}(Ab+B)+i(A+b)}K_{\left(\frac{i}{p},0\right)}(p\tau).
  $$
  Since $ad-bc=1$, then if $a$ even, then $A$ is odd and $B$ is odd, and thus $Ab+B$ and $A+b$ are even; also, if $a$ is odd, then $A$ is even, and so $i^{2}(Ab+B)$ and $i(A+b)$ have the same parity since $p$ is an odd prime. These show that $i^{2}(Ab+B)+i(A+b)$ must be even, and therefore
   $$
  K_{\left(\frac{ia}{p},ib\right)}(p\tau)=K_{\left(\frac{i}{p},0\right)}(p\tau). 
  $$
Hence,  $K_{(\frac{i}{p},0)}(p\tau)$ is weakly modular of  weight $-1$
  on $\Gamma(p)$. 
    Note that $\Gamma(p)$ is a normal subgroup of $\Gamma_{1}(p)$ with coset
  representatives $\left (\begin{smallmatrix}1&a\\0&1\end{smallmatrix} \right
)$ for
  $a=0,\ldots, p-1$. Thus, if $
  \prod_{i=1}^{(p-1)/2}K_{\left(\frac{i}{p},0\right)}(p\tau)^{a_{i}}$ is
  weakly modular on $\Gamma_{1}(p)$, it is enough for one to have that
  \begin{align*}
   \prod_{i=1}^{(p-1)/2}K_{(\frac{i}{p},0)}(p(\tau+1))^{a_{i}}&=e^{2\pi i\left(-\sum_{i=1}^{(p-1)/2}\frac{p}{2}\left 
  \langle\frac{i}{p} \right \rangle\left(1-\left \langle\frac{i}{p}
  \right \rangle\right)a_{i}\right)}q^{-\sum_{i=1}^{(p-1)/2}\frac{p}{2}\left 
  \langle\frac{i}{p} \right \rangle\left(1-\left \langle\frac{i}{p}
  \right \rangle\right)a_{i}}+\cdots \\ &=\prod_{i=1}^{(p-1)/2}K_{(\frac{i}{p},0)}(p\tau)^{a_{i}},
  \end{align*}
  and therefore \eqref{eq:7} is satified. In fact, this is enough for
  $\prod_{i=1}^{(p-1)/2}K_{(i/p,0)}(p\tau)^{a_{i}}$ to
  be weakly  modular on $\Gamma_{1}(p)$ since by definition one can see that
  $$
  \prod_{i=1}^{(p-1)/2}K_{\left(\frac{i}{p},0\right)}(p\tau)^{a_{i}}=q^{ -\sum_{i=1}^{(p-1)/2}\frac{p}{2}\left
  \langle\frac{i}{p} \right \rangle\left(1-\left \langle\frac{i}{p}
  \right \rangle\right)a_{i}}\left(1+q\mathbb{Z}[[q]]\right).
  $$ The weak modularity of $f_{a_{0}, \ldots, a_{(p-1)/2}}(\tau)$
  follows from the preceding discussion and the fact that the
  multiplier of $\eta(\tau)$ is a $24$-th root of unity. To prove the Corollary, it remains to show holomorphicity at  each cusp. 
 A set of inequivalent cusps for $X_{1}(p)$ is $$\{1/p, 2/p,\cdots (p-1)/(2p),1, 1/2,
1/3,\ldots, 1/((p-1)/2)\}.$$ 
By Lemma~\ref{klf}, one can show that the divisor of $K_{p,i}^{-1}$ is well defined on $X_{1}(p)$, and that the orders are supported by the cusps $[a/p]$, $1\le a \le (p-1)/2$. In particular, 
\begin{align*}
  {\rm Div}_{X_{1}(p)}(K_{p,i}^{-1})
  &=\sum_{a=1}^{(p-1)/2}\frac{p}{2}\left
    \langle\frac{ia}{p}\right \rangle\left(1-\left
    \langle\frac{ia}{p}\right \rangle\right)[a/p].
\end{align*}
Therefore,
\begin{align*}
  {\rm ord}_{a/p} \left (
  \prod_{i=1}^{(p-1)/2}K_{p,i}^{a_{i}} \right ) = -\sum_{i=1}^{(p-1)/2}\frac{p}{2}\left
  \langle\frac{ia}{p} \right \rangle\left(1-\left \langle\frac{ia}{p}
  \right \rangle\right)a_{i}.
\end{align*}
Note also that $\eta(p\tau)^{a_{0}+2\sum_{i=1}^{(p-1)/2}a_{i}}$ has order supported on cusps that are $\Gamma_{0}(p)$-equivalent to $0$. Therefore, the eta factor is holomorphic at each cusp and satisfies the requisite transformation formulas with respect to $\Gamma_{1}(p)$
if and only if $$24 \mid a_{0}+2\sum_{i=1}^{(p-1)/2}a_{i}\geq
0.$$
Thus,
\begin{align*}
  {\rm ord}_{a/p} \left (
  f_{a_{0},\ldots,a_{(p-1)/2}}\right ) = \frac{p}{24}\left (a_{0}+2\sum_{i=1}^{(p-1)/2}a_{i} \right )-\sum_{i=1}^{(p-1)/2}\frac{p}{2}\left
  \langle\frac{ia}{p} \right \rangle\left(1-\left \langle\frac{ia}{p}
  \right \rangle\right)a_{i}.
\end{align*}
Hence, the conditions in \eqref{ineq} assure
$f_{a_{0},\ldots,a_{(p-1)/2}}(\tau)$ is
holomorphic at each cusp.
\end{proof}

 For $p = 5$ and each fixed $a_{0}$, the system \eqref{ineq} subsumes a
 finite $2$-dimensional $\Bbb Z$-polyhedron
 \begin{align}\label{aa5}
 a_{0}+2a_{1}+2a_{2} \geq0,
    \quad  25 a_0+2 a_1-22 a_2\ge 0, \quad 
  25 a_0-22 a_1+2 a_2 \ge 0.
 \end{align}
The bounded characterization turns out to be more generally true and implies that the set of
Klein form quotients of the form \eqref{faaa} that are modular with
respect to $\Gamma_{1}(p)$ is enumerable. The inequalities of \eqref{ineq} indeed define a bounded polyhedron whose lattice points
may be systematically mined for modular forms of the form \eqref{eq:p} satisfying Ramanujan type congruences.

\begin{Corollary} \label{finitez}
For a given  positive even integer $a_{0}$, there are a finite number of
subscripts $(a_{1}, \ldots, a_{(p-1)/2})$ such that
 $ f_{a_{0},\ldots,a_{(p-1)/2}}(\tau)$
is a modular
form of weight~$\frac{a_{0}}{2}$ for the subgroup $\Gamma_{1}(p)$. 
\end{Corollary}

\begin{proof}
For a given $a_{0}\in 2\Bbb Z^{+}$, the system of inequalities
\eqref{ineq} represents a polyhedron of dimension~$(p-1)/2$ over
$\mathbb{R}$. It is not hard to show that  the system does not contain a ray by noting that the system 
\begin{align*}
\nonumber {2\sum_{i=1}^{(p-1)/2}a_{i}}\geq0,\qquad 
\nonumber 
  \frac{p}{12}\sum_{i=1}^{(p-1)/2}a_{i}-\sum_{i=1}^{(p-1)/2}\frac{p}{2}\left 
  \langle\frac{ia}{p}\right \rangle\left(1-\left \langle\frac{ia}{p}
  \right \rangle\right)a_{i}\geq0\quad\mbox{for $a=1,\ldots,(p-1)/2$}
\end{align*}
only admits the zero vector, and thus it is a bounded polyhedron \cite{Gr}. Therefore, it is a bounded convex polytope of full dimension and contains finitely many lattice points.
\end{proof}


The lattice points described by the system \eqref{ineq} may be parameterized.
\begin{Lemma} \label{lemp}
For primes $p \ge 5$ and each fixed positive even integer $a_0$, let $L(a_0)$ denote the set of solutions to the system \eqref{ineq}.
\begin{enumerate}
    \item For $p=5$ and $b, t \in \Bbb Z$,
\begin{align*}
     L(a_0) = \bigg \{ (a_0, a_1, a_2) \ \bigg | \begin{array}{l}   a_{1} = -\frac{3}{2} a_{0}+11 b+5t, \\ a_{2} = a_0 + b- 5t,\end{array} \begin{array}{l} 0 \le b \le \left \lfloor \frac{a_0}{4} \right \rfloor, \\  0 \le t \le \frac{a_{0} - 4b}{2}\end{array} \bigg \}.
\end{align*}
\item For $p=7$ and $b, t_1, t_2 \in \Bbb Z$,
\begin{align*}
     L(a_0) = \left \{ (a_0, a_1, a_2, a_3) \left |  \begin{array}{l}  a_{1}	=-\frac{a_{0}}{2}+5b+2t_{1}-t_{2},\\
 a_{2}	=-\frac{3a_{0}}{2}+8b+t_{1}+3t_{2},\\
\nonumber a_{3} =\frac{3a_{0}}{2}-b-3t_{1}-2t_{2},\end{array} \right. \begin{array}{l} 0	\le b\le\left\lfloor  \frac{a_{0}}{3}\right\rfloor , \\
\label{r7} 0	\le t_{1}\le a_{0}-3b, \\ 	0 \le t_{2}\le a_{0}-3b-t_{1}\end{array} \right \}.
\end{align*}
\item For $p=11$, $L(a_0)$ is parameterized by 
\begin{equation}\label{pa}
\vec{x} =\begin{pmatrix}
 a_1 \\ a_2 \\ a_3 \\ a_4 \\ a_5
\end{pmatrix}= 
\frac{1}{5}\left(
\begin{array}{cccccc}
 -16 & 27 & 7 & 4 & 1 & 3 \\
 19 & -8 & -3 & -6 & -4 & -7 \\
 -11 & 22 & 2 & -1 & 6 & 3 \\
 4 & 7 & -3 & -1 & -4 & 3 \\
 -1 & 12 & -3 & 4 & 1 & -2 \\
\end{array}
\right) \begin{pmatrix}
 a_0/2 \\ b \\ t_1 \\ t_2 \\ t_3 \\ t_4
\end{pmatrix},
\end{equation}
where 
\begin{align}
0 \le b \le a_0/2, \qquad 0 \le t_i \le \frac{5a_0}{2} - 5b - \sum_{j=1}^{i-1} t_i\quad \text{for }1\le i \le 4,
\end{align}
and 
\begin{align}
  a_0/2-12 b+3 t_1-4 t_2- t_3+2 t_4 &\equiv 0 \pmod{5}.
\end{align}
\item For each prime $p\ge 5$, $L(a_0)$ is parameterized, for appropriate bounds on $b, t_{i}$, by \begin{align*}
 \label{pa}
\begin{pmatrix}
 a_1 \\ a_2 \\ \vdots\\ a_{\frac{p-1}{2}}
\end{pmatrix}   =\left(\begin{array}{cccc}
1 & 1 & \cdots & 1\\
\lambda_{1,1} & \lambda_{1,2} & \cdots & \lambda_{1,\frac{p-1}{2}}\\
\lambda_{2,1} & \lambda_{2,2} & \cdots & \lambda_{2,\frac{p-1}{2}}\\
\vdots & \vdots & \vdots & \vdots\\
\lambda_{\frac{p-3}{2},1} & \lambda_{\frac{p-3}{2},2} & \cdots & \lambda_{\frac{p-3}{2},\frac{p-1}{2}}
\end{array}\right)^{-1}\left(\begin{array}{cccccc}
-1 & 12 & 0 & 0& \cdots & 0\\
-p^{2} & 0 & 12p & 0 & \cdots & 0\\
-p^{2} & 0 & 0 & 12p & \ddots & \vdots\\
\vdots & \vdots & \vdots & \ddots & \ddots & 0\\
-p^{2} & 0 & 0 & \cdots & 0 & 12p
\end{array}\right)\begin{pmatrix}
 a_0/2 \\ b \\ t_1 \\ t_2 \\ \vdots \\ t_{\frac{p-3}{2}}
\end{pmatrix},
\end{align*}
where, for $\alpha$, the minimal positive primitive root modulo $p$,
\begin{align}
    \lambda_{m,i}=p^{2}\left(1+6\left\langle \frac{\alpha^{m-1}i}{p}\right\rangle -6\left\langle \frac{\alpha^{m-1}i}{p}\right\rangle^{2} \right),
\end{align}
 and where a finite set of linear congruence conditions is assumed to ensure integrality. For $p=13, 17, 19$ the conditions required to ensure integrality are, respectively 
    $$6a_0+25 b-11 t_1+9 t_2+2 t_3-6 t_4-7 t_5 \equiv 0 \pmod{19},$$
   $$ \frac{1013a_0}{2}+152 b-253 t_1-198 t_2-83 t_3 -108 t_4-t_5-202 t_6+121 t_7 \equiv 0 \pmod{584},$$
   $$ \frac{1823 a_0}{2}+ 4458 b-1145 t_1+1097 t_2+489 t_3 -89 t_4 -898 t_5+510 t_6-689 t_7-661 t_8 \equiv 0 \pmod{4383}.$$

\end{enumerate}
\end{Lemma}
\begin{proof}
We begin with a proof for the case $p=5$. To satisfy the congruences required by \eqref{ineq}, we make the substitution
\begin{align*}
    a_{0}+2a_{1}+2a_{2}	=24\beta,\qquad 25a_{0}+2a_{1}-22a_{2}	=5r,\qquad b,r \in \Bbb N\cup \{0\}.
\end{align*}
Thus, it follows that 
\begin{align*}
  a_{1} = -\frac{3}{2} a_{0}+11 b+\frac{5 r}{24}, \qquad a_{2} = a_0 + b - \frac{5 r}{24}.
\end{align*}
The ordered pair $(a_1,a_2)$ is integral if and only if $5r/24\in \Bbb Z$. Replacing $r$ by $24t$ gives the parametrization
\begin{align*}
  a_{1} = -\frac{3}{2} a_{0}+11 b + 5t, \qquad a_{2} = a_0 + b - 5t.
\end{align*}
Substituting the formulas for $a_{1}$ and $a_{2}$ into \eqref{aa5} results in
\begin{align*}
  b &\ge 0, \qquad t \ge 0, \qquad \frac {a_{0}}{2}-2b -t \ge 0.
\end{align*}
This proves the claimed result. A similar strategy allows us to derive the level $7$ parameterization. To generalize the parameterizations to higher levels, standard tools may be applied to solve linear Diophantine systems. We demonstrate this in detail for the case $p=11$. Here \eqref{ineq} becomes
\begin{align} \label{111}
    0 \le a_{0}+2a_{1}+2a_{2}+2a_{3}+2a_{4}+2a_{5} &\equiv 0 \pmod {24}, \\ \label{112} 0 \le  121a_{0}+122a_{1}+26a_{2}-46a_{3}-94a_{4}-118a_{5} &\equiv 0 \pmod {11}, \\ \label{113} 0 \le 121a_{0}+26a_{1}-94a_{2}-118a_{3}-46a_{4}+122a_{5} &\equiv 0 \pmod {11}, \\ 0 \le 121a_{0}-94a-46a_{2}+122a_{3}-118a_{4}+26a_{5} &\equiv 0 \pmod {11}, \\0 \le 121a_{0}-46a_{1}-118a_{2}+26a_{3}+122a_{4}-94a_{5} &\equiv 0 \pmod {11}, \\ \label{116} 0 \le 121a_{0}-118a_{1}+122a_{2}-94a_{3}+26a_{4}-46a_{5} &\equiv 0 \pmod {11}.
\end{align}
We parameterize the first $5$ congruences via
\begin{align} \label{11mat}
 \left(
\begin{array}{ccccc}
 1 & 1 & 1 & 1 & 1 \\
 61 & 13 & -23 & -47 & -59 \\
 13 & -47 & -59 & -23 & 61 \\
 -47 & -23 & 61 & -59 & 13 \\
 -23 & -59 & 13 & 61 & -47 \\
\end{array}
\right) \begin{pmatrix}
 a_1 \\ a_2 \\ a_3 \\ a_4 \\ a_5
\end{pmatrix}  = \frac{1}{2} \left(
\begin{array}{c}
 24 b- a_0 \\
  11\cdot 24 t_1-121 a_0 \\
 11\cdot 24  t_2-121 a_0 \\
 11\cdot 24  t_3-121 a_0 \\
 11\cdot 24  t_4-121 a_0 \\
\end{array}
\right).
\end{align}
If we denote the system in \eqref{11mat} as $A\vec{x} = \vec{b}$, then the Smith normal form of $A$ is $D = LAR$, where $D$ is the diagonal matrix with entries $( 1, 12, 132, 132, 660)$, and
\begin{align*}
  L=\left(
\begin{array}{ccccc}
 1 & 0 & 0 & 0 & 0 \\
 95 & -1 & 1 & 1 & 0 \\
 154 & -2 & 3 & 2 & -1 \\
 121 & -3 & 3 & 0 & -1 \\
 132 & -3 & 4 & 1 & -2 \\
\end{array}, 
\right), \quad 
R = \left(
\begin{array}{ccccc}
 1 & -1 & 9 & 4 & -34 \\
 0 & 1 & -10 & -5 & 41 \\
 0 & 0 & 1 & 0 & -4 \\
 0 & 0 & 0 & 1 & -4 \\
 0 & 0 & 0 & 0 & 1 \\
\end{array}
\right). 
\end{align*}
To determine all integer solutions to the linear system $A\vec{x} = \vec{b}$ in \eqref{11mat}, we note that the system is equivalent to $D\vec{y} = \vec{c}$, where $\vec{y} = R^{-1}\vec{x}$ and 
\begin{align}
    \vec{c} =L\vec{b} = \left ( \begin{array}{c}  12 b-a_0/2 \\ -12 (9 a_0-95 b+11 t_1-11 t_2-11 t_3) \\ -66 (3 a_0-28 b+4 t_1-6 t_2-4 t_3+2 t_4) \\ 132 (11 b-3 t_1+3 t_2-t_4) \\ -66 (a_0-24 b+6 t_1-8 t_2-2 t_3+4 t_4) \end{array} \right).
\end{align}
The system $D\vec{y} = \vec{c}$ has a solution over $\Bbb Z^{5}$ if and only if the diagonal entries of $D$ divide the respective entries of $\vec{c}$. This is clearly satisfied except for divisibility of the last entries 
\begin{align} \label{go}
    660 \mid -66 (a_0-24 b+6 t_1-8 t_2-2 t_3+4t_4) .
\end{align}
The divisibility requirement in \eqref{go} is satisfied if and only if 
\begin{align}
    a_0/2-12 b+3 t_1-4 t_2- t_3+2 t_4\equiv 0 \pmod{5}.
\end{align}
Therefore, \eqref{go} holds precisely when the parameterization $(a_{1}, a_{2}, \ldots, a_{5})^{T} = R \vec{y}$ is integral. This is equivalent to the claim in the lemma for the case $p=11$. The bounds for the parameters $b, t_i$, $1 \le i \le 4$, are obtained by substituting \eqref{11mat} into the inequalities \eqref{111}--\eqref{116}.
\end{proof}

The polytope parameterizations permit an enumeration of the lattice points satisfying \eqref{ineq}. 
\begin{Corollary}
Let $L(a_0)$ be the number of lattice points in the system determined by \eqref{ineq}. 
\begin{enumerate}
    \item For $p=5$, $|L(a_0)| = \left\lfloor \frac{\left(a_{0}+4\right)^{2}}{16}\right\rfloor.$
    \item For $p=7$, $|L(a_0)| = \left\lfloor \frac{(a_{0}+2)^{3}+3(a_{0}+2)^{2}}{18}\right\rfloor.$
\end{enumerate} 
\end{Corollary}
\begin{proof}
A proof for the case $p=5$ is given that demonstrates the ideas needed for the $p=7$ case. From the parameterization for the lattice points in the polytope from Lemma \ref{lemp}, we have 
\begin{align}
  \nonumber   L(a_{0})	&=\sum_{b=0}^{\lfloor\frac{a_{0}}{4}\rfloor}\frac{a_{0}-4b+2}{2} \\ &=
\nonumber 	\frac{a_{0}+2}{2}\sum_{b=0}^{\left\lfloor\frac{a_{0}}{4}\right\rfloor}1-2\sum_{b=0}^{\left\lfloor\frac{a_{0}}{4}\right\rfloor}b \\ \nonumber 
	&=\left(\left\lfloor\frac{a_{0}}{4}\right\rfloor+1\right)\left(\frac{a_{0}+2}{2}-\left\lfloor\frac{a_{0}}{4}\right\rfloor\right) \\ \label{lasteq}
	&=\left\lfloor \frac{\left(a_{0}+4\right)^{2}}{16}\right\rfloor .
\end{align}
 Equality \eqref{lasteq} can be proven in cases. 
For $a_{0}=4k+2, k\in\mathbb{Z}$, note that
\begin{align*}
\left\lfloor \left(\frac{a_{0}+4}{4}\right)^{2}\right\rfloor &=\left\lfloor \left(\frac{4k+2+4}{4}\right)^{2}\right\rfloor =\left\lfloor \left(k+1+\frac{1}{2}\right)^{2}\right\rfloor \\ 
&=\left\lfloor \left(k+1\right)^{2}+(k+1)+\frac{1}{4}\right\rfloor \\ &=k^{2}+3k+2 \\ &=\left(\left\lfloor\frac{a_{0}}{4}\right\rfloor+1\right)\left(\frac{a_{0}+2}{2}-\left\lfloor\frac{a_{0}}{4}\right\rfloor\right).
\end{align*}
Equality \eqref{lasteq} may be proven for the case $a_{0} \equiv 0 \pmod{4}$ similarly. 
\end{proof}

\section{Construction of relevant vector spaces in terms of Klein forms}
Products of Klein forms generate the graded algebras of modular forms relevant to this work. We next establish explicit bases in terms of products of the form \eqref{eq:p} for the vector spaces
$M_{1}(\Gamma_{1}(p))$. These serve to generate the associated graded ring of modular forms. For primes $5 \le p \le 19$, the space of weight $1$ cusp forms is trivial (see, e.g., \cite{DS}). Therefore, a natural way to construct a basis is in terms of Eisenstein series. This construction along with the next lemma yields a parallel decomposition in terms of quotients of Klein forms
\begin{align} \label{decompg1}
    M_{1}(\Gamma_{1}(p)) = \bigoplus_{n=0}^{(p-1)/2} \Bbb C f_{\sigma_{p}^{n}(2,a_{1}, \ldots, a_{(p-1)/2})},
\end{align}
for the permutation $\sigma_{p}$ defined by~\eqref{sigmap}. 
Throughout this work, for a fixed $\alpha$ such that $(\mathbb{Z}/p\mathbb{Z})^{\times}/\{\pm1\}$, write $\gamma_{p}$ for some chosen $\left (\begin{smallmatrix}\tilde{\alpha}&p\\p\beta&\delta\end{smallmatrix} \right )\in\Gamma_{0}(p)$ with $\tilde{\alpha}$ the inverse of $\alpha\pmod{p}$.
\begin{Lemma} \label{lema}
Suppose that $(\mathbb{Z}/p\mathbb{Z})^{\times}/\{\pm1\}=\langle\alpha\rangle$, and suppose that $a_{1},\ldots,a_{(p-1)/2}\in\mathbb{Z}$ such that \eqref{ineq} is satisfied and 
 $$
 \frac{p}{24}\left (2+2\sum_{i=1}^{(p-1)/2}a_{i} \right )-\sum_{i=1}^{(p-1)/2}\frac{p}{2}\left
  \langle\frac{i\alpha^{m}}{p} \right \rangle\left(1-\left \langle\frac{i\alpha^{m}}{p}
  \right \rangle\right)a_{i}
$$
has a unique minimum as $m$ ranges over $0,\ldots,(p-3)/2$.  Then for any odd prime $p\leq 19$, 
$$
\{f_{\sigma_{p}^{i}(2,a_{1},\ldots,a_{(p-1)/2})}|\,i=1,\ldots,(p-1)/2\}
$$
is a basis for $M_{1}(\Gamma_{1}(p))$.

\end{Lemma}

\begin{proof}
Clearly, by Lemma~\ref{m1r}, $f_{2,a_{1},\ldots,a_{(p-1)/2}}$ lies in $M_{1}(\Gamma_{1}(p))$. By Lemma~\ref{klf} (2), one can deduce that
$$
f_{\sigma_{p}(2,a_{1},\ldots,a_{(p-1)/2})}=\pm \left.f_{2,a_{1},\ldots,a_{(p-1)/2}}\right|_{1}\gamma_{p},
$$
and also lies in $M_{1}(\Gamma_{1}(p))$ since $\Gamma_{1}(p)$ is a normal subgroup of $\Gamma_{0}(p)$. Also, by Lemma~\ref{klf} (3), one can show that
\begin{align*}
&{\rm Div}_{X_{1}(p)}\left(f_{\sigma_{p}^{j}(2,a_{1},\ldots,a_{(p-1)/2})}\right)\\
&=\sum_{m=0}^{\frac{p-3}{2}}\left(\frac{p}{24}\left (2+2\sum_{i=1}^{(p-1)/2}a_{i} \right )-\sum_{i=1}^{(p-1)/2}\frac{p}{2}\left
  \langle\frac{i\alpha^{m+j}}{p} \right \rangle\left(1-\left \langle\frac{i\alpha^{m+j}}{p}
  \right \rangle\right)a_{i}\right)\left[\frac{\alpha^{m}}{p}\right].
  \end{align*}
  Now without loss of generality, assume that 
\begin{align*}
&\frac{p}{24}\left (2+2\sum_{i=1}^{(p-1)/2}a_{i} \right )-\sum_{i=1}^{(p-1)/2}\frac{p}{2}\left
  \langle\frac{i}{p} \right \rangle\left(1-\left \langle\frac{i}{p}
  \right \rangle\right)a_{i}\\
  &<\frac{p}{24}\left (2+2\sum_{i=1}^{(p-1)/2}a_{i} \right )-\sum_{i=1}^{(p-1)/2}\frac{p}{2}\left
  \langle\frac{i\alpha^{m}}{p} \right \rangle\left(1-\left \langle\frac{i\alpha^{m}}{p}
  \right \rangle\right)a_{i}
  \end{align*}
  for all $m\geq1$. 
Therefore, the leading coefficient of 
$$
c_{0}f_{2,a_{1},\ldots,a_{(p-1)/2}}+\cdots+c_{(p-3)/2}f_{\sigma_{p}^{(p-3)/{2}}(2,a_{1},\ldots,a_{(p-1)/2})}
$$
at the cusp $\left[\frac{\alpha^{-j}}{p}\right]$ is exactly $\pm c_{j}$, and thus
$$
c_{0}f_{2,a_{1},\ldots,a_{(p-1)/2}}+\cdots+c_{(p-3)/2}f_{\sigma_{p}^{({p-3})/{2}}(2,a_{1},\ldots,a_{(p-1)/2})}=0
$$
yields that $c_{j}=0$ for all $j=0,\ldots,(p-3)/2$, that is, 
$$
\{f_{\sigma_{p}^{i}(2,a_{1},\ldots,a_{(p-1)/2})}|\,i=1,\ldots,(p-1)/2\}
$$
 is a linearly independent set. Finally, for an odd prime $p\leq 19$, it is known that ${\rm dim}(M_{1}(\Gamma_{1}(p))=\frac{p-1}{2}$, and this together with the conclusion above justifies the lemma.
\end{proof}

A set of subscripts $a_{1}, \ldots, a_{(p-1)/2}$ for \eqref{decompg1} saisfying the requirements of Lemma \ref{lema} may be canonically determined for each prime $5\le p \le 19$ by the fact that the normalized sum of the weight one Hecke Eisenstein series is a product of the form \eqref{eq:p}. For odd
Dirichlet character $\chi$ modulo $p$, the Hecke Eisenstein series are defined by
\begin{align*}
   E_{\chi,k}(\tau) = 1 +  \frac{2}{L(1 - k, \chi)} \sum_{n=1}^{\infty} \chi(n) \frac{n^{k-1} q^{n}}{1 - q^{n}},
\end{align*}
where $L(1  - k, \chi)$ is the analytic continuation of the associated
Dirichlet $L$-series.
\begin{Theorem}{\cite[Theorem 3.2]{HLM}}
  For primes $p$ with $5 \le p \le 19$ there exist $\vec{a}_{p}=(a_{0},\ldots,a_{(p-1)/2}) \in \{0,\pm 1,\pm 2, \pm 3\}^{(p+1)/2}$ such that 
\begin{align} \label{eiss}
  \sum_{\chi(-1) = -1} E_{\chi,1}(\tau) =
 \frac{p-1}{2} f_{\vec{a}_{p}}.
\end{align}
\end{Theorem} For instance, if the odd Dirichlet characters modulo $p$ are defined in terms of their respective primitive roots by $\chi_{n,5}(2) =(-1)^{n-1}i$, $\langle \chi_{n,7}(3) \rangle_{n=1}^{3} = \langle e^{\pi i/3}, -1,  e^{-\pi i/3}\rangle$, $\langle \chi_{n,11}(2) \rangle_{n=1}^{5} = \langle e^{\pi i/5}, e^{3\pi i/5}, -1, e^{-3\pi i/5}$, $e^{-\pi i/5} \rangle$, then \begin{align*}
&E_{\chi_{1,5},1}(\tau) + E_{\chi_{2,5},1}(\tau) =
  2\frac{K_{5,2}^{2}}{K_{5,1}^{3}}, \qquad E_{\chi_{1,7},1}(\tau) + E_{\chi_{2,7},1}(\tau) + E_{\chi_{3,7},1}(\tau) =
  3\frac{K_{7,3}}{K_{7,1}^{2}}, \\ 
  &E_{\chi_{1,11},1}(\tau) + E_{\chi_{2,11},1}(\tau) + E_{\chi_{3,11},1}(\tau)  + E_{\chi_{4,11},1}(\tau)  + E_{\chi_{5,11},1}(\tau) =
  5\frac{K_{11,4}}{K_{11,1}K_{11,2}}
\end{align*}
Similar formulas hold for $p=13, 17, 19$.
The next corollary specializes~\eqref{decompg1} to the cases $5\leq p\leq 19$ by generating the vector spaces $M_{1}(\Gamma_{1}(p))$ via the Klein form representations for Eisenstein sums. Permutations of exponents of the products generate the graded algebra of modular forms for $\Gamma_{1}(p)$.
\begin{Lemma}[Rustom \cite{Ru}] \label{rus}
 For $N\geq5$, the graded $\mathbb{C}$-algebra $\bigoplus_{k=2}^{\infty}M_{k}(\Gamma_{1}(N))$ is generated in weight at most~$3$.
\end{Lemma}

\begin{Corollary} \label{tf}
For each prime $5 \le p \le 19$, let $\sigma_{p}$ be defined as in \eqref{sigmap} and let $$\vec{a}_{5}=(2,-3,2),\ \vec{a}_{7}=(2,-2,0,1),\ \vec{a}_{11}=(2,-1,-1,0,1,0),$$  $$\vec{a}_{13}=(2,-1, 0, -1, 0, 0, 1),\ \vec{a}_{17}=(2,0,-1,-1,0,0,0,0,1),\ \vec{a}_{19}=(2,0, 0, -1, -1, -1, 0, 0, 1, 1).$$
Then\begin{align}
   \label{d1} M_{1}(\Gamma_{1}(p))&=\bigoplus_{n=0}^{\frac{p-3}{2}}\mathbb{C}f_{\sigma_{p}^{n}(\vec{a}_{p})}, \\ 
   \label{d2} \bigoplus_{k=1}^{\infty}M_{k}(\Gamma_{1}(p)) &=\left  \langle f_{\sigma_{p}^{n}(\vec{a}_{p})}\right \rangle_{n=0,\mathbb{C}}^{\frac{p-3}{2}},
\end{align}
as graded rings.
\end{Corollary}

\begin{proof}
The direct sum decompositions in \eqref{d1} follow from Lemmas~\ref{klf} and~\ref{lema}. To show that the weight one forms from \eqref{d1} generate the graded algebra of integer weight modular forms for $\Gamma_{1}(p)$ we establish that for primes $5\leq p\leq 19$, the weight one forms $f_{\sigma_{p}^{n}(\vec{a}_{p})}$ generate bases for $M_{k}(\Gamma_{1}(p))$ for $k=2$ and~3, and thus, by Lemma \ref{rus} and Corollary \ref{tf}, the products generate the graded $\mathbb{C}$-algebra of modular forms for $\Gamma_{1}(p)$. For $p=5$, let $n=0,1$, and define $x_{n} = f_{\sigma_{5}^{n}(2,-3,2)}$. Then, for $1\leq k\leq 3$, the set of $k+1$ monomials in $x_{0}$ and $x_{1}$ of degree $k$ form an ordered basis for $M_{k}(\Gamma_{1}(5))$. Similarly, if $x_{n}$ for $0\leq n\leq 2$ is defined by $x_{n} = f_{\sigma_{7}^{n}(2,-2,0,1)}$, for $1\leq k\leq 3$, 
\begin{align} \label{7b}
 \{x_{0}^{k},x_{1}^{k}&,x_{2}^{k}\}\sqcup\bigcup_{j=1}^{\lfloor\frac{k}{2}\rfloor}\{(x_{0}x_{2})^{j}x_{0}^{k-2j},(x_{1}x_{2})^{j}x_{1}^{k-2j}\}\sqcup\bigcup_{j=1}^{\lfloor\frac{k-1}{2}\rfloor}\{(x_{0}x_{2})^{j}x_{2}^{k-2j},(x_{1}x_{2})^{j}x_{2}^{k-2j}\}    
\end{align}
is an ordered basis for $M_{k}(\Gamma_{1}(7))$. Likewise, if $x_{n}=f_{\sigma_{11}^{n}(2,-1, -1,0,1,0)}$ for $0\leq n\leq 4$, then 
\begin{align}
\{x_{0}^{k},\ldots,x_{4}^{k}\}&\sqcup\bigcup_{j=1}^{\lfloor\frac{k}{2}\rfloor}\{(x_{0}x_{3})^{j}x_{0}^{k-2j},(x_{2}x_{4})^{j}x_{2}^{k-2j},(x_{3}x_{1})^{j}x_{3}^{k-2j},(x_{1}x_{2})^{j}x_{1}^{k-2j},(x_{4}x_{0})^{j}x_{4}^{k-2j}\} \nonumber \\
&\sqcup\bigcup_{j=1}^{\lfloor\frac{k-1}{2}\rfloor}\{(x_{0}x_{2})^{j}x_{1}^{k-2j},(x_{2}x_{5})^{j}x_{3}^{k-2j},(x_{3}x_{1})^{j}x_{2}^{k-2j},(x_{1}x_{2})^{j}x_{4}^{k-2j},(x_{4}x_{1})^{j}x_{0}^{k-2j}\}
\end{align}
is an ordered basis for $M_{k}(\Gamma_{1}(11))$ for $1 \le k \le 5$. For level $13$, if $x_{n} =f_{\sigma_{13}^{n}(2,-1, 0, -1, 0, 0, 1)}$, $0 \le n \le 5$, then a basis for $M_{2}(\Gamma_{1}(13))$ is \begin{align*}
  \{x_0^2,x_0 x_1,x_1^2,x_0 x_2,x_1 x_2,x_2^2,x_0 x_3,x_1 x_3,x_3^2,x_0 x_4,x_4^2,x_1 x_5,x_5^2\},
\end{align*}
and a basis for $M_{3}(\Gamma_{1}(13))$ is
\begin{align*}
         \left  \{ \begin{array}{c}  
    x_0^3,\,x_0^2 x_1,\,x_0 x_1^2,\,x_1^3,\,x_0^2 x_2,\,x_0 x_1 x_2,\,x_1^2 x_2,\,x_0 x_2^2,\,x_1 x_2^2,\,x_2^3, \\ x_0^2 x_3,\,x_1^2 x_3,\,x_0 x_3^2,\,x_3^3,\,x_0^2 x_4,\,x_0 x_4^2,\,x_4^3,\,x_1^2 x_5,\,x_1 x_5^2,\,x_5^3 \end{array} \right \}.
\end{align*}
Likewise if $x_{n} =f_{\sigma_{17}^{n}(2,0, -1, -1, 0, 0, 0, 0, 1)}$, $0 \le n \le 7$, then a basis for $M_{2}(\Gamma_{1}(17))$ is
\begin{align*}
  \{ x_0^2,\,x_0 x_1,\,x_1^2,\,x_0 x_2,\,x_1 x_2,\,x_2^2,\,x_0 x_3,\,x_1 x_3,\,x_2 x_3,\,x_3^2,\,x_0 x_4,\,x_1 x_4,\,x_3 x_4,\,x_4^2,\,x_0 x_5,\,x_5^2,\,x_0 x_6,\,x_6^2,\,x_1 x_7,\,x_7^2\},
\end{align*}
and a basis for $M_{3}(\Gamma_{1}(17))$ is
\begin{align*}
 \left  \{ \begin{array}{c}  
     x_0^3,\,x_0^2 x_1,\,x_0 x_1^2,\,x_1^3,\,x_0^2 x_2,\,x_0 x_1 x_2,\,x_1^2 x_2,\,x_0 x_2^2,\,x_1 x_2^2,\,x_2^3,\,x_0^2 x_3,\,x_0 x_1 x_3,\,x_1^2 x_3,\,x_0 x_2 x_3,\,x_1 x_2 x_3, \\ x_2^2 x_3,\,x_0 x_3^2,\,x_1 x_3^2,\,x_3^3,\,x_0^2 x_4,\,x_1^2 x_4,\,x_1 x_4^2,\,x_4^3,\,x_0^2 x_5,\,x_0 x_5^2,\,x_5^3,\,x_0^2 x_6,\,x_0 x_6^2,\,x_6^3,\,x_1^2 x_7,\,x_1 x_7^2,\,x_7^3
    \end{array} \right \}
\end{align*}
Finally, if $x_{n} = f_{\sigma_{19}^{n}(2,0, 0, -1, -1, -1, 0, 0, 1, 1)}$, $0\le n \le 8$, then a basis for $M_{2}(\Gamma_{1}(19))$ is 
\begin{align*}
  \left  \{ \begin{array}{c}   
x_0^2,\,x_0 x_1,\,x_1^2,\,x_0 x_2,\,x_1 x_2,\,x_2^2,\,x_0 x_3,\,x_1 x_3,\,x_2 x_3,\,x_3^2,\,x_0 x_4,\,x_1 x_4, \\  x_2 x_4,\,x_3 x_4,\,x_4^2,\,x_0 x_5,\,x_1 x_5,\,x_5^2,\,x_0 x_6,\,x_6^2,\,x_0 x_7,\,x_7^2,\,x_1 x_8,\,x_8^2 
 \end{array} \right \},
\end{align*}  
and a basis for  $M_{3}(\Gamma_{1}(19))$ is 
\begin{align*}
 \left  \{ \begin{array}{c} 
 x_0^3,\,x_0^2 x_1,\,x_0 x_1^2,\,x_1^3,\,x_0^2 x_2,\,x_0 x_1 x_2,\,x_1^2 x_2,\,x_0 x_2^2,\,x_1 x_2^2,\,x_2^3,\,x_0^2 x_3,\,x_0 x_1 x_3,\,  x_1^2 x_3,\\  x_0 x_2 x_3,\,x_1 x_2 x_3, \,x_2^2 x_3,\,x_0 x_3^2,\,x_1 x_3^2,\,x_2 x_3^2,\,x_3^3,\,x_0^2 x_4,\,x_0 x_1 x_4,\,x_1^2 x_4,\,x_0 x_2 x_4,\,x_2^2 x_4, \\  x_0 x_4^2,\, x_4^3,\,x_1^2 x_5,\,_1 x_5^2,\,x_5^3,\,x_0^2 x_6,\, x_0 x_6^2,\,x_6^3,\,x_0^2 x_7,\,x_0 x_7^2,\,x_7^3,\,x_1^2 x_8,\,x_1 x_8^2,\,x_8^3
 \end{array} \right \}.
\end{align*}
\end{proof}

In the proof of Lemma~\ref{lema}, we note that $f_{2,\ldots,a_{(p-1)/2}}|_{1}\gamma_{p}$ and $f_{\sigma_{p}(2,a_{1},\ldots,a_{(p-1)/2})}$ differ by a sign. The sign is explicitly determined in the following lemma.
\begin{Lemma}
\label{flf}
  For a product $f_{a_{0},\ldots,a_{(p-1)/2}}(\tau)$, and $\gamma_{p}=\begin{pmatrix}\tilde{\alpha}&p\\p\beta&\delta\end{pmatrix}\in\Gamma_{0}(p)$, supposing $\tilde{\alpha}i=q_{i}p+r_{i}$ for $-\frac{p-1}{2}\leq r_{i}\leq \frac{p-1}{2}$, one has that
  \begin{enumerate}
      \item 
      $$
 \left.f_{a_{0},\ldots,a_{(p-1)/2}}(\tau/p)\right|_{a_{0}/2}\gamma_{p}=\lambda(\gamma_{p},\vec{a}_{p})f_{\sigma_{p}(a_{0},a_{1},\ldots,a_{(p-1)/2})}(\tau/p),
 $$
 where  for $\vec{a}_{p}=(a_{0},\ldots,a_{(p-1)/2})$,
 $$
 \lambda(\gamma_{p},\vec{a}_{p})=\prod_{i=1}^{\frac{p-1}{2}}{\rm sgn}(r_{i})^{a_{i}}(-1)^{A(\gamma,\vec{a}_{p})}
 $$
 with
 $$
 A(\gamma,\vec{a}_{p})=\sum_{i=1}^{\frac{p-1}{2}}\left(q_{i}i+i+q_{i}+\frac{r_{i}i}{p}\right)a_{i},
 $$
 \item 
      $$
 \left.f_{a_{0},\ldots,a_{(p-1)/2}}(\tau)\right|_{a_{0}/2}\gamma_{p}=\lambda(\gamma_{p},\vec{a}_{p})^{p}f_{\sigma_{p}(a_{0},\ldots,a_{(p-1)/2})}(\tau),
 $$
 \item  when $f_{a_0,\ldots,a_{(p-1)/2}}(\tau)\in M_{k}(\Gamma_{1}(p))$,  $\lambda(\gamma_{p},\vec{a}_{p})=\pm1$.
  \end{enumerate}
\end{Lemma}

\begin{proof}
By Lemma~\ref{klf} one has that
 \begin{align*}
     \left.K_{\left(\frac{i}{p},0\right)}(\tau)\right|_{a_{0}/2}\gamma_{p}&=K_{\left(\frac{\tilde{\alpha}i}{p},i\right)}(\tau)\\
     &=(-1)^{q_{i}i+q_{i}+i+\frac{i^{2}}{p}}K_{\left(\frac{r_{i}}{p},0\right)}(\tau)\\
     &={\rm sgn}(r_{i})(-1)^{q_{i}i+q_{i}+i+\frac{r_{i}i}{p}}K_{\left(\frac{r_{i}}{{\rm sgn}(r_{i})p},0\right)}(\tau),
 \end{align*}
 and thus,
 $$
 \left.f_{a_{0},\ldots,a_{(p-1)/2}}(\tau/p)\right|_{a_{0}/2}\gamma_{p}=\lambda(\gamma_{p},\vec{a}_{p})f_{\sigma_{p}(a_{0},a_{1},\ldots,a_{(p-1)/2}}(\tau/p),
 $$
 where 
 $$
 \lambda(\gamma_{p},\vec{a}_{p})=\prod_{i=1}^{\frac{p-1}{2}}{\rm sgn}(r_{i})^{a_{i}}(-1)^{A(\gamma,\vec{a}_{p})},
 $$
 with
 $$
 A(\gamma,\vec{a}_{p})=\sum_{i=1}^{\frac{p-1}{2}}\left(q_{i}i+i+q_{i}+\frac{r_{i}i}{p}\right)a_{i}.
 $$
 
Proof of the second assertion is similar, and we leave the details to the reader. For the last assertion, note that $$\sum_{i=1}^{(p-1)/2}\frac{r_{i}ia_{i}}{p}\equiv  \tilde{\alpha}\sum_{i=1}^{(p-1)/2}\frac{i^{2}a_{i}}{p} \pmod{p},$$ and by the assumption that $f_{a_{0},\ldots,a_{(p-1)/2}}(\tau)$ is a holomorphic modular for $\Gamma_{1}(p)$, Lemma~\ref{m1r} indicates that
$$
-\sum_{i=1}^{(p-1)/2}\frac{p}{2}\left
  \langle\frac{i}{p}\right \rangle\left(1-\left \langle\frac{i}{p}
  \right \rangle\right)a_{i}\in\mathbb{Z},
  $$
  and this implies
 $\sum_{i=1}^{(p-1)/2}\frac{i^{2}a_{i}}{p}\in\mathbb{Z}$. Therefore,
 $$
 \sum_{i=1}^{\frac{p-1}{2}}\frac{r_{i}ia_{i}}{p}\in\mathbb{Z},
 $$
and 
the assertion follows.
\end{proof}


In Lemma \ref{fgh}, we trace the path of dissection components of products $f_{\vec{a}_{p}}(\tau/p)$ from \eqref{eq:p} under permutations of the subscripts $\vec{a}_{p}$. 
\begin{Lemma} 
\label{fgh}
Suppose that $\alpha$ is the least positive integer such that $(\mathbb{Z}/p\mathbb{Z})^{\times}/\{\pm1\}=\langle \alpha\rangle$. Let $\sigma_{p}$ by defined by~\eqref{sigmap}. 
Denote by $\tilde{\alpha} \equiv \alpha^{-1} \pmod{p}$. For
$f_{a_{0}, \ldots, a_{(p-1)/2}} \in M_{a_{0}/2}(\Gamma_{1}(p))$,
suppose
$$
f_{a_{0},a_{1},\ldots,a_{(p-1)/2}}(\tau/p)=\sum_{r=0}^{p-1}g_{r}(\tau), \quad g_{r}(\tau)\in q^{\frac{r}{p}}\mathbb{Z}[[q]].
$$
Then
$$
f_{\sigma_{p}(a_{0},a_{1},\ldots,a_{(p-1)/2})}(\tau/p)=\sum_{r=0}^{p-1}h_{s(r)}(\tau), \quad h_{s(r)}(\tau)\in q^{\frac{s(r)}{p}}\mathbb{Z}[[q]],
$$
with $s(r)\equiv\tilde{\alpha}^{2}r\pmod{p}$, and
$$
h_{s(r)}(\tau)=\lambda(\gamma_{p},\vec{a}_{p})^{-1}\left.g_{r}\right|_{a_{0}/2}\gamma_{p},\quad i.e.,\quad U_{p,s(r)}(f_{\sigma_{p}(\vec{a}_{p})})=\lambda(\gamma_{p},\vec{a}_{p})^{-1}U_{p,r}\left.(f_{\vec{a}_{p}})\right|_{a_{0}/2}\gamma_{p},
$$
where $\vec{a}_{p}=(a_{0},\ldots,a_{(p-1)/2})$, and
$
 \lambda(\gamma_{p},\vec{a}_{p})
 $ is defined as in Lemma~\ref{flf}.
In particular, when $r=0$, one has that
$$
U_{p}\left(f_{\sigma_{p}(\vec{a}_{p})}\right)=\lambda(\gamma_{p},\vec{a}_{p})^{-1}U_{p}\left(f_{\vec{a}_{p}}\right)|_{a_{0}/2}\gamma_{p},
$$
where $U_{p}(f):=U_{p,0}$ is defined by \eqref{defu}.

\end{Lemma}

\begin{proof}
  The preceding lemma implies that
$$
f_{\sigma_{p}(a_{0},a_{1},\ldots,a_{(p-1)/2}}(\tau/p)=\lambda(\gamma_{p},\vec{a}_{p})^{-1}\sum_{r=0}^{p-1}g_{r}|_{a_{0}/2}\gamma_{p}.
 $$
It remains to show that $g_{r}|_{\frac{a_0}{2}}\gamma_{p}$   is the $\tilde{\alpha}^{2}r$-component of the $p$-dissection of $f_{\sigma_{p}(a_{0},a_{1},\ldots,a_{(p-1)/2}}(\tau/p)$. First note that
$$
\left.g_{r}\right|_{a_{0}/2}\begin{pmatrix}1&1\\0&1\end{pmatrix}=\zeta_{p}^{r}g_{r}.
$$
  It is not hard to see that $g_{r}(\tau)$ is a linear combination of $f_{a_{0},\ldots,a_{(p-1)/2}}\left(\frac{\tau+\alpha^{m}}{p}\right)$ for $0 \le m \le p-1$. Since $f_{a_{0},\ldots,a_{(p-1)/2}}(\tau) \in M_{a_{0}/2}(\Gamma_{1}(p))$, and $\Gamma(p)$ is a normal subgroup of $\Gamma_{1}(p)$, then $g_{r}(\tau)$ is a modular form of weight~$a_{0}/2$ for $\Gamma(p)$. Note that
$$
\gamma_{p}\begin{pmatrix}1&1\\0&1\end{pmatrix}
=\kappa\begin{pmatrix}1&\tilde{\alpha}^{2}\\0&1\end{pmatrix}\gamma_{p},
$$
where $\kappa=\gamma_{p}\begin{pmatrix}1&1\\0&1\end{pmatrix}\gamma_{p}^{-1}\begin{pmatrix}1&-\tilde{\alpha}^{2}\\0&1\end{pmatrix}\in\Gamma(p)$ since $\tilde{\alpha}\delta\equiv1\pmod{p}$ by assumption. Then
\begin{align*}
&\left.\left(\left.g_{r}\right|_{a_{0}/2}\gamma_{p}\right)\right|_{a_{0}/2}\begin{pmatrix}1&1\\0&1\end{pmatrix}\\
&=\left.g_{r}\right|_{a_{0}/2}\kappa \begin{pmatrix}1&\tilde{\alpha}^{2}\\0&1\end{pmatrix}\gamma_{p}\\
&=\left.g_{r}\right|_{a_{0}/2}\begin{pmatrix}1&\tilde{\alpha}^{2}\\0&1\end{pmatrix}\gamma_{p}\\
&=\zeta_{p}^{\tilde{\alpha}^{2}r}\left.g_{r}\right|_{a_{0}/2}\gamma_{p},
\end{align*}
and thus $\left.g_{r}\right|_{a_{0}/2}\gamma_{p}\in q^{\frac{s}{p}}\mathbb{Z}[[q]]$ with $s=\tilde{\alpha}^{2}r\pmod{p}$.
Therefore, if one writes
$$
f_{\sigma_{p}(a_{0},a_{1},\ldots,a_{(p-1)/2})}(\tau/p)=\sum_{s=0}^{p-1}h_{s}(\tau),
$$
where $h_{s}(\tau)\in q^{\frac{s}{p}}\mathbb{Z}[[q]]$, then 
$$
h_{s}(\tau)=\lambda(\gamma_{p},\vec{a}_{p})^{-1}\left.g_{r}\right|_{a_{0}/2}\gamma_{p},
$$
for $s\equiv\tilde{\alpha}^{2}r\pmod{p}$.
\end{proof}

The operator $U_{p}$ defined as in Lemma~\ref{fgh} indeed acts linearly on the vector space $M_{k}(\Gamma_{1}(p))$. See, e.g., \cite{AL}.

\begin{Corollary}
Let $f_{\vec{a}_{p}}\in M_{a_{0}/2}(\Gamma_{1}(p))$ be a product defined as in~\eqref{faaa}. Then $f_{\vec{a}}$ is an eigenvector of $U_{p}$ if and only if $f_{\sigma_{p}(\vec{a}_{p})}$ is an eigenvector of $U_{p}$.
\end{Corollary}

\begin{proof}
This follows immediately from Lemma~\ref{fgh}.
\end{proof}

\begin{Corollary}
\label{ufuf}
  Let $f_{\vec{a}_{p}}=f_{a_{0},\ldots,a_{(p-1)/2}}(\tau)\in M_{a_{0}/2}(\Gamma_{1}(p))$, and $\gamma_{p}=\begin{pmatrix}\tilde{\alpha}&p\\p\beta&\delta\end{pmatrix}\in\Gamma_{0}(p)$. Then one has that
  $$
U_{p}^{n}\left(f_{\sigma_{p}(\vec{a}_{p})}\right)=\lambda(\gamma_{p},\vec{a}_{p})^{-1}\left.U_{p}^{n}\left(f_{\vec{a}_{p}}\right)\right|_{a_{0}/2}\gamma_{p}.
$$
\end{Corollary}

\begin{proof}
By Lemma~\ref{fgh},
$$
U_{p}\left(f_{\sigma_{p}(\vec{a}_{p})}\right)=\lambda(\gamma_{p},\vec{a}_{p})^{-1}\left.U_{p}\left(f_{\vec{a}_{p}}\right)\right|_{a_{0}/2}\gamma_{p},
$$
where $U_{p}(f_{\vec{a}_{p}})$ lies in $M_{k}(\Gamma_{1}(p))$.
By Corollary~\ref{tf}, $U_{p}(f_{\vec{a}_{p}})$ is a linear combination of products $f_{\vec{a}_{p,j} } \in M_{a_0/2}(\Gamma_{1}(p))$, i.e.,
$$
U_{p}(f_{\vec{a}_{p}})=\sum_{j}c_{j}f_{\vec{a}_{p,j}},
$$
and thus, by Lemma~\ref{flf},
\begin{align*}
\left.U_{p}(f_{\vec{a}_{p}})\right|_{a_{0}/2}\gamma_{p}&=\sum_{j}c_{j}\left.f_{\vec{a}_{p,j}}\right|_{a_{0}/2}\gamma_{p}=\sum_{j}c_{j}\lambda(\gamma_{p},\vec{a}_{p,j})^{p}f_{\sigma_{p}(\vec{a}_{p,j})}.
\end{align*}
Applying $U_{p}$ to both sides, 
\begin{align*}
U_{p}(U_{p}(f_{\vec{a}_{p}})|_{\frac{a_0}{2}}\gamma_{p})
&=\sum_{j}c_{j}\lambda(\gamma_{p},\vec{a}_{p,j})^{p}U_{p}(f_{\sigma_{p}(\vec{a}_{p,j})})
\\
&=\sum_{j}c_{j}\lambda(\gamma_{p},\vec{a}_{p,j})^{p-1}\left.U_{p}(f_{\vec{a}_{p,j}})\right|_{a_{0}/2}\gamma_{p}\\
&=\left.U_{p}\left(\sum_{j}c_{j}f_{\vec{a}_{p,j}}\right)\right|_{a_{0}/2}\gamma_{p}\\
&=\left.U_{p}^{2}(f_{\vec{a}_{p}})\right|_{a_{0}/2}\gamma_{p},
\end{align*}
where the cancellation of $\lambda(\gamma_{p},\vec{a}_{p,j})^{p-1}$ follows from Lemma~\ref{flf}. Therefore, 
$$
U_{p}^{2}(f_{\sigma_{p}(\vec{a}_{p})})=\left.\lambda(\gamma_{p},\vec{a}_{p})^{-1}U_{p}^{2}(f_{\vec{a}_{p}})\right|_{a_{0}/2}\gamma_{p},
$$
and inductively, one obtains the desired relation.
\end{proof}

The preceding results allow us to provide a proof of Theorem~\ref{z1}, rephrased as Corollary~\ref{ctt}.
\begin{Corollary} \label{ctt}
For primes $5 \le p \le 19$, action of the permutation $\sigma_{p}$ on products $f_{\vec{a}} \in M_{k}(\Gamma_{1}(p))$ preserves Ramanujan type congruences. 
\end{Corollary}
\begin{proof}
For $5\leq p\leq19$ and $k\geq2$, with the aid of SageMath, one can obtain Fourier coefficients for generators of the graded $\mathbb{Z}$-algebra of holomorphic modular forms of weight~$k\geq2$ for $\Gamma_{1}(p)$, which is the union of bases for the $\mathbb{Z}$-module of holomorphic modular forms of weights~2 and~3 for $\Gamma_{1}(p)$, say, $\{h_{1},\ldots,h_{m}\}$. By Corollary~\ref{tf}, one can also obtain explicit $\mathbb{C}$-bases for modular forms of weights~2 and~3 for $\Gamma_{1}(p)$ in terms of Klein form products $f_{\vec{b}}$. Suppose the union of these two bases is $\{f_{\vec{b}_{1}},\ldots,f_{\vec{b}_{m}}\}$. One can check that for each $5\leq p\leq19$ and any $1\leq i\leq m$,
$$
h_{i}=\sum_{j=1}^{m}c_{j}f_{\vec{b}_{j}},
$$
where $c_{j}\in \mathbb{Z}_{(p)}$, and $\mathbb{Z}_{(p)}$ denotes the localization of $\mathbb{Z}$ at the ideal $(p)$, and thus, by Lemma~\ref{flf}, one can deduce that for $k=2$ or~3, 
\begin{equation}
    \label{hff}
h_{i}|_{k}\gamma_{p}=\sum_{j=1}^{m}c_{j}\left.f_{\vec{b}_{j}}\right|_{k}\gamma_{p}=
\sum_{j=1}^{m}c_{j}\lambda(\gamma_{p},\vec{b}_{j})f_{\sigma_{p}(\vec{b}_{j})}\in \mathbb{Z}_{(p)}[[q]].
\end{equation}

Now suppose that $f_{\vec{a}_{p}}$ with $\vec{a}_{p}=(a_{0},\ldots,a_{(p-1)/2})$ and $a_{0}\geq4$ satisfies a Ramanujan type congruence modulo~$p$. Then, for any $n\geq0$,
$$
\frac{1}{p^{\alpha n}}U_{p}^{n}(f_{\vec{a}_{p}})\in \mathbb{Z}[[q]]
$$
for some $\alpha\geq1$.
Since $U_{p}$ is a linear operator on $M_{k}(\Gamma_{1}(p))$, it follows that $\frac{1}{p^{\alpha n}}U_{p}^{n}(f_{\vec{a}_{p}})$ lies in the $\mathbb{Z}$-module of holomorphic modular forms of weight~$a_{0}/2$ for $\Gamma_{1}(p)$, and thus, is a $\mathbb{Z}$-linear combination of monomials in $h_{i}$. Hence,
$$
\frac{1}{p^{\alpha n}}U_{p}^{n}(f_{\vec{a}_{p}})=\sum_{(i_{1},\ldots,i_{l})}d_{i_{1},\ldots,i_{l}}h_{i_{1}}\cdots h_{i_{l}},
$$
where $d_{i_{1},\ldots,i_{l}}\in\mathbb{Z}$. By~\eqref{hff}, one can further deduce that
\begin{align}
\frac{1}{p^{\alpha n}}\left.U_{p}^{n}(f_{\vec{a}_{p}})\right|_{a_{0}/2}\gamma_{p}&=\sum_{(i_{1},\ldots,i_{l})}d_{i_{1},\ldots,i_{l}}(h_{i_{1}}\cdots h_{i_{l}})|_{a_{0}/2}\gamma_{p}\nonumber\\
\label{uhhh}
&=\sum_{(i_{1},\ldots,i_{l})}d_{i_{1},\ldots,i_{l}}(h_{i_{1}}|_{k_{1}}\gamma_{p})\cdots (h_{i_{l}}|_{k_{l}}\gamma_{p})\in \mathbb{Z}_{(p)}[[q]],
\end{align}
where $k_{\ell}=2$ or~3 such that $k_{1}+\cdots+k_{l}=a_{0}/2$.
Finally, by Corollary~\ref{ufuf}, 
$$
\lambda(\gamma_{p},\vec{a}_{p})\frac{1}{p^{\alpha n}}U_{p}^{n}(f_{\sigma_{p}(\vec{a}_{p})})=\frac{1}{p^{\alpha n}}\left.U_{p}^{n}(f_{\vec{a}_{p}})\right|_{a_{0}/2}\gamma_{p}\in \mathbb{Z}[1/p][[q]]\cap \mathbb{Z}_{(p)}[[q]]=\mathbb{Z}[[q]].
$$
Therefore, $\frac{1}{p^{\alpha n}}U_{p}^{n}(f_{\sigma_{p}(\vec{a}_{p})})\in \mathbb{Z}[[q]]$, since $\lambda(\gamma_{p},\vec{a}_{p})=\pm1$ by Lemma~\ref{flf}.
\end{proof}

\begin{Remark}
 Theorem~\ref{z1} implies that if $f(\tau)\in \mathbb{Z}[[q]]\cap M_{k}(\Gamma_{1}(p))$ for $5\leq p\leq19$ and $k\geq2$, then the expansions of $f(\tau)$ at the cusps $r/p$ with $(r,p)=1$ are in $\mathbb{Z}[[q]]$. This is not the case for general $\Gamma_{1}(N)$. For example, one can check that the expansion at the cusp $\frac{3}{16}$ of $\eta(4\tau)^{6}$ lies in $\mathbb{Z}[i][[q]]$, while the modular form itself lies in $\mathbb{Z}[[q]]\cap M_{3}(\Gamma_{1}(16))$.
\end{Remark}

The decomposition of $f_{\sigma_{p}(\vec{a})}(\tau /p)$ is determined by the action of $\sigma_{p}$ on each component.

\begin{Corollary}
\label{fgh2}
Suppose that
\begin{align*}
  f_{a_{0}, \ldots, a_{(p-1)/2}} \in {M}_{k}(\Gamma_{1}(p))
\end{align*}
has the decomposition
$$
f_{a_{0},a_{1},\ldots,a_{(p-1)/2}}(\tau/p)=\sum_{r=0}^{p-1}g_{r}(\tau),
$$
where 
$
g_{r}(\tau)\in q^{\frac{r}{p}}\mathbb{Z}[[q]]$
is a $\mathbb{Q}$-linear combination of $f_{a_{0},a_{1,r},\ldots,a_{(p-1)/2,r}}(\tau)\in q^{\frac{r}{p}}\mathbb{Z}[[q]]$.
Then
$$
f_{\sigma_{p}(a_{0},a_{1},\ldots,a_{(p-1)/2})}(\tau/p)=\lambda(\gamma_{p},\vec{a}_{p})^{-1}\sum_{s(r)=0}^{p-1}h_{s(r)}(\tau),
$$
where
$
h_{s(r)}(\tau)\in q^{\frac{s(r)}{p}}\mathbb{Z}[[q]]
$
is a $\mathbb{Q}$-linear combination of
$
\lambda(\gamma_{p},\vec{a}_{p,r})^{p} f_{\sigma_{p}(\vec{a}_{p,r})}(\tau)\in q^{\frac{s(r)}{p}}\mathbb{Z}[[q]]
$
with $\vec{a}_{p,r}=(a_{0},a_{1,r},\ldots,a_{(p-1)/2,r})$ and $s(r)\equiv \tilde{\alpha}^{2}r\pmod{p}$. 
 In particular, for $5\leq p\leq 19$ a prime, $g_{r}(\tau)\equiv0\pmod{p}$ if and only if $h_{s(r)}(\tau)\equiv0\pmod{p}$.
\end{Corollary}

\begin{proof}
The first assertion follows from Lemmas~\ref{flf} and~\ref{fgh}, and the second one follows from the arguments used in proof of Theorem~\ref{z1} and the fact that $g_{r}(\tau)^{2p}$ lies in the $\mathbb{Z}$-module of modular forms for $\Gamma_{1}(p)$.
\end{proof}

 \section{Eigenspace characterization for Ramanujan type congruences} \label{eigen}
Since the operator $U_{N}$
acts as  a linear transformation on
$M_{k}(\Gamma_{1}(N))$, we can explicitly construct the diagonalizable matrix of
linear transformation for $U_{p}$ in terms of ordered bases. For the case $p=5$ and weight $2$, if we act on the left, the matrix of linear transformation for $U_{5}$ with respect to the ordered basis $\langle x_{0}^{2}, x_{0}x_{1}, x_{1}^{2} \rangle$ for $M_{2}(\Gamma_{1}(5))$, where $x_{n} = f_{\sigma_{5}^{n}(2,-3,2)}$, is $$
A_{2,5}=\left(
\begin{array}{ccc}
 1 & 0 & 0 \\
 22 & 5 & -22 \\
 0 & 0 & 1 \\
\end{array}
\right).
$$
Write $E_{k,\lambda}$ for the eigenspace of eigenvalue~$\lambda$, the eigenspace decomposition with respect to $U_{5}$ is given by $M_{2}(\Gamma_{1}(5))=E_{2,1}\oplus E_{2,5}$, where
the eigenspace $E_{2,1}$ is spanned by
$\begin{pmatrix}1&0&1\end{pmatrix}^{t}$ and
$\begin{pmatrix}-2&11&0\end{pmatrix}^{t}$, and $E_{2,5}$ is spanned by
$\begin{pmatrix}0&1&0\end{pmatrix}^{t}$. As a result, 
$$
U_5 \left ( q\frac{(q^{5}; q^{5})_{\infty}^{5}}{(q;q)_{\infty}} \right ) = U_{5}(f_{4,-1,-1})=U_{5}(x_0x_1)=5f_{4,-1,-1}.
$$
It follows that $f_{4,-1,-1}$ satisfies a Ramanujan type congruence modulo $5$. Similarly, we may show that $f_{6,1,-4}$ is an eigenfunction for $U_{5}$ with eigenvalue $5^{2}$. This method does not apply universally since not all products satisfying Ramanujan type congruences modulo $p$ are eigenfunctions for $U_{p}$. The lowest weight example involving a quotient of Klein forms at level $5$ is $f_{8,3,-7}$. Lemma \ref{upeee} provides a sufficient condition for $f\in M_{k}(\Gamma_{1}(p))$ to satisfy a Ramanujan type congruence modulo~$p$.
\begin{Lemma}
\label{upeee}
Suppose that $M_{k}(\Gamma_{1}(p))$ has a decomposition with respect to $U_{p}$ 
{containing the subspace} $\bigoplus_{j\geq1}E_{k,p^{j}\lambda_{j}}$, where $E_{k,\alpha}$ is the eigenspace associated to eigenvalue $\alpha$, and $\lambda_{j}$ are algebraic integers, and in particular, $E_{k,p^{j}\lambda_{j}}$ is spanned by elements in $\mathbb{Z}_{(p)}[\lambda_{j}][[q]]$, with $\mathbb{Z}_{(p)}$ denoting the localization of $\mathbb{Z}$ at the prime ideal $(p)$. Then for any $f\in M_{k}(\Gamma_{1}(p))$ with $f\in\mathbb{Z}[[q]]$, if
$$
f=\sum_{j\geq1}\sum_{h_{j}\in {\rm basis}(E_{k,p^{j}\lambda_{j}})}c_{h_{j}}h_{j},
$$
with $c_{h_{j}}\in\mathbb{Z}_{(p)}[\lambda_{j}]$, one has that $U_{p}^{n}(f)\equiv0\pmod{p^{n}}$ for any nonnegative integer~$n$.

\end{Lemma}

\begin{proof}
Assume that 
$$
f=\sum_{j\geq1}\sum_{h_{j}\in {\rm basis}(E_{k,p^{j}\lambda_{j}})}c_{h_{j}}h_{j},
$$
with $c_{h_{j}}\in\mathbb{Z}_{(p)}[\lambda_{j}]$. Then under the assumptions, one has that for any nonnegative integer $n$,
\begin{align*}
U_{p}^{n}(f)&=U_{p}^{n}\left(\sum_{j\geq1}\sum_{h_{j}\in {\rm basis}(E_{k,p^{j}\lambda_{j}})}c_{h_{j}}h_{j}\right)\\
&=\sum_{j\geq1}p^{nj}\lambda_{j}^{n}\sum_{h_{j}\in {\rm basis}(E_{k,p^{j}\lambda_{j}})}c_{h_{j}}h_{j}.
\end{align*}
Since $U_{p}^{n}(f)\in\mathbb{Z}[[q]]$, while the double sum lies in $p^{n}\mathbb{Z}_{(p)}[\lambda_{j}][[q]]$, with $\lambda_{j}$ an algebraic integer, then $\mathbb{Z}\cap p^{n}\mathbb{Z}_{(p)}[\lambda_{j}]=p^{n}\mathbb{Z}$, and thus, $U_{p}^{n}(f)\in p^{n}\mathbb{Z}[[q]]$.
\end{proof}

To prove that each product from Theorem \ref{mainth1} satisfies a Ramanujan type congruence modulo $p$, Lemma \ref{upeee} may be applied in a similar way. We demonstrate the technique for the cases of Theorem \ref{mainth1} in which $p=13$.



  \begin{Lemma} \label{eigv} Write $E_{k,\lambda}$ for the eigenspace of eigenvalue $\lambda$ of $M_{k}(\Gamma_{1}(p))$ with respect to $U_{p}$. Then for $5 \le p \le 11,$
  \begin{align*}
       M_{2}(\Gamma_{1}(p))&=E_{2,1}\oplus E_{2,p},
  \end{align*}
  and for $13 \le p \le 19$
  \begin{align*}
  M_{2}(\Gamma_{1}(13)) &= E_{2,1}\oplus E_{2,13} \bigoplus_{j=1}^{2} E_{2,\lambda_{j,13}}, \\ 
      M_{2}(\Gamma_{1}(17)) &= E_{2,1}\oplus E_{2,17} \bigoplus_{j=1}^{4} E_{2,\lambda_{j,17}}, \\ 
      M_{2}(\Gamma_{1}(19)) &= E_{2,1}\oplus E_{2,19}\bigoplus_{j=1}^{6} E_{2,\lambda_{j,19}},
  \end{align*}
    where $\lambda_{j,13}$ are the roots of $x^{2} + 5x + 13$ and $\lambda_{j,17}$, $\lambda_{j,19}$ are respectively distinct roots of irreducible quartic and sextic polynomials over $\Bbb Z$.
     Likewise, 
     \begin{align*}
     M_{3}(\Gamma_{1}(7))&=E_{3,1}\oplus E_{3,-7}\oplus E_{3,7^{2}},\\
      M_{3}(\Gamma_{1}(11)) &=E_{3,1}\oplus E_{3,-11} \oplus E_{3,11^{2}} \oplus \bigoplus_{i=1}^{4}E_{3,\lambda_{i,11}},\\
      M_{3}(\Gamma_{1}(13)) &= E_{3,1}\oplus E_{3,13^{2}} \bigoplus_{i=1}^{2} E_{3,\lambda_{i,13}}\bigoplus_{j=1}^{4} E_{3,\delta_{j,13}}, \\
      M_{3}(\Gamma_{1}(17)) &= E_{3,1}\oplus E_{3,17^{2}} \bigoplus_{i=1}^{8} E_{3,\lambda_{i,17}}\bigoplus_{j=1}^{8} E_{3,\delta_{j,17}},
      \end{align*}
      \begin{align*}
      M_{3}(\Gamma_{1}(19)) &= E_{3,1}\oplus E_{3,-19} \oplus E_{3,19^{2}}\bigoplus_{i=1}^{2} E_{3,\lambda_{i,19}}\bigoplus_{j=1}^{6} E_{3,\delta_{j,19}}\bigoplus_{k=1}^{12} E_{3,\kappa_{k,19}},
     \end{align*}
      where  $\lambda_{i,11}$ are roots of a quartic polynomial, $\lambda_{i,13}$ are the roots of $x^{2} + 13x + 13^2$; $\delta_{j,13}$ are the roots of $x^4 - 8x^3 + 2^{3}\cdot13x^2 - 2^{3}\cdot 13^{2}x + 13^{4}$; 
    and where the remaining eigenvalues 
      parameterize distinct roots of irreducible factors of the characteristic polynomial for $U_{p}$.
  \end{Lemma}
 \begin{Remark}
    Eigenvalues not explicitly given in Lemma \ref{eigv} may be computed from the the matrix of linear transformation for $U_p$. This matrix may be constructed directly from the Fourier expansions of the corresponding basis appearing in the proof of Corollary \ref{tf}. If $B_{i}(q)=\sum_{n=0}^{\infty} b_{i}(n)q^{n}$ is the $i^{\text{th}}$ element of the ordered basis of dimension $d$, then the matrix of linear transformation equals $m^{-1}w$, where $m$ is the matrix $m=(b_{j}(n))_{j,n}^{T}$ and $w=(b_{j}(pn))_{j,n}^{T}$, $1\ \le j \le d$, and $0\le n \le d-1$.
    
\end{Remark}
\subsection{A proof of Theorem \ref{mainth1}}  
  A proof is now given for Theorem~\ref{mainth1} that makes use of the structure from Lemma \ref{eigv}.
 We prove the congruences for the modulus $p=13$. The congruences for the other cases are similarly shown. By Theorem~\ref{z1}, to prove congruences associated with $\sigma_{p}^{j}(\vec{a}_{p})$, one only has to verify for the case $j=0$.  For 
$\vec{a}_{13}=(a_{0},\ldots,a_{6})=(6,-3,-1,2,-3,2,0),$
  one can check that $f_{\vec{a}_{13}}$ lies in $E_{3,13^{2}}\oplus E_{3,13\xi}\oplus E_{3,13\overline{\xi}}$, where $\xi = e^{2 \pi i/3}$. With respect to the ordered basis given in Corollary~\ref{tf}, the space $E_{3,13^{2}}$ is spanned by ${\vec{v}_{i}\! }^{t}$ where 
  $$
  \begin{pmatrix}
   \vec{v}_{1} \\ \vec{v}_{2} \\ \vec{v}_{3} \\ \vec{v}_{4} \\ \vec{v}_{5} \\ \vec{v}_{6}
  \end{pmatrix} =
  \left(
\begin{array}{cccccccccccccccccccc}
 0 & 31 & 13 & 0 & 8 & -4 & -6 & -8 & 0 & 0 & 0 & 0 & -4 & 0 & 0 & -2 & 0 & 0 & 1 & 0 \\
 0 & -17 & -6 & 0 & -5 & -1 & 0 & 4 & 0 & 0 & 0 & 0 & 2 & 0 & 0 & 1 & 0 & 1 & 0 & 0 \\
 0 & -6 & -1 & 0 & -2 & 3 & 0 & 2 & 0 & 0 & 0 & 0 & 1 & 0 & 1 & 0 & 0 & 0 & 0 & 0 \\
 0 & -5 & -2 & 0 & -2 & -2 & -2 & 1 & 0 & 0 & 0 & 1 & 0 & 0 & 0 & 0 & 0 & 0 & 0 & 0 \\
 0 & -1 & 0 & 0 & -1 & 1 & 0 & 1 & 0 & 0 & 1 & 0 & 0 & 0 & 0 & 0 & 0 & 0 & 0 & 0 \\
 0 & -2 & -2 & 0 & -1 & -2 & 1 & 0 & 1 & 0 & 0 & 0 & 0 & 0 & 0 & 0 & 0 & 0 & 0 & 0 \\
\end{array}
\right)$$ 
The spaces $E_{3,13\xi}$ and $E_{3,13\overline{\xi}}$ are spanned respectively by ${\vec{v}_{i}\! }^{t}$, $7 \le i \le 10$, where $\vec{v}_{9}, \vec{v}_{10}$ are the complex conjugates, of $$\vec{v}_{7} =\left ( \begin{array}{l} 0,\ -4+8 i \sqrt{3},\ 2+4 i \sqrt{3},0,\ 1+3 i \sqrt{3},\ 7+i \sqrt{3},\ -3-4 i \sqrt{3},\ -\frac{1}{2} i \left(\sqrt{3}-5 i\right), \\  -10-6 i \sqrt{3},\ 0,\ \frac{1}{2} \left(1+i \sqrt{3}\right),\ -\frac{1}{2} i \left(\sqrt{3}-i\right),\ \frac{1}{2} i \left(\sqrt{3}+i\right),\ 0,\ \frac{1}{2} \left(1-9 i \sqrt{3}\right),\\  \frac{1}{2} \left(-1-3 i \sqrt{3}\right),\ 0,\ 0,\ 1,\ 0 \end{array} \right ),$$
$$\vec{v}_{8} =\left ( \begin{array}{l} 0,\ 3 i \sqrt{3},\ 1+2 i \sqrt{3},\ 0,\frac{1}{2} \left(-1+3 i \sqrt{3}\right),\ 1+i \sqrt{3},\ \frac{1}{2} \left(1-3 i \sqrt{3}\right),\ -1,-3-2 i \sqrt{3},\ 0,\\ \frac{1}{2} \left(3-i \sqrt{3}\right),\ -1,\ 0,\ 0,\frac{1}{2} \left(-1-3 i \sqrt{3}\right),\ -\frac{1}{2} i \left(\sqrt{3}-i\right),\ 0,\ 1,\ 0,\ 0\end{array} \right ).$$
The basis representation for $f_{\vec{a}_{13}}$ is a linear combination of the eigenvectors  $\sum_{i=1}^{10} c_{i} \vec{v}_{i},$
  where
$$(c_{i})_{i=1}^{6} = \left (-\frac{28}{183}, -\frac{98}{183}, \frac{86}{183}, \frac{124}{183}, -\frac{131}{183}, -\frac{212}{183}
\right) \in \Bbb Z_{(13)}^{6}$$ and, with $\xi = e^{2\pi i/3}$, $(c_{i})_{i=7}^{10}$ are
\begin{align*}
%
\left (\frac{2}{183} \left(7-2 i \sqrt{3}\right), \frac{1}{366} \left(-85+33 i \sqrt{3}\right), \frac{2}{183} \left(7+2 i \sqrt{3}\right), \frac{1}{366} \left(-85-33 i \sqrt{3}\right) \right)\in
                                                                             (\Bbb
                                                                             Z_{(13)}[\xi])^{4}.                                                   
  \end{align*}
 The desired conclusion follows from Lemma~\ref{upeee}.
 
To prove the Ramanujan type congruence for $\vec{a}_{p}=(a_{0},\ldots,a_{(p-1)/2})=(6,1,\vec{\bf 0}_{(p-5)/2},-4)$, one can check that $f_{\vec{a}_{p}}$ lies in $E_{3,p^{2}}$, from which the corresponding Ramanujan type congruence  follows. For the remaining  $\vec{a}_{p}$'s, one can check that each of the corresponding $f_{\vec{a}_{p}}$ lies in $E_{k,p}$ for $k=2,3$.

After generating a complete set of lattice points satisfying the system in \eqref{ineq}, one can show directly that products with exponent vectors in the polytope not indicated in Theorem \ref{mainth1} fail to satisfy a Ramanujan type congruence. This concludes the proof of Theorem \ref{mainth1} for the case $p=13$.

\section{$p$-Dissections and Bases for level $p$ subgroups} \label{Spdissect}

In Lemma \ref{mdecomp}, the vector space $M_{k}(\Gamma(p))$ is decomposed into subspaces indexed by $r$ modulo $p$, namely
\begin{align}
M_{k}(\Gamma(p))=\bigoplus_{r=0}^{p-1}M_{k,r}(\Gamma(p)), \quad M_{k,r}(\Gamma(p))=\{g\in M_{k}(\Gamma(p)):\,\,g\in q^{\frac{r}{p}}\mathbb{Z}[[q]]\},\label{deco}    
\end{align}
and where the subspaces are representable in terms of quotients of Klein forms. Since, for primes $5 \le p \le 19$, the graded algebra of modular forms for $\Gamma(p)$ is generated by weight one forms \cite{KM}, decomposition formulas of weight one determine explicit expansions for dissections of quotients of Klein forms via
\begin{align} \label{f1}
    U_{p,r} \left ( \sum_{i=1}^{\dim (M_{k}(\Gamma(p))} c_{i}f_{a_{0,i}, a_{1,i}, \ldots,  a_{(p-1)/2,i}} \right ) = \sum_{p\cdot \ell(a_{0,i}, a_{1,i}, \ldots,  a_{(p-1)/2,i}) \equiv r \mod{p}} c_{i}f_{a_{0,i}, a_{1,i}, \ldots,  a_{(p-1)/2,i}}.
\end{align}
The map $\tau \mapsto \tau/p$ sends a modular form for $\Gamma_{1}(p)$ to a corresponding modular form of the same weight for $\Gamma(p)$. Therefore, representations for the image of basis elements $x_{i} = f_{\vec{a_{i}}}$ of $M_{1}(\Gamma_{1}(p))$ in $M_{1}(\Gamma(p))$ under the map $\tau \mapsto \tau/p$ permit an application of \eqref{f1} to derive explicit representations for dissections of products $f\in M_{k}(\Gamma_{1}(p))$ of the form \eqref{faaa}. The resulting identities witness congruences for products of Klein forms in $M_{k}(\Gamma_{1}(p))$ that are not apparent from an analysis of the the eigenstructure of the Hecke matrices in the last section. In this section, we derive the required weight one decomposition identities for prime levels $5 \le p \le 19$.

Several features in each of the next dissection lemmas are characteristic of the decompositions of higher level. Subspaces corresponding to indices that are quadratic resides or non-residues, respectively, have the same form. Moreover, coefficients in the decomposition formulas for weight one generators are permuted in each formula up to a change in sign. This will allow us to present corresponding dissection formulas for larger primes in an abbreviated form.

\begin{Lemma}[5-Dissection] \label{basis5a}
Define  for $0 \le j
\le 5$, $ v_{j} =  K_{5,1}^{j-3}K_{5,2}^{2-j}$.
Then
  \begin{enumerate}
\item  $\langle v_{0}, v_{5} \rangle$ forms a basis for $M_{1}(\Gamma_{1}(5))$. \label{two}
\item  Components of $\langle v_{j} \rangle_{j=0}^{5} $ form a basis for
  $M_{1}(\Gamma(5))$, and  $v_{j}\in q^{j/5}\mathbb{Z}[[q]]$.\label{5one}
  \item  Components of $\langle v_{0}, v_{5}\rangle$ generate the graded
  algebra of modular forms for $\Gamma_{1}(5)$, and components $\langle v_{j} \rangle_{j=0}^{5}$ generate the graded
  algebra of modular forms for $\Gamma(5)$.
\item The following decomposition formulas hold
  \begin{align*}
  v_{0}(\tau/5) &= v_{0}(\tau) +3v_{1}(\tau) +4v_{2}(\tau) + 2v_{3}(\tau)+ v_{4}(\tau),  \\
  v_{5}(\tau/5) &= v_{5}(\tau) +v_{1}(\tau) -2v_{2}(\tau) + 4v_{3}(\tau)-3v_{4}(\tau).     
  \end{align*}
\end{enumerate}

\end{Lemma}

\begin{Lemma}[7-Dissection] \label{decomp}
Let $x_{n}$ for $0\leq n\leq2$ be defined by $x_{n} = f_{\sigma_{7}^{n}(2,-2,0,1)}$.
Define 
\begin{align*} 
& \vec{v}_{0} = \left\langle\begin{array}{ccc}
  x_{0}&x_{1}&x_{2}\end{array}\right\rangle,\quad \vec{v}_{1} =\left\langle
\frac{K_{7,2}}{K_{7,1}^{2}}\quad\frac{K_{7,1}}{K_{7,2}K_{7,3}}\right\rangle,\quad \vec{v}_{2} =\left\langle
  \frac{K_{7,1}}{K_{7,3}^{2}}\quad\frac{K_{7,3}}{K_{7,1}K_{7,2}}\right\rangle, \\
 &             \vec{v}_{3} = \left\langle \frac{1}{K_{7,1}}\right\rangle,\quad \vec{v}_{4}=\left\langle
    \frac{K_{7,3}}{K_{7,2}^{2}}\quad\frac{K_{7,2}}{K_{7,1}K_{7,3}}\right\rangle,\quad \vec{v}_{5}=\left\langle \frac{1}{K_{7,2}}\right\rangle,\quad
   \vec{v}_{6} =\left\langle \frac{1}{K_{7,3}}\right\rangle.
\end{align*}
\begin{enumerate}
\item Components of $\vec{v}_{0}$ form a basis for $M_{1}(\Gamma_{1}(7))$. \label{two}
\item For $0\le j \le 6$, components of $\vec{v}_{j}$ are elements of $q^{j/7}\mathbb{Z}[[q]]$ and form a basis for
  $M_{1}(\Gamma(7))$. \label{one}
  \item  Components of $\vec{v}_{0}$ generate the graded
  algebra of modular forms for $\Gamma_{1}(7)$, and the components of $\vec v_{j}$ generate the graded
  algebra of modular forms for $\Gamma(7)$.
 \item  If we write $ (b_{1}, \ldots, b_{k})_{j} =
  \vec{v}_{j} \begin{pmatrix}b_{1}& \cdots&b_{k}\end{pmatrix}^{t}$, then
  \begin{align*}
        x_{0}(\tau/7) &=
                                            (1,0,0)_{0}+(2,0)_{1}+(1,3)_{2}+(3)_{3}+(1,1)_{4}+(1)_{5}+(2)_{6}, \\
    x_{1}
    (\tau/7)&=
    (0,1,0)_{0}+(1,1)_{1}+(-2,0)_{2}+(-1)_{3}+(-1,3)_{4}+(2)_{5}+(-3)_{6},
    \\
x_{2}(\tau/7)           &=
    (0,0,1)_{0}+(1,-3)_{1}+(1,-1)_{2}+(2)_{3}+(-2,0)_{4}+(3)_{5}+(-1)_{6}.
  \end{align*}
\end{enumerate}
\end{Lemma}

\begin{Lemma}[11-Dissection]
\label{decomp11}
Let $x_{n}$ for $0\leq n\leq4$ be defined by $x_{n}=f_{\sigma_{11}^{n}(2,-1, -1,0,1,0)}$.
%
Define the following subscripted vectors:
\begin{align*}
  \begin{array}{c}
  \vec{v}_{0}=\left\langle x_{0}\quad x_{1}\quad x_{2}\quad x_{3}\quad x_{4}\right\rangle, \quad
  \vec{v}_{1}=\left  \langle\frac{  K_{11,3} K_{11,5}}{K_{11,1} K_{11,2} K_{11,4}}\quad K_{11,3}^{-1}\quad \frac{
             K_{11,1}^{2}}{K_{11,5}^{2}K_{11,4}}\right \rangle,  \\
   \vec{v}_{2}=\left  \langle \frac{K_{11,1} }{K_{11,4}K_{11,5}}\quad \frac{K_{11,3}}{K_{11,1}K_{11,2}}
                \right \rangle, \quad  \vec{v}_{3}=\left  \langle K_{11,4}^{-1}\quad \frac{ K_{11,3}K_{11,4}}{K_{11,1}K_{11,2}K_{11,5}}\quad \frac{ K_{11,5}^{2}}{K_{11,2}K_{11,3}^{2}} \right \rangle,   \\
  \vec{v}_{4}=\left  \langle K_{11,5}^{-1}\quad \frac{ K_{11,1}K_{11,5}}{K_{11,2}K_{11,3}K_{11,4}}\quad \frac{K_{11,2}^{2}}{K_{11,1}^{2}K_{11,4} }\right \rangle,   \quad\vec{v}_{5}=\left  \langle K_{11,1}^{-1}\quad  \frac{ K_{11,1}K_{11,2}}{K_{11,3}K_{11,4}K_{11,5}}\quad \frac{K_{11,4}^{2} }{K_{11,2}^{2}K_{11,5}}\right \rangle,   \\
  \vec{v}_{6}=\left  \langle \frac{K_{11,5} }{K_{11,2}K_{11,3}}\quad \frac{K_{11,4} }{K_{11,1}K_{11,5}} \right \rangle,  \quad \vec{v}_{7}=\left  \langle \frac{K_{11,3} }{K_{11,1}K_{11,4}}\quad \frac{K_{11,2} }{K_{11,3}K_{11,5}}\right \rangle,   \\
  \vec{v}_{8}=\left  \langle \frac{K_{11,5} }{K_{11,2}K_{11,4}}\quad \frac{K_{11,2}}{K_{11,1}K_{11,3}}\right \rangle, \quad \vec{v}_{9}=\left  \langle \frac{K_{11,2}K_{11,4} }{K_{11,1}K_{11,3}K_{11,5} }\quad K_{11,2}^{-1}\quad \frac{K_{11,3}^{2}}{K_{11,1}^{2}K_{11,4}^{2}}\right \rangle,   \\
       \vec{v}_{10}=\left  \langle \frac{K_{11,4} }{K_{11,2}K_{11,5}}\quad \frac{K_{11,1}}{K_{11,2}K_{11,4}}
                     \right \rangle.   
\end{array}
\end{align*}
\begin{enumerate}
\item Components of $\vec{v}_{0}$ form a basis for $M_{1}(\Gamma_{1}(11))$. \label{two}
\item \label{three} For $0\le j \le 10$, components of $\vec{v}_{j}$ are elements of $q^{j/11}\mathbb{Z}[[q]]$ and form a basis for
  $M_{1}(\Gamma(11))$.
  \item  Components of $\vec{v}_{0}$ generate the graded
  algebra of modular forms for $\Gamma_{1}(11)$, and the components of $\vec v_{j}$  generate the graded
  algebra of modular forms for $\Gamma(11)$.
\item If we write $ (b_{1}, \ldots, b_{k})_{j} =
  \vec{v}_{j} \begin{pmatrix}b_{1}& \cdots&b_{k}\end{pmatrix}^{t}$,
  then
  \begin{align*}
    x_{0}(\tau/11)	&=	(1,0,0,0,0)_{0}+(1,2,1)_{1}+(-1,2)_{2}+(2,1,1)_{3}+(0,-1,2)_{4}
		+(2,0,0)_{5} \\ & \qquad +(2,0)_{6}+(1,0)_{7}+(1,0)_{8}+(1,0,0)_{9}
                                  +(2,1)_{10}, \\
                    x_{1} (\tau/11)   &=	(0,1,0,0,0)_{0}+(1,0,2)_{1}+(0,1)_{2}+(2,-1,1)_{3}+(-2,0,0)_{4}
		+(0,-1,0)_{5} \\ & \qquad +(-1,2)_{6}+(2,0)_{7}+(-1,2)_{8}+(-1,2,-1)_{9}
                                   +(-1,0)_{10}, \\
x_{2}(\tau/11)	&=	(0,0,0,0,1)_{0}+(0,-2,0)_{1}+(-2,1)_{2}+(0,1,-2)_{3}+(0,-1,0)_{4}
		+(2,-1,-1)_{5} \\ & \qquad  +(0,-1)_{6}+(1,-2)_{7}+(0,1)_{8}+(-1,-2,1)_{9}
                   +(2,0)_{10}, \\
    x_{3}(\tau/11)	&=	(0,0,0,1,0)_{0}+(1,0,0)_{1}+(1,0)_{2}+(-2,0,0)_{3}+(-2,1,1)_{4}
		+(2,1,-1)_{5}\\ &\qquad +(-2,1)_{6} +(0,-1)_{7}+(0,2)_{8}+(-1,0,2)_{9}
  +(-1,-2)_{10}, \\
  x_{4}(\tau/11)	&=	(0,0,1,0,0)_{0}+(1,-2,-1)_{1}+(-2,0)_{2}+(0,1,0)_{3}+(2,1,-1)_{4}
		+(0,-1,2)_{5} \\ & \qquad +(-1,0)_{6}+(2,-1)_{7}+(-2,1)_{8}+(0,2,0)_{9}
                   +(0,1)_{10}.
\end{align*}
\end{enumerate}
\end{Lemma}

\begin{proof}[Proofs of Lemmas~$\ref{basis5a}$--$\ref{decomp11}$]
Assertions (1)--(2) are straightforward and follow from known dimension formulas for the vector spaces. Claim (3) is equivalent to Corollary \ref{tf} and follows from (1), (2) and \cite{KM}.
 Since $x_{n}(\tau)$ are modular forms of weight~1 for $\Gamma_{1}(p)$, then $x_{n}(\tau/p)$ lie in $M_{1}(\Gamma(p))$ and must be a linear combination of the components of $\vec v_{j}$. The last assertion follows.
\end{proof}  

The next lemmas generalize the preceding results for larger primes. We start by demonstrating that the space of weight one modular forms for $\Gamma(p)$ decomposes as a direct sum of subspaces of series of the form $q^{\frac{r}{p}}\mathbb{Z}[[q]]$, for each congruence class $r$ modulo $p$. Moreover, subspaces corresponding to quadratic residues (and non-residues, respectively) are isomorphic.
\begin{Lemma} Define $M_{k,r}(\Gamma(p))$ as in \eqref{deco}. Then 
\label{mdecomp}
$$
M_{k}(\Gamma(p))=\bigoplus_{r=0}^{p-1}M_{k,r}(\Gamma(p)).
$$

In particular, $$M_{k,0}(\Gamma(p))=M_{k}(\Gamma_{1}(p)).$$ Moreover, if $(\mathbb{Z}/p\mathbb{Z})^{\times}=\langle\alpha\rangle$, then for $r$ such that $\left(\frac{r}{p}\right)=1$,
$$
M_{k,r}(\Gamma(p))\cong M_{k,1}(\Gamma(p)),
$$
and for $r$ such that $\left(\frac{r}{p}\right)=-1$,
$$
M_{k,r}(\Gamma(p))\cong M_{k,\tilde{\alpha}}(\Gamma(p)),
$$
where $\tilde{\alpha}$ denotes the inverse of $\alpha\pmod{p}$.
\end{Lemma}

\begin{proof}
Clearly, the sum is direct whose summands are all nonempty. For any $f\in M_{k}(\Gamma(p))$, writing 
$$
f(\tau)=\sum_{r=0}^{p-1}\phi_{r}(\tau)
$$
for which $\phi_{r}(\tau)\in q^{\frac{r}{p}}\mathbb{Z}[[q]]$, one can note that
 $$
 \begin{pmatrix}f\left(\tau\right)\\f\left(\tau+1\right)\\f\left(\tau+2\right)\\
 \vdots\\f\left(\tau+p-1\right)\end{pmatrix}=\begin{pmatrix}1&1&1&\cdots&1\\
 1&\zeta_{p}&\zeta_{p}^{2}&\cdots&\zeta_{p}^{p-1}\\ 1&\zeta_{p}^{2}&\left(\zeta_{p}^{2}\right)^{2}&
 \cdots&\left(\zeta_{p}^{p-1}\right)^{2}\\\vdots&\vdots&\vdots&&\vdots\\
 1&\zeta_{p}^{p-1}&\left(\zeta_{p}^{2}\right)^{p-1}&
 \cdots&\left(\zeta_{p}^{p-1}\right)^{p-1}
 \end{pmatrix}
  \begin{pmatrix}\phi_{0}\\\phi_{1}\\\phi_{2}\\\vdots\\\phi_{p-1}\end{pmatrix}
 $$
 where the matrix is the Vandermonde matrix and thus invertible. This implies that $\phi_{m}(\tau)$ are all linear combinations of $f\left(\tau+r\right)$ which are all in $M_{k}(\Gamma(p))$ since $\Gamma(p)$ is a normal subgroup of ${\rm SL}_{2}(\mathbb{Z})$. This proves the first claim. 
For the second assertion, recall from the proof of Lemma~\ref{fgh} that the map $g_{r}\to g_{r}|_{k}\gamma_{p}$ is an isomorphism from $M_{k,r}(\Gamma(p))$ to $M_{k,\tilde{\alpha}^{2}r}(\Gamma(p))$, where $\gamma_{p}=\left ( \begin{smallmatrix}\tilde{\alpha}&p\\p\beta&\delta\end{smallmatrix} \right) $. The second assertion follows.
\end{proof}

Corollary \ref{mdecomp} may be applied to represent the decomposition components for $M_{k}(\Gamma(p))$ in terms of permutations of a generating set of Klein form quotients. For nonzero values of $r$, elements of $M_{k,r}(\Gamma(p))$ are generated under permutative action by two vectors with respective indices corresponding to a quadratic residue and non-residue modulo $p$.

\begin{Corollary}
If $M_{k,1}(\Gamma(p))={\rm Span}\{f_{a_{0},\ldots,a_{(p-1)/2}}\}$ and $M_{k,\tilde{\alpha}}(\Gamma(p))={\rm Span}\{g_{a_{0},\ldots,a_{(p-1)/2}}\}$, then
$$
M_{k,\tilde{\alpha}^{2i}}(\Gamma(p))=\sigma_{p}^{i}M_{k,1}(\Gamma(p)),\qquad M_{k,\tilde{\alpha}^{2i+1}}(\Gamma(p))=\sigma_{p}^{i}M_{k,\tilde{\alpha}}(\Gamma(p)).
$$
\end{Corollary}

\begin{proof}
These follow immediately from Corollary~\ref{fgh2} and Lemma~\ref{mdecomp}.
\end{proof}

The results above may be used to formulate maximally linearly independent sets of elements in the subspaces  $M_{1,r}(\Gamma(p))$ in terms of quotients of products defined by \eqref{faaa}. These products correspond to a subset of the lattice points meeting the requirements of Lemma \ref{m1r}. We will see subsequently that each subspace $M_{1,r}(\Gamma(p))$ of $M_{1}(\Gamma(p))$ may be given in terms of permutations of a pair of vectors $V_{p,1},V_{p,n}$, where $n$ is a quadratic nonresidue modulo $p$.

\begin{Corollary} \label{tf1}
Define $\vec{a}_{p}$ as in Corollary \ref{tf} and denote
\begin{align*}
  V_{13,1} &=  \Bigl \langle  (2,-1, 1, 0, -1, 1, -1), (2,-2, 1, 0, 0, 0, 0), (2,-1, -1, 1, 0, -1, 1)\Bigl
    \rangle, \\
V_{13,7} &=\Bigl \langle (2,1, 0, 0, -1, 0, -1), (2,-2, 2, -1, 1, -1, 0), (2,1, -1, -1, 0, 0, 0)\Bigl
    \rangle,
  \end{align*}
   \begin{align*}
  V_{17,1}&= \left  \langle  \begin{array}{c} (2,0, -1, 0, 0, 0, 1, 0, -1),  (2,0, 0, -1, 0, 0, -1, 0, 1), (2,-1, -1, 1, 0, 0, 0, -1, 1),   \\ [-10pt] (2,-2, 0, 2, -1, 0, 1, -2, 1) \end{array} 
   \right \rangle\\
     V_{17,6} &=\left \langle  \begin{array}{c}  (2,-1, 0, 0, 0, 0, -1, 1, 0), (2,0, 0, 1, -1, 0, 0, -1, 0), (2,0, 0, 0, 1, \
-1, 0, 0, -1),  \\ [-10pt]   (2,0, -1, -1, 0, 1, 0, 0, 0) \end{array} 
   \right  \rangle,
  \end{align*}
   \begin{align*}
  V_{19,1} &= \left   \langle  \begin{array}{c}  (2,0, 0, 0, -1, -1, 0, 0, 0, 1), (2,-1, 1, 0, -1,
    0, 0, -1, 1, 0),   (2,0, -1, 0, 0, 0, 0, 1, 0, -1),  \\ [-10pt]   (2,0, 0, -1, 0, -1, 1, 0, 0, 0), (2,-1, -1, 1, 0, 0, 0, 0, -1, 1) \end{array} 
   \right  \rangle, \\
    V_{19,10} &=  \left  \langle  \begin{array}{c}  (2,-1, 0, 0, 0, -1, 2, 0, -1, 0), (2,0, -1, 0, 1, 0, -1, 1, 0, -1), (2,-1, \
0, 1, 0, 0, 0, 0, -1, 0), \\ [-10pt]   (2,0, -1, 0, 0, 1, 1, -1, -1, 0) \end{array} 
   \right  \rangle.
  \end{align*}
If $[V]_{j}$ represents the $j$-th component of $V$, then the orbits distribute over $q^{r/p}\Bbb Z[[q]]$ as follows:
\begin{align*}
 f_{\sigma_{13}^{n}[{V_{13,1}]_{j}}} &\in M_{1,7 ^{2n}}(\Gamma(13)), \quad f_{\sigma_{13}^{n}[{V_{13,7}]_{j}}} \in M_{1,7^{2n+1}}(\Gamma(13)),\quad 0 \le n \le 5, \quad 1 \le j \le 3, \\ 
  f_{\sigma_{17}^{n}[{V_{17,1}]_{j}}} &\in M_{1,6^{2n}}(\Gamma(17)), \quad f_{\sigma_{17}^{n}[{V_{17,6}]_{j}}}  \in M_{1,6^{2n+1}}(\Gamma(17)),\quad 0 \le n \le 7, \quad 1 \le j \le 4,\\
   f_{\sigma_{19}^{n}[{V_{19,1}]_{j}}}&\in M_{1,10^{2n}}(\Gamma(19)), \quad f_{\sigma_{19}^{n}[{V_{19,10}]_{k}}} \in M_{1,10^{2n+1}}(\Gamma(19)),\quad 0 \le n \le 8,\ 1\le j \le 5,\ 1 \le k \le 4.
\end{align*}
\end{Corollary}

Each containment in the corollary follows directly from the definition of $f_{\vec{a}}$ in \eqref{faaa}. Since the vectors $V_{p,1}$ and $V_{\tilde{\alpha}, p}$ are mapped to $M_{1,\tilde{\alpha}^{2n}}(\Gamma(p))$ and $M_{1,\tilde{\alpha}^{2n+1}}(\Gamma(p))$, respectively, under $\sigma^{n}$, we define the following notation for each subspace generated by the vectors in Corollary \ref{tf1}.
\begin{Definition}  \label{defvrp}  Let $\alpha$ denote the least positive primitive root modulo $p$ and suppose $\tilde{\alpha}$ is the inverse of $\alpha$ modulo $p$. For $0\le n \le (p-3)/2$, define the vectors $V_{p,n}$ by
$$V_{p,[\tilde{\alpha}^{2n}]} = \sigma_{p}^{n}V_{p,1}, \qquad V_{p,[\tilde{\alpha}^{2n+1}]} = \sigma_{p}^{n}V_{p,\tilde{\alpha}},$$
where $[m]=m\pmod{p}$ denotes the least positive residue.
\end{Definition}
The prior corollary and next lemma are needed to determine the decomposition of $M_{1}(\Gamma(p))$ into the subspaces spanned by the corresponding quotients of Klein forms. 
\begin{Lemma}
For $p\leq 19$ an odd prime, one has that
$$
{\rm dim}(M_{1}(\Gamma(p)))=\frac{p^{2}-1}{4}.
$$
\end{Lemma}

\begin{proof}
One can first note that $\Gamma(p)$ is isomorphic to $\Gamma_{1}(p)\cap \Gamma_{0}(p^{2})$ via
$$
\gamma\to \begin{pmatrix}\frac{1}{p}&0\\0&1\end{pmatrix}\gamma  \begin{pmatrix}p&0\\0&1\end{pmatrix},
$$
and thus 
$$
M_{1}(\Gamma(p))\cong M_{1}(\Gamma_{1}(p)\cap\Gamma_{0}(p^{2})).
$$
Clearly, $\Gamma_{1}(p^{2})\subset\Gamma_{1}(p)\cap \Gamma_{0}(p^{2})$, and so, $S_{1}(\Gamma_{1}(p)\cap \Gamma_{0}(p^{2}))\subset S_{1}(\Gamma_{1}(p^{2}))$. Looking up the L-function and Modular Form Database (LMFDB), one can find that $S_{1}(\Gamma_{1}(p^{2}))=\{0\}$ for  all odd primes $p\leq 19$. Therefore, $S_{1}(\Gamma_{1}(p)\cap\Gamma_{0}(p^{2}))=\{0\}$, so that $S_{1}(\Gamma(p))=\{0\}$, which implies that
$$
{\rm dim}(M_{1}(\Gamma(p)))={\rm dim}(E_{1}(\Gamma(p)))=\frac{p^{2}-1}{4}.
$$
\end{proof}

\begin{Lemma}[Khuri-Makdisi \cite{KM}]
  For primes $5\le p \le 19$, the graded ring of modular forms for $\Gamma(p)$ is generated by a basis for $M_{1}(\Gamma(p))$.
\end{Lemma}

The dimension formula for $M_{1}(\Gamma(p))$ along with the previous corollaries allow us to construct product bases for the principal congruence subgroups of level $p$. Each vector space may be described in terms of 2 vector collections $V_{p,1}, V_{p,\alpha}$ from Corollary \ref{tf1} and the vector $\vec{a}_{p}$ from Corollary \ref{tf}.
\begin{Corollary} \label{tf2}
 Let $[\vec{v}]_{i}$ represent the $i$th component of the vector $\vec{v}$. With the notation of Corollaries \ref{tf}, \ref{tf1}, and Definition \ref{defvrp},
\begin{align*}
    M_{1}(\Gamma(13)) &=\bigoplus_{n=0}^{5} \left  ( \Bbb C f_{\sigma_{13}^{n}(\vec{a}_{13})} \oplus\bigoplus_{j=1}^{3} \left(\Bbb C f_{\sigma_{13}^{n}([V_{13,1}]_{j})}\oplus  \Bbb C f_{\sigma_{13}^{n}([V_{13,7}]_{j})}\right) \right ), \\
        M_{1}(\Gamma(17)) &=\bigoplus_{n=0}^{7} \left  ( \Bbb C f_{\sigma_{17}^{n}(\vec{a}_{17})}\oplus \bigoplus_{j=1}^{4} \left(\Bbb C f_{\sigma_{17}^{n}([V_{17,1}]_{j})}\oplus  \Bbb C f_{\sigma_{17}^{n}([V_{17,6}]_{j})}\right) \right ), \\
            M_{1}(\Gamma(19)) &=\bigoplus_{n=0}^{8} \left  ( \Bbb C f_{\sigma_{19}^{n}(\vec{a}_{19})} \oplus \bigoplus_{j=1}^{5} \Bbb C f_{\sigma_{19}^{n}([V_{19,1}]_{j})} \oplus \bigoplus_{k=1}^{4}\Bbb C f_{\sigma_{19}^{n}([V_{19,10}]_{k})} \right ).
\end{align*}
In particular, 
\begin{align*}
     M_{1,7^{2n}}(\Gamma(13))=\bigoplus_{j=1}^{3} \Bbb C f_{\sigma_{13}^{n}([V_{13,1}]_{j})},\quad M_{1,7^{2n+1}}(\Gamma(13))=\bigoplus_{j=1}^{3} \Bbb C f_{\sigma_{13}^{n}([V_{13,7}]_{j})}\\
    M_{1,6^{2n}}(\Gamma(17))=\bigoplus_{j=1}^{4} \Bbb C f_{\sigma_{17}^{n}([V_{17,1}]_{j})},\quad M_{1,6^{2n+1}}(\Gamma(17))=\bigoplus_{j=1}^{4} \Bbb C f_{\sigma_{17}^{n}([V_{17,6}]_{j})}\\
M_{1,10^{2n}}(\Gamma(19))=\bigoplus_{j=1}^{5} \Bbb C f_{\sigma_{19}^{n}([V_{19,1}]_{j})},\quad M_{1,10^{2n+1}}(\Gamma(19))=\bigoplus_{j=1}^{4} \Bbb C f_{\sigma_{19}^{n}([V_{19,10}]_{j})}.
\end{align*}

Moreover, for $p=13, 17, 19$, the graded rings of modular forms for $\Gamma(p)$ are generated by the basis components in the decomposition of $M_{1}(\Gamma(p))$ appearing above.
\end{Corollary}

In the next lemma, we find representations in terms of the basis elements from Corollary~\ref{tf2} for bases from Corollary \ref{tf} under the map $\tau \mapsto \tau/p$. This extends Lemmas \ref{basis5a} to \ref{decomp11} for $p=13,17,19$.
\begin{Lemma}[$p=13,17,19$ Dissection Formulas] \label{pdiss}
  Let $\vec{a}_{p}$ and $\sigma_{p}$ be defined as in Corollary \ref{tf}. Let $V_{p,0} = \langle \sigma_{p}^{n}(f_{\vec{a}_{p}}) \rangle_{n=0}^{(p-1)/2}$. For $0\le n \le (p-3)/2$, let $V_{p,n}$ be as defined in Definition \ref{defvrp}, let $[\vec{V}]_{i}$ represent the $i$th component of the vector $\vec{v}$, and denote $\mathcal{V}_{p,i}^{T} = \langle f_{[V_
  {p,i}]_{j}} \rangle_{j=1}^{|V_{p,i}|}.$ Write the standard basis vector of length $(p-1)/2$ as $e_{1} = (1,0,0,\ldots,0)$. Then 
 \begin{align*}
         f_{\vec{a}_{13}}(\tau/13) = e_{1} \mathcal{V}_{13,0} &+ (0, 1, 0)\mathcal{V}_{13,1} +(3, 1,-2)\mathcal{V}_{13,2}+(2, 2, 1)\mathcal{V}_{13,3}+(1, 1, 1)\mathcal{V}_{13,4} \\ &+ (0, 1, 2)\mathcal{V}_{13,5}+ (1, 0, 2)\mathcal{V}_{13,6}+ (1, 1, 0)\mathcal{V}_{13,7}+ (0, 1, 0)\mathcal{V}_{13,8} \\ &+ (0, 0, 1)\mathcal{V}_{13,9} + (-1, -1, 1)\mathcal{V}_{13,10}+(1, 0, 0)\mathcal{V}_{13,11}+(0, -1, 2)\mathcal{V}_{13,12}, \\
 f_{\vec{a}_{17}}(\tau/17) = e_{1} \mathcal{V}_{17,0} &+ (-1, 1, -2,2)\mathcal{V}_{17,1} +(0,1, 0, 0)\mathcal{V}_{17,2}+(1, 1, 0, 0)\mathcal{V}_{17,3} \\ &+(3, 1, -3, -2)\mathcal{V}_{17,4}+(1, 1, 1, 1)\mathcal{V}_{17,5}+(0, -1, 1, 2)\mathcal{V}_{17,6} \\ &+(0, 1, 1, 0)\mathcal{V}_{17,7}+(0,
        1, 1, 0)\mathcal{V}_{17,8}+(1, 0,2, -1)\mathcal{V}_{17,9} \\ &+(1, 2, 0, -1)\mathcal{V}_{17,10}+(1, 0, 0, 1)\mathcal{V}_{17,11}+(1, 1, 1, -1)\mathcal{V}_{17,12} \\ &+(1,0, -1, 1)\mathcal{V}_{17,13}+(-1, -1, 0, 2)\mathcal{V}_{17,14}+(0, 0, 2, -1)\mathcal{V}_{17,15}+(0, 0, 1, 1)\mathcal{V}_{17,16}, \\
     f_{\vec{a}_{19}}(\tau/19)= e_{1} \mathcal{V}_{19,0} &+ (0, 1, 0, 1, 0)\mathcal{V}_{19,1} +(1, 2, -1, -2)\mathcal{V}_{19,2} +(0, 0, 1,
                                     1)\mathcal{V}_{19,3} \\ &+ (-1, 1, 1, 0, 0)\mathcal{V}_{19,4} +
                                                                       (1,
                                                                       -1,
                                                                       1,
                                                                       -1,
                                                                       0)\mathcal{V}_{19,5}
      + (2, -1, 0, -1,
                                                                       1)\mathcal{V}_{19,6} \\ &+ (0, 1, 0, 1,
                                               0)\mathcal{V}_{19,7} +(1, 0, 0,
             0)V_{19,8}+ (1,
                                               0, 1, -1,
                                               1)\mathcal{V}_{19,9} \\ &+(1, 0, 0, -1)\mathcal{V}_{19,10}+ (0, 1, 0, 0, 1)\mathcal{V}_{19,11}+(0, 1, -2, 0)\mathcal{V}_{19,12} \\ &+(0, 0, 1, 0)\mathcal{V}_{19,13}+(0, 0, 0,
             1)\mathcal{V}_{19,14}+(1, 0, -1,
                                     -1)\mathcal{V}_{19,15} \\ & + (1, -2, 1, 0, -1)\mathcal{V}_{19,16}   + (0,
                                               0, 0, 1,
                                               0)\mathcal{V}_{19,17} 
         +(0, 1, 0, 0)\mathcal{V}_{19,18}.
\end{align*}      
With $\vec{v}_{r}$ as defined above, namely $ f_{\vec{a}_{p}}(\tau /p) = e_{1}\mathcal{V}_{p,0}+ \sum_{r=1}^{p-1} \vec{v}_{r} V_{p,r}$, then for $0\le n \le (p-3)/2$,
\begin{align*}
    f_{\sigma_{p}^{n}(\vec{a}_{p})}(\tau/p)&=\left ( \prod_{i=0}^{n-1}\lambda(\gamma_{p},\sigma_{p}^{i}(\vec{a}_{p}))^{-1} \right )\sum_{r=0}^{p-1}\vec{v}_{r} \left \langle \left (\prod_{i=0}^{n-1}\lambda^{p}(\gamma_{p},\sigma_{p}^{i}([V_{p,r}]_{d})) \right ) f_{\sigma_{p}^{n}([V_{p,r}]_{d})}(\tau)\right \rangle_{1 \le d \le \text{length}(V_{p,r})}^{T},
\end{align*}
where 
 $$
 \lambda(\gamma_{p},\vec{a}_{p})=\prod_{i=1}^{\frac{p-1}{2}}{\rm sgn}(r_{i})^{a_{i}}(-1)^{A(\gamma_{p},\vec{a}_{p})}, \quad  A(\gamma_{p},\vec{a}_{p})=\sum_{i=1}^{\frac{p-1}{2}}\left(q_{i}i+i+q_{i}+\frac{r_{i}i}{p}\right)a_{i}.
 $$
 Here, we define, for $\alpha$ the least positive primitive root modulo $p$, $\tilde{\alpha} \equiv \alpha^{-1}\pmod{p}$ and
 $$
 \gamma_{p} =\begin{pmatrix}
  \tilde{\alpha}&p\\p\beta&\delta 
 \end{pmatrix}\in{\rm SL}_{2}(\mathbb{Z}),
 $$
 and we write $\tilde{\alpha}i=q_{i}p+r_{i}$ with $-\frac{p-1}{2}\leq r_{i}\leq \frac{p-1}{2}$
In particular, $\gamma_{p}$ may be selected as follows: $$\gamma_{13} = \left(
\begin{array}{cc}
 7 & 13 \\
78 & 145 \\
\end{array}
\right), \quad \gamma_{17} = \left(
\begin{array}{cc}
 6 & 17 \\
 85 & 241 \\
\end{array}
\right), \quad \gamma_{19} = \left(
\begin{array}{cc}
 10 & 19 \\
 171 & 325 \\
\end{array}
\right). $$
\end{Lemma}

\begin{proof}

For $x_{n}(\tau)=f_{\sigma_{p}^{n}(a_{0},a_{1},\ldots,a_{(p-1)/2})}(\tau)\in M_{1}(\Gamma_{1}(p))$, suppose that
 $$
 x_{0}(\tau/p)=\sum_{r=0}^{p-1}\vec{v}_{r}V_{r},
 $$
 where
 $$
 V_{r}=(g_{\vec{b}_{p}}(\tau))\subset M_{1,r}(\Gamma(p)).
 $$
By the proofs of Lemma~\ref{fgh} and Corollary~\ref{fgh2},
 one has the following equivalent equations.
 \begin{align*}
 \left.x_{0}(\tau/p)\right)|_{a_{0}/2}^{n}\gamma_{p}&=\sum_{r=0}^{p-1}\vec{v}_{r}\left.V_{r}\right|_{a_{0}/2}^{n}\gamma_{p},\\
 \left.f_{a_{0},a_{1},\ldots,a_{(p-1)/2}}(\tau/p)\right|_{a_{0}/2}^{n}\gamma_{p}&=\sum_{r=0}^{p-1}\vec{v}_{r}\left(\left.g_{\vec{b}_{p}}(\tau)\right|_{a_{0}/2}^{n}\gamma_{p}\right),\\
 \prod_{i=0}^{n-1}\lambda(\gamma_{p},\sigma_{p}^{i}(\vec{a}_{p}))f_{\sigma_{p}^{n}(a_{0},a_{1},\ldots,a_{(p-1)/2})}(\tau/p)&=\sum_{r=0}^{p-1}\vec{v}_{r}\left(\prod_{i=0}^{n-1}\lambda(\gamma_{p},\sigma_{p}^{i}(\vec{b}_{p}))g_{\sigma_{p}^{n}(\vec{b}_{p})}(\tau)\right),\\
 x_{n}(\tau/p)&=\prod_{i=0}^{n-1}\lambda(\gamma_{p},\sigma_{p}^{i}(\vec{a}_{p}))^{-1}\sum_{r=0}^{p-1}\vec{v}_{r}\left(\prod_{i=0}^{n-1}\lambda(\gamma_{p},\sigma^{i}(\vec{b}_{p}))g_{\sigma_{p}^{n}(\vec{b}_{p})}(\tau)\right). 
 \end{align*}
\end{proof}

One can obtain novel congruence between the arithmetic functions $P_{a_{0},\ldots,a_{(p-1)/2}}(n)$ from Lemma~\ref{pdiss} and the preceding results on graded algebras with knowledge of a sufficient set of algebraic relations. 
For the level $5$ case, the $q$-expansions may be used to verify that
the ring of modular forms generated by the
weight one forms in the proof of Corollary~\ref{tf} is algebraically
independent. At levels $7$ and $11$, the relevant relations date to the work of
Klein \cite{IB,K1,K2}.

\begin{Lemma}[Klein's Relations Modulo $7$] \label{Klein7}
\label{kr7}
Let $\vec{a}_{7}$ and $\sigma_{7}$ be as in Corollary \ref{tf}. For $0\le n \le 2$, define $x_{n} = \sigma_{7}^{n-1}(f_{\vec{a}_{7}}).$ Then 
\begin{align}
  \label{k7}
x_{0}x_{1}-x_{1}x_{2}-x_{0}x_{2}=0.
\end{align}
\end{Lemma}

\begin{Lemma}[Klein's Relations Modulo $11$] \label{Klein11}
Let $\vec{a}_{11}$ and $\sigma_{11}$ be as in Corollary \ref{tf}. For $0\le n \le 4$, define $x_{n} = \sigma_{11}^{n-1}(f_{\vec{a}_{11}}).$ Then 
\begin{align}
  \label{eq:1}
    x_0 x_1-x_1 x_4-x_0 x_3 &=0, \\
  x_{0}x_{4}  - x_{0}x_{3} - x_{2}x_{4} &=  0, \\
 x_{1}x_{4} - x_{0}x_{2}  - x_{1}x_{2} &=  0, \\
  x_{1}x_{3}  - x_{0}x_{2} - x_{3}x_{2} &=  0, \\
  x_{1}x_{3}  - x_{4}x_{3} - x_{2}x_{4} &=  0.
\end{align}
\end{Lemma}

For $\Gamma_{1}(13)$ we can show all degree $3$ and $4$ relations among the generators are a consequence of the quadratic relations listed in Lemma \ref{lev13}. The relations may be proven from the Sturm bound. 

\begin{Lemma} \label{lev13} Let $\vec{a}_{13}$ and $\sigma_{13}$ be as in Corollary \ref{tf}. For $0\le n \le 5$, define $x_{n} = \sigma_{13}^{n-1}(f_{\vec{a}_{13}}).$ Then
\begin{align} 
x_{1}x_{3}-x_{3}x_{5}-x_{1}x_{5} &= 0, \\
x_{3}x_{4}-x_{4}x_{5}-x_{0}x_{5} &=0, \\ 
x_{0}x_{4}-x_{0}x_{2}-x_{2}x_{4}&= 0, \\
x_{1}x_{4}+x_{2}x_{4}-x_{3}x_{4}-x_{1}x_{5}-x_{2}x_{5} &=0, \\ 
x_{0}x_{3}+x_{3}x_{5}-x_{0}x_{4}+x_{2}x_{4}-x_{2}x_{5}-x_{4}x_{5} &= 0.
\end{align}
\end{Lemma}

\begin{Remark}
In fact, the geometry behind the relations given in Lemmas~\ref{Klein7}--\ref{lev13} is that they respectively represent the defining equations for the images of the embeddings of their associated modular curves into certain projective space, i.e.,
$$
\Gamma_{1}(7)\backslash\left(\mathbb{H}\cup\mathbb{Q}\cup\{i\infty\}\right)\cong {\rm Proj}\mathbb{C}[x_{0},x_{1},x_{2}],
$$
$$
\Gamma_{1}(11)\backslash\left(\mathbb{H}\cup\mathbb{Q}\cup\{i\infty\}\right)\cong {\rm Proj}\mathbb{C}[x_{0},x_{1},x_{2},x_{3},x_{4}],
$$
and 
$$
\Gamma_{1}(13)\backslash\left(\mathbb{H}\cup\mathbb{Q}\cup\{i\infty\}\right)\cong {\rm Proj}\mathbb{C}[x_{0},x_{1},x_{2},x_{3},x_{4},x_{5}].
$$
A standard argument for these can be found in \cite[p. 8]{I} specialized to the case of the line bundle $\mathcal{L}$ associated with the usual factor of automorphy of weight~1.  By the Riemann--Hurwitz formula, one can find that $\Gamma_{1}(N)\backslash\left(\mathbb{H}\cup\mathbb{Q}\cup\{i\infty\}\right)$ for $N=7,11,13$ are respectively of genus~0,~1 and~2, and on the other hand, 
by Plucker's formula and SAGE, the projective plane curve ${\rm Proj}\mathbb{C}[x_{0},x_{1},x_{2}]$ and the projective curves ${\rm Proj}\mathbb{C}[x_{0},x_{1},x_{2},x_{3},x_{4}],$ and ${\rm Proj}\mathbb{C}[x_{0},x_{1},x_{2},x_{3},x_{4},x_{5}]$  are respectively of genus~0,~1 and~2 as well. 
Moreover, by specializing \cite[Corollary 1]{I} to our case mentioned above, one can tell that relations between $x_{i}$ must come from relations of degree less than or equal to~4. One can show that the relations of degrees~3 and~4 all come from the degree-2 relations in the lemmas.
\end{Remark}

\section{Applications to dissections}
In this section, we employ the decomposition lemmas established in the last section to study dissections of modular forms of levels $5$ and higher whose Fourier coefficients
encode interesting arithmetic information. We first apply the dissection operator to identify all quotients $f$ of the form \eqref{faaa} such that
$U_{5,r}(f) \equiv 0 \pmod{5}$. Corollary \ref{fgh2} implies that the
products are permuted cyclically
according to whether $r$ is a quadratic residue modulo $5$. 
Each dissection is accomplished by expanding the product of
appropriate basis elements and using Lemma \ref{basis5a} to express
the result as
   a linear combination of quotients
    $\prod_{i=1}^{2}K_{5,i}^{a_{1}}\in M_{2}(\Gamma(5))$ and selecting those quotients with
    \begin{align*}
     5\cdot {\rm ord}_{i\infty} \left (\prod_{i=1}^{2}K_{5,i}^{a_{1}}
      \right) = 2 a_{1} + 3 a_{2} \equiv r \pmod{5}.
    \end{align*}

\begin{Corollary} \label{5s}
The following is a complete list of subscripts defining $f_{a_{0},a_{1},a_{2}}\in M_{2}(\Gamma_{1}(5))$ such
that for $0 \le r \le 4$, $$U_{5,r}(f_{a_{0},a_{1},a_{2}}) \equiv 0
\pmod{5}.$$ In particular, for 
$(a_{0}, a_{1},a_{2})$, we have
\begin{align}
  \label{eq:5}
  P_{a_{0},a_{1},a_{2}}(5n-\ell + r) \equiv 0 \pmod{5}.
\end{align}

\begin{center}
\begin{tabular}{|c|c|} \hline
  $r$ & $(a_{0}, a_{1},a_{2})$\\ \hline 
$0$ & $(4,-1,-1)$\\ \hline
  $1$ & $(4,4,-6), (4,5,5)$ \\\hline
  $2$ &  $(4,5,5)$ \\\hline
  $3$  &  $(4,5,5)$ \\\hline
  $4$ &  $(4,-6,4), (4,5,5)$ 
  \\ \hline
\end{tabular}
\end{center}

\end{Corollary}
The table in Theorem \ref{5s} belies additional Fourier coefficient
divisibility. Lemma \ref{basis5a} implies the following quotients involving Rogers-Ramanujan series $R(z) = \sum_{n=0}^{\infty} \frac{q^{n^{2}} z^{n}}{(q;q)_{n}}$ (see \cite{R3})
\begin{align*}
 q^{2} \sum_{n=0}^{\infty} P_{4,4,-6}(n)q^{n} = q^{2}\frac{R^{4}(1)}{R^{6}(q)} (q^{5}; q^{5})_{\infty}^{4}
\quad  \mbox{and}\quad \sum_{n=0}^{\infty} P_{4,-6,4}(n)q^{n} =q^{2} \frac{R^{4}(q)}{R^{6}(1)} (q^{5}; q^{5})_{\infty}^{4},
\end{align*}
satisfy congruences modulo $30$ witnessed by the relations
\begin{align*}
   \sum_{n=2}^{\infty} P_{4,4,-6}(5n-1)q^{(5n+1)/5} =
  \frac{30}{K_{5,2}} \quad\mbox{and}\quad  \sum_{n=1}^{\infty}
  P_{4,-6,4}(5n-1)q^{(5n-1)/5} =  \frac{30}{K_{5,1}}, 
\end{align*}
respectively. 
Corresponding congruences for higher weight forms may
be deduced by applying the dissection formulas from Lemma \ref{basis5a}.

We may similarly determine all quotients $f$ of the form \eqref{faaa} such that
$U_{7,r}(f) \equiv 0 \pmod{7}$. In particular, the smallest weight for
which such a congruence is satisfied is $2$ (i.e., $a_{0} =4$). As in
the level $5$ table, the algebraic symmetry of the
congruences implies that quotients satisfying congruences are permuted
cyclically and separated into equivalence classes
according to whether $r$ is a quadratic residue modulo $7$.

\begin{Theorem} \label{7s}
Let $\sigma=(1,2,3)\in S_{3}$. The following is a complete list of subscripts defining all 
quotients $f_{a_0,a_{1},a_{2},a_{3}}\in M_{k}(\Gamma_{1}(7))$,
$k=1,2,3$ such
that $$U_{7,r}(f_{a_0,a_{1},a_{2},a_{3}}) \equiv 0
\pmod{7}.$$ In particular, for 
$(a_{0},a_{1},a_{2},a_{3})$, we have
\begin{align}
  \label{eq:5}
  P_{a_{0},a_{1},a_{2},a_{3}}(7n-\ell + r) \equiv 0 \pmod{7}.
\end{align}

\begin{center}
\begin{tabular}{|c|c|} \hline
  $r$ & $(a_{0}, a_{1},a_{2}, a_{3})$\\ \hline 
$0$ & \text{See Theorem~$\ref{mainth1}$}  \\ \hline 
  $1$ & $(4,3, 0, -5), (6,-7, 3, 1), (6,-3, 5, -5), (6,-2, 2, -3), (6,0, 3, -6),
        (6,3, -6, 0),$  \\ & $(6,7, 7,7)$ \\ \hline 
  $2$ & $\{\sigma (a_{0},a_{1}, a_{2}, a_{3})) \mid
 U_{7,1}(f_{a_0,a_{1},a_{2},a_{3}}) \equiv 0
\pmod{7} \}$      
  \\ \hline 
  $3$  & $(4,2, 5, 3), (4,3, 2, 5), (4,5, 3, 2), (4,6, -2, -6), (6,-6, 7, -4), (6,-3, -9, 9)$, \\
  & $(6,0, -4, 1), (6,3, 1, -7), (6,3, 3, 3), (6,4, -2, \
         -5), (6,5, -5, -3), (6,6, -1, -8),$ \\
  &  $(6,7, -4, -6), (6,7, 7, 7), (6,9,
    -3, -9), (4,0, 2, -4)$ \\\hline
  $4$ & $\{ \sigma^{2}(a_{0},a_{1}, a_{2}, a_{3}) )\mid 
U_{7,1}(f_{a_0,a_{1},a_{2},a_{3}}) \equiv 0
\pmod{7} \}$   
   \\ \hline 
  $5$ & 
  $\{ \sigma^{2} (a_{0},a_{1}, a_{2}, a_{3})) \mid
U_{7,3}(f_{a_0,a_{1},a_{2},a_{3}}) \equiv 0
\pmod{7} \}$    
  \\ \hline 
  $6$ & 
 $\{\sigma (a_{0},a_{1}, a_{2}, a_{3})) \mid
U_{7,3}(f_{a_0,a_{1},a_{2},a_{3}}) \equiv 0
\pmod{7} \}$   
  \\ \hline
\end{tabular}
\end{center}
\end{Theorem}

\begin{proof}
These follow from the basis representations for $M_{k}(\Gamma_{1}(7))$ exhibited in Corollary \ref{tf} as well as Lemmas~\ref{fgh}, \ref{decomp} and~\ref{kr7}.
\end{proof}

We have seen in Theorem~\ref{mainth1} that $f_{6,3,3,3}(\tau)$
satisfies a Ramanujan type congruence modulo~$7$. The
quotient $f_{6,3,3,3}(\tau)$ agrees with the eta product $\eta(\tau)^{3}\eta(7\tau)^{3}=\sum_{n=1}^{\infty}a(n)q^{n}$. In fact, it is known \cite{Sc} that $\eta(\tau)^{3}\eta(7\tau)^{3}$ is the unique CM modular form of weight~3 over $\mathbb{Q}$ with CM by $K=\mathbb{Q}(\sqrt{-7})$ originally defined via
$$
\sum_{\substack{\mathfrak{a}\in I(\mathcal{O}_{K})\\ integral }}\psi_{2}(\mathfrak{a})q^{{\rm N}(\mathfrak{a})},
$$
where $\psi_{l}$ is the Hecke character of trivial conductor defined
in \cite[Lemma 4.1]{Sc}. This alternative representation yields
certain arithmetic information on $a(n)=P_{6,3,3,3}(n-1)$. For
example, it is clear by the Hecke theta series that $a(n)=0$ for $n$
such that $\left(\frac{n}{7}\right)=-1$.  The dissection made possible
by Lemma~\ref{decomp} allows us to prove this directly.

\begin{Corollary}
\label{ee7}
Let $a(n)$ be defined as above. Then for $n\equiv3,5,6\mod{7}$, $a(n)=0$, that is, $P_{6,3,3,3}(n-1)=0$.
\end{Corollary}

\begin{proof}
These follow from the fact
$$
\eta(\tau)^{3}\eta(7\tau)^{3}=x_{2} x_{0}^2-6 x_{2}^2x_{0}-7x_{1} x_{2}^2-x_{1}^2 x_{2}, \quad x_{n} = f_{\sigma_{7}^{n}(2,-2,0,1)}.
$$
After replacing $q$ by $q^{1/7}$ on the right side of the identity above, expanding through Lemma~\ref{decomp}, and simplifying through Lemma~\ref{kr7}, the result follows.
\end{proof}

By applying the dissection
  identities, selecting terms with index in certain nonzero
  congruence classes, and simplifying using Klein's relations, the
  following explicit representation for applications of $U_{p,r}$ result. Similar identities may be derived for permutations of the subscripts.
\begin{Corollary} \label{wt3}
Let $P_{a_{0},\ldots, a_{5}}(n)$ be defined by Definition~$\ref{defcw}$. Then 
  \begin{align*}
  & \sum _{n=0}^{\infty} P_{6,-4,1,0,0,0}(11^{j+1}n+11^{j} \cdot 10
   -1)q^{(11n+10)/11}  = \\
  &  11^{2j+1}\Bigg (  \frac{K_{11,4}^3 K_{11,5}^4 K_{11,3}^5}{K_{11,1}^{10} K_{11,2}^5}+\frac{2 K_{11,4}^2
       K_{11,5}^4 K_{11,3}^5}{K_{11,1}^{14}}-\frac{11 K_{11,2}^2 K_{11,5}^5
       K_{11,3}^5}{K_{11,1}^{15}}-\frac{3 K_{11,5}^4 K_{11,3}^5}{K_{11,1}^{11} K_{11,2}}-\frac{5 K_{11,4}
       K_{11,5}^5 K_{11,3}^5}{K_{11,1}^{11} K_{11,2}^3}\\
       &\qquad\qquad-\frac{14 K_{11,2}^4 K_{11,5}^4
       K_{11,3}^5}{K_{11,1}^{15} K_{11,4}}+\frac{3 K_{11,2}^2 K_{11,4}^2 K_{11,5}^5
       K_{11,3}^4}{K_{11,1}^{16}}+\frac{6 K_{11,4}^2 K_{11,5}^4 K_{11,3}^4}{K_{11,1}^{12} K_{11,2}}
  +\frac{17 K_{11,2}^4 K_{11,4} K_{11,5}^4 K_{11,3}^4}{K_{11,1}^{16}}\\
  &\qquad\qquad-\frac{3 K_{11,2}^4 K_{11,4}^3 K_{11,5}^4 K_{11,3}^3}{K_{11,1}^{17}}+\frac{6 K_{11,2}^3 K_{11,4}^3 K_{11,5}^4 K_{11,3}^2}{K_{11,1}^{15}}-\frac{3 K_{11,2}^4 K_{11,5}^5 K_{11,3}^2}{K_{11,1}^{13} K_{11,4}}\\
  &\qquad\qquad+\frac{15 K_{11,2}^4 K_{11,4}^3 K_{11,5}^2}{K_{11,1}^{10} K_{11,3}^2}-\frac{11 K_{11,2}^4 K_{11,4}^3 K_{11,5}^5}{K_{11,1}^{15}}-\frac{3 K_{11,2}^3 K_{11,4}^3 K_{11,5}^5}{K_{11,1}^{13} K_{11,3}}-\frac{3 K_{11,2}^4 K_{11,4}^2 K_{11,5}^5}{K_{11,1}^9 K_{11,3}^5} \Bigg).
  \end{align*}
\end{Corollary}
\begin{proof}
 Since $f_{6,-4,1,0,0,0}(\tau)$ is an
    eigenfunction for $U_{11}$ with eigenvalue $11^{2}$,
    \begin{align*}
(U_{11,10}\circ U_{11}^{j})(f_{6,-4,1,0,0,0} ) 
 &= 11^{2j}U_{11,10}( x_3 
    x_0^2+2 x_1 x_3 x_0+x_1 x_2 x_4),   
    \end{align*}
where $x_{n}=f_{\sigma_{11}^{n}(2,-1, -1,0,1,0)}$. The image under $U_{11,10}$ may be obtained by replacing $q$ by $q^{1/11}$ in $x_3 
    x_0^2+2 x_1 x_3 x_0+x_1 x_2 x_4$ and expanding the resulting expression via Lemma \ref{decomp11} as
   a linear combination of quotients
    $\prod_{i=1}^{5}K_{11,i}^{a_{1}}$. By selecting those quotients from each term with
    \begin{align*}
     11\cdot {\rm ord}_{i\infty} \left (\prod_{i=1}^{5}K_{11,i}^{a_{1}} \right) = 5 a_1+9 a_2+12 a_3+14 a_4+15 a_5 \equiv 10 \pmod{11}
    \end{align*}
and simplifying through Lemma~\ref{Klein11}, the claimed identity follows. 
\end{proof}

We next list all quotients $f$ up to weight $6$ of the form \eqref{faaa} such that
$U_{11,r}(f) \equiv 0 \pmod{11}$. In particular, the smallest weight for
which such a congruence is satisfied is $3$ (i.e., $a_{0} =6$). We
list the congruences corresponding to $r\equiv 1$ and $r\equiv 2$
modulo $11$ and apply the symmetries in Lemma \ref{fgh} to list
the remaining classes. 

\begin{Theorem}
\label{residuer}
Let $\sigma = (1,2,4,3,5) \in  S_{5}$. The following is a complete list of subscripts defining all 
forms $f_{a_{0},a_{1},a_{2},a_{3}a,_{4},a_{5}} \in
M_{k}(\Gamma_{1}(11))$ for $1\le k \le 3$ such
that $$U_{11,r}(f_{a_{0},a_{1},a_{2},a_{3}a,_{4},a_{5}}) \equiv 0
\pmod{11}.$$ In particular, for $1\le a_{0}/2 \le 3$, we have
\begin{align}
  \label{eq:5}
  P_{a_{0},a_{1},a_{2},a_{3},a_{4},a_{5}}(11n-\ell + r) \equiv 0 \pmod{11}.
\end{align}

\begin{center}
\begin{longtable}{|c|c|} \hline 
  $r$ & $(a_{0}, a_{1},a_{2},a_{3},a_{4},a_{5})$\\ \hline 
 $0$ & See Theorems~$\ref{mainth1}$ and \ref{chi11}\\ \hline
  $1$ & $(6,-3, 7, -4, -1, -2), (6,-2, 2, -5, 0, 2), (6, 1, 4, 0, 2,
         2), (6,7, -2, -4, -3, -1) $  
\\ \hline
  $2$ &$(6, -3, 6, -2, -1, -3), (6, 0, -7, 0, 8, -4), (6, 0, -4, 4, 0, -3), 
(6, 0, 0, -4, 0, 1)$ \\ & $(6, 1, -6, 3, 5, -6), (6, 2, 3, 3, 0, 1), (6, 3, -4, -5, -3, 6), (6, 3, -4, 0, 4, -6)$ \\ &   $(6, 3, 2, 1, 2, 1), (6, 3, 4, 1, 1, 0)$
\\ \hline
  $3$ & $\{\sigma (a_{0},a_{1}, a_{2}, a_{3}, a_{4}, a_{5}))
        \mid
U_{11,1}(f_{a_{0},a_{1},a_{2},a_{3}a,_{4},a_{5}}) \equiv 0
\pmod{11} \}$
  \\
  \hline 
  $4$ & $\{ \sigma^{4} (a_{0},a_{1}, a_{2}, a_{3}, a_{4}, a_{5}) )\mid
U_{11,1}(f_{a_{0},a_{1},a_{2},a_{3}a,_{4},a_{5}}) \equiv 0
\pmod{11} \}$
  \\
  \hline
  $5$ & $\{ \sigma^{3} (a_{0},a_{1}, a_{2}, a_{3}, a_{4}, a_{5}))
        \mid
U_{11,1}(f_{a_{0},a_{1},a_{2},a_{3}a,_{4},a_{5}}) \equiv 0
\pmod{11}\}$
  \\ \hline
  $6$ &  $\{ \sigma (a_{0},a_{1}, a_{2}, a_{3}, a_{4}, a_{5}))
        \mid
U_{11,2}(f_{a_{0},a_{1},a_{2},a_{3}a,_{4},a_{5}}) \equiv 0
\pmod{11}\}$
  \\
  \hline
  $7$ &  $\{ \sigma^{2} (a_{0},a_{1}, a_{2}, a_{3}, a_{4}, a_{5}))
        \mid
U_{11,2}(f_{a_{0},a_{1},a_{2},a_{3}a,_{4},a_{5}}) \equiv 0
\pmod{11}\}$
  \\ \hline
  $8$ &  $\{ \sigma^{4} (a_{0},a_{1}, a_{2}, a_{3}, a_{4}, a_{5})) \mid
U_{11,2}(f_{a_{0},a_{1},a_{2},a_{3}a,_{4},a_{5}}) \equiv 0
\pmod{11}\}$
  \\
  \hline
  $9$ &  $\{ \sigma^{2} (a_{0},a_{1}, a_{2}, a_{3}, a_{4}, a_{5})) \mid
U_{11,1}(f_{a_{0},a_{1},a_{2},a_{3}a,_{4},a_{5}}) \equiv 0
\pmod{11}\}$
  \\ \hline
  $10$ &  $\{ \sigma^{3} (a_{0},a_{1}, a_{2}, a_{3}, a_{4}, a_{5})
         )\mid U_{11,2}(f_{a_{0},a_{1},a_{2},a_{3}a,_{4},a_{5}}) \equiv 0 \pmod{11} \}$ \\ \hline
\end{longtable}
\end{center}
\end{Theorem}
\begin{proof}
These follow from the basis provided for $M_{k}(\Gamma_{1}(11))$ in the proof of Corollary \ref{tf} as well as Lemmas~\ref{fgh}, \ref{decomp11}, and~\ref{Klein11}.
\end{proof}

The unique modular forms over $\mathbb{Q}$ of weights~3 and~5 with CM by $K=\mathbb{Q}(\sqrt{-11})$ are given by
\begin{align*}
g_{11,3}(\tau)&=\sum_{\substack{\mathfrak{a}\in I(\mathcal{O}_{K})\\ integral }}\psi_{2}(\mathfrak{a})q^{{\rm N}(\mathfrak{a})} =\sum_{n=1}^{\infty}a_{11,3}(n)q^{n},
\end{align*}
and 
\begin{align*}
g_{11,5}(\tau)&=\sum_{\substack{\mathfrak{a}\in I(\mathcal{O}_{K})\\ integral }}\psi_{4}(\mathfrak{a})q^{{\rm N}(\mathfrak{a})} =\sum_{n=1}^{\infty}a_{11,5}(n)q^{n},
\end{align*}
where $\psi_{l}$ is the Hecke character of trivial conductor defined in \cite[Lemma 4.1]{Sc}. From the definitions, it can be deduced that for $n$ such that $\left(\frac{n}{11}\right)=-1$, $a_{11,3}(n)=a_{11,5}(n)=0$. We may re-derive these results using our dissection lemma, and moreover show that both of these CM modular forms satisfy Ramanujan type congruence modulo~11.

\begin{Corollary}
\label{a1135}
Let $a_{11,3}(n)$ and $a_{11,5}(n)$ be defined as above. Then the following hold.
\begin{enumerate}
\item for $n$ such that $\left(\frac{n}{11}\right)=-1$, $a_{11,3}(n)=a_{11,5}(n)=0$. 

\item $a_{11,3}(11^{j}n)\equiv0\pmod{11^{j}}$ and $a_{11,5}(11^{j}n)\equiv0\pmod{11^{2j}}.$

\end{enumerate}
\end{Corollary}

\begin{proof}
By employing the basis for $M_{5}(\Gamma_{1}(11))$ from the proof of Corollary \ref{tf}, one can derive basis representations whose dissections through Lemma~\ref{decomp11} imply the first claim. The Ramanujan type congruences in the second claim may be deduced upon using the basis representations to show that $g_{11,3} \in E_{3,-11}$ and $g_{11,5}\in E_{5,11^{2}}$.
\end{proof}

Lemma~\ref{decomp11}, implies congruences for the newform $\eta(\tau)^{2}\eta(11\tau)^{2}=\sum_{n=1}^{\infty}P_{4,2,2,2,2,2}(n-1)q^{n}$.

\begin{Corollary}
\label{ue11}
Let $U_{11,r}$ be defined by \eqref{defu}. 
\begin{enumerate}
\item If $r$ is a nonquadratic residue modulo~$11$, then
\begin{align*}
      U_{11,r}(\eta_{1}^{2}\eta_{11}^{2}) &\equiv 0 \pmod{2}.
\end{align*}

\item If $r$ is a nontrivial quadratic residue modulo~$11$, then
$$
U_{11,r}(\eta_{1}^{4}\eta_{11}^{4}) \equiv 0 \pmod{2}.
$$
\end{enumerate}
\end{Corollary}

Corollary~\ref{ue11} may be used to study the parity of the number of solutions of an elliptic curve of conductor~11 over finite fields.  Let $E_{11}$ denote the elliptic curve of conductor~11 defined by $y^{2}+y=x^{3}-x^{2}-10x-20$. One of the important arithmetic information on the elliptic curve $E_{11}$ is the number $|E_{11}(\mathbb{F}_{q})|$ of solutions of $E_{11}$ over the finite field $\mathbb{F}_{q}$ for a prime $q$. 


\begin{Corollary}
\label{e11q}
For any odd prime $q$ with $\left(\frac{q}{11}\right)=-1$, $|E_{11}(\mathbb{F}_{q})|$ is odd.
\end{Corollary}

\begin{proof}
Write $\eta(\tau)^{2}\eta(11\tau)^{2}=\sum_{n=1}^{\infty}a(n)q^{n}$, which is known \cite{MO} to be the newform of weight~2 and level $\Gamma_{0}(11)$ associated to $E_{11}$, and thus by the Modularity Theorem, for any prime $q$,
$$
|E_{11}(\mathbb{F}_{q})|=q-a(q).
$$
By Corollary~\ref{ue11} (1), for $r$ such that $\left(\frac{r}{11}\right)=-1$,
$$
 U_{11,r}(\eta_{1}^{2}\eta_{11}^{2}) \equiv 0 \pmod{2}.
 $$
 Therefore, $a(11n+r)\equiv0\pmod{2}$, and the corollary follows.
\end{proof}

For higher levels, Lemma~\ref{pdiss} may be applied to derive similar expansions.   
\begin{Theorem} \label{rdiss23}
The following congruences hold:
\begin{align*}
  U_{13,6}(f_{4,-6, 6, -2, 0, 0, 0}) &\equiv  U_{13,7}(f_{4,-6, 6, -2, 0, 0, 0}) \equiv 0 \pmod{2}, \\
    U_{17,6}(f_{4,-6, 6, -2, 0, 0, 0, 0, 0}) &\equiv  U_{17,11}(f_{4,-6, 6, -2, 0, 0, 0, 0, 0}) \equiv 0 \pmod{2}, \\
    U_{19,6}(f_{4,-6, 6, -2, 0, 0, 0, 0, 0,0}) &\equiv  U_{19,13}(f_{4,-6, 6, -2, 0, 0, 0, 0, 0,0}) \equiv 0 \pmod{2},\\
     U_{13,3}(f_{4,-2, 0, -2, -1, 3, 0}) &\equiv  U_{13,10}(f_{4,-2, 0, -2, -1, 3, 0}) \equiv 0 \pmod{3}, \\
      U_{17,8}(f_{4,-6, 6, -2, 0, 0, 0, 0, 0}) &\equiv  U_{17,9}(f_{4,-6, 6, -2, 0, 0, 0, 0, 0}) \equiv 0 \pmod{3}, \\
       U_{19,5}(f_{4,-4, 1, 1, 1, -1, 0, 0, 0, 0}  ) &\equiv   U_{19,14}(f_{4,-4, 1, 1, 1, -1, 0, 0, 0, 0}  ) \equiv 0 \pmod{3}, \\
   U_{13,5}(f_{4,-4, 1, 1, 1, -1, 0}) &\equiv  U_{13,8}(f_{4,-4, 1, 1, 1, -1, 0}) \equiv 0 \pmod{4}, \\
   U_{13,4}(f_{4,-5, 3, 1, -1, 0, 0}) &\equiv  U_{13,9}(f_{4,-5, 3, 1, -1, 0, 0}) \equiv 0 \pmod{6}, \\
      U_{17,4}(f_{4,-1,0,0,-5,1,0,0,3}) &\equiv  U_{17,13}(f_{4,4,-1,0,0,-5,1,0,0,3}) \equiv 0 \pmod{6}.
      \end{align*}
\end{Theorem}
\begin{proof}
A constructive proof through Lemma \ref{pdiss} is provided here for the first congruence. The others may be proven similarly. With $x_{n} =f_{\sigma_{13}^{n}(2,-1, 0, -1, 0, 0, 1)}$, basis expansions from the proof of Corollary \ref{tf}  imply \begin{align}
    f_{4,-6,6,-2,0,0,0}(\tau) = x_{0}^2+4 x_{0} x_{1}+2 x_{0}x_{3}+2x_{0} x_{4}+x_{1}^2+x_{2}^2+x_{4}^2+2 x_{1} x_{3}.\label{brep1} 
\end{align}
Apply Lemma \ref{pdiss} to the right side of the \eqref{brep1} and extract terms of the form $\prod_{i=1}^{6}K_{13,i}^{c_{i}}$, such that $-6 c_{1} - 11 c_{2}- 15 c_{3}- 
 18 c_{4}- 20 c_{5} - 21 c_{6} \equiv 6 \pmod{13}.$
This results in an expansion
\begin{align*}
    U_{13,6}(f_{4,-6,6,-2,0,0,0}) &= 2 \left ( \frac{K_{13,3}^3}{K_{13,1} K_{13,2} K_{13,5}^3}+\frac{2 K_{13,4}}{K_{13,2}^2 K_{13,3}}+\frac{16 K_{13,2}}{K_{13,1}^2 K_{13,4}}+\frac{13 K_{13,1}}{K_{13,3} K_{13,4}^2} + \cdots \right ).
\end{align*}
 The claimed congruence follows.  
\end{proof}

The congruences derived here constitute a small number of those satisfied with respect to both prime and composite moduli for products $f_{\vec{a}} \in M_{k}(\Gamma_{1}(p))$. These congruences can be found mining the lattice \eqref{ineq} for products potentially satisfying congruences and proven through dissections or by using the appropriate Sturm bound. Several examples are listed in the next theorem corresponding to modular forms of level $13$ and $19$.
\begin{Theorem}
  \begin{align*}
    U_{13,4}(f_{6,-8,6,0,-1,0,0}) \equiv  U_{13,10}(f_{6,-8,6,0,-1,0,0}) \equiv U_{13,3}(f_{4,-4, 2, -1, 1, 1, -1}) &\equiv 0 \pmod{13}, \\
     U_{13,7}(f_{4,3, 0, 0, 1, 0, 0}) \equiv 0 \pmod{26}, \qquad 
       U_{19,18}(f_{2,-3, 0, 0, 0, 3, -2, 0, 0, 0}) &\equiv 0 \pmod{19}.
  \end{align*}
\end{Theorem}

We conclude this section with an
example that illustrates how the dissection procedure  can be extended to levels that
are composite multiples of the prime levels considered here. The
construction is made possible by the fact that $\Gamma(N)$ is normal in ${\rm
  SL}_{2}(\mathbb{Z})$, and
$f(\tau/\ell)\in M_{k}(\Gamma(\ell p))$ if $f(\tau)\in \Gamma(p)$. The example is
illustrative of the general procedure for obtaining decomposible bases
for composite levels. For multiples of $7$, one must consider
$3$-dissections for at least $\dim M_{1}(\Gamma(21)) =96$ series, and
for level 11, one must formulate 5-dissections. These cases turn out to be
much more technical than level $10$, but still computationally feasible.

In \cite{Ga}, Garvan establishes the following congruences for certain $p$-core partition functions.
 \begin{Theorem}[Garvan]
 \label{gar}
 If $p\in\{5,7,11,13,17,19,23\}$, then for each prime divisor $\ell$ of $(p-1)/2$, 
 $$
 P_{(p-1),(-1)_{(p-1)/2}}\left(n-\frac{p^{2}-1}{24}\right)\equiv0\pmod{\ell}
 $$
 whenever $\left(\frac{n}{p}\right)=\epsilon_{p}$ and $\ell\nmid n$, $\epsilon_{11}=\epsilon_{13}=1$, and  $\epsilon_{5}=\epsilon_{7}=\epsilon_{17}=\epsilon_{19}=\epsilon_{23}=-1$.
 \end{Theorem}
One can  discover Garvan's congruences by extending the dissection method
presented in this work to composite levels. We illustrate this with a formulation of expansions for relevant modular forms of level $10$ that prove the case $p=5$ in Theorem~\ref{gar}, which states in our notation
$$
U_{10,3}(f_{4,-1,-1})\equiv U_{10,7}(f_{4,-1,-1})\equiv0\pmod{2}.
$$

Recall that by Lemma~\ref{basis5a}, 
$$
\{u_{j}\}_{0\leq j\leq 5}=\left\{ K_{5,1}^{j-3}K_{5,2}^{2-j}\right\}_{0\leq j\leq 5}
$$
is a basis for $M_{1}(\Gamma(5))$, and $f_{4,-1,-1}=u_{0}u_{5}$. Now let
\begin{align*}
w_{1}&=\frac{1}{2}\left(u_{0}\left(\frac{\tau}{2}\right)+u_{0}\left(\frac{\tau+1}{2}\right)\right),\quad w_{2}=\frac{1}{2}\left(u_{0}\left(\frac{\tau}{2}\right)-u_{0}\left(\frac{\tau+1}{2}\right)\right),\\
w_{3}&=\frac{1}{2}\left(u_{5}\left(\frac{\tau}{2}\right)+u_{5}\left(\frac{\tau+1}{2}\right)\right),\quad w_{4}=\frac{1}{2}\left(u_{5}\left(\frac{\tau}{2}\right)-u_{5}\left(\frac{\tau+1}{2}\right)\right),\\
w_{5}&=\frac{1}{2}\left(u_{1}\left(\frac{\tau}{2}\right)+\zeta_{10}^{-1}u_{1}\left(\frac{\tau+1}{2}\right)\right),\quad w_{6}=\frac{1}{2}\left(u_{1}\left(\frac{\tau}{2}\right)-\zeta_{10}^{-1}u_{1}\left(\frac{\tau+1}{2}\right)\right),\\
w_{7}&=\frac{1}{2}\left(u_{2}\left(\frac{\tau}{2}\right)+\zeta_{10}^{-2}u_{2}\left(\frac{\tau+1}{2}\right)\right),\quad w_{8}=\frac{1}{2}\left(u_{2}\left(\frac{\tau}{2}\right)-\zeta_{10}^{-2}u_{2}\left(\frac{\tau+1}{2}\right)\right),\\
w_{9}&=\frac{1}{2}\left(u_{3}\left(\frac{\tau}{2}\right)+\zeta_{10}^{-3}u_{3}\left(\frac{\tau+1}{2}\right)\right),\quad w_{10}=\frac{1}{2}\left(u_{3}\left(\frac{\tau}{2}\right)-\zeta_{10}^{-3}u_{3}\left(\frac{\tau+1}{2}\right)\right),\\
w_{11}&=\frac{1}{2}\left(u_{4}\left(\frac{\tau}{2}\right)+\zeta_{10}^{-4}u_{4}\left(\frac{\tau+1}{2}\right)\right),\quad w_{12}=\frac{1}{2}\left(u_{4}\left(\frac{\tau}{2}\right)-\zeta_{10}^{-4}u_{4}\left(\frac{\tau+1}{2}\right)\right).
\end{align*}
One can easily show that
$$
w_{1},w_{3}\in\mathbb{Z}[[q]],\quad w_{2},w_{4}\in q^{\frac{5}{10}}\mathbb{Z}[[q]],\quad w_{5}\in q^{\frac{1}{10}}\mathbb{Z}[[q]],\quad w_{6}\in q^{\frac{6}{10}}\mathbb{Z}[[q]],\quad w_{7}\in q^{\frac{2}{10}}\mathbb{Z}[[q]],
$$
$$
w_{8}\in q^{\frac{7}{10}}\mathbb{Z}[[q]],\quad w_{9}\in q^{\frac{3}{10}}\mathbb{Z}[[q]],\quad w_{10}\in q^{\frac{8}{10}}\mathbb{Z}[[q]],\quad w_{11}\in q^{\frac{4}{10}}\mathbb{Z}[[q]],\quad w_{12}\in q^{\frac{9}{10}}\mathbb{Z}[[q]],
$$
and verify that $\{w_{j}\}_{1\leq j\leq 12}$ is a basis for $M_{1}(\Gamma(10))$. Thus,
\begin{align*}
u_{0}(\tau/10)&=w_{1}+w_{2}+w_{5}+w_{6}-2w_{7}-2w_{8}+4w_{9}+4w_{10}-3w_{11}-3w_{12},\\
u_{5}(\tau/10)&=w_{3}+w_{4}+3w_{5}+3w_{6}+4w_{7}+4w_{8}+2w_{9}+2w_{10}+w_{11}+w_{12},
\end{align*}
and therefore, one deduces that
\begin{align*}
U_{10,3}(f_{4,-1,-1})&=2\left(w_{1}w_{9}+w_{2}w_{10}+2w_{4}w_{10}-3w_{11}w_{12}+2w_{3}w_{9}-w_{5}w_{7}-w_{6}w_{8}\right),\\
U_{10,7}(f_{4,-1,-1})&=2\left(2w_{1}w_{8}-w_{10}w_{12}-w_{9}w_{11}+2w_{2}w_{7}-w_{3}w_{8}-w_{4}w_{7}+3w_{5}w_{6}\right),
\end{align*}
which prove Garvan's congruences for $p=5$.

\section{Higher Levels, Orbits, and Chimeral Ramanujan Type Congruences}
This section addresses congruences for products under consideration modulo primes larger than those considered so far. We also characterize the orbit classes of products under the permutations $\sigma_{p}$ and comment on chimeral Ramanujan type congruences for higher prime levels. We begin with a proof of Theorem~\ref{61004}.
 

\begin{proof}[Proof of Theorem~\ref{61004}]
A proof is given for the first claim that is representative of each of the other cases. For the only if-part, assume that $U_{p}(f_{\vec{a}_{p}})=p^{2}f_{\vec{a}_{p}}$. By Lemma~\ref{fgh}, one deduces that
$$
U_{p}(f_{\sigma_{p}^{j}(\vec{a}_{p})})=\pm p^{2}f_{\sigma_{p}^{j}(\vec{a}_{p})},
$$
and thus, $U_{p}(f_{\sigma_{p}^{j}(\vec{a}_{p})})=\pm p^{2}q^{\ell(\sigma_{p}^{j}(\vec{a}_{p}))}+O(q^{\ell(\sigma_{p}^{j}(\vec{a}_{p}))+1})$. Therefore, by the definition of $U_{p}$, one concludes that for $s<\ell(\sigma_{p}^{j}(\vec{a}_{p}))$ and $j\leq\frac{p-1}{2}$ the coefficients $b(ps)$ of $f_{\sigma_{p}^{j}(\vec{a}_{p})}$  must all vanish.

For the if-part, first of all, simple computations show that the $p$-th coefficient of $f_{\vec{a}_{p}}$ is $p^{2}$. Then $h_{p}(\tau)=U_{p}(f_{\vec{a}_{p}})/p^{2}f_{\vec{a}_{p}}=1+O(q)$. Now assume that $b(ps)=0$ for any $s<\ell(\sigma_{p}^{j}(\vec{a}_{p}))$ and $2\leq j\leq\frac{p-1}{2}$. Since $f_{\vec{a}_{p}}\in M_{3}(\Gamma_{1}(p))$ by Lemma~\ref{m1r}, and $U_{p}$ is a linear operator on $M_{3}(\Gamma_{1}(p))$, $U_{p}(f_{\vec{a}_{p}})\in M_{3}(\Gamma_{1}(p))$, and thus, the orders of vanishing of $U_{p}(f_{\vec{a}_{p}})$ at any point of $X(\Gamma_{1}(p))$ are all non-negative. On the other hand, by the definition of $f_{\vec{a}_{p}}$ and Lemma~\ref{klf}, one can check that the orders of vanishing of $f_{\vec{a}_{p}}$ at points of $X(\Gamma_{1}(p))$ except for the cusps $r/p$ for $1\leq r\leq\frac{p-1}{2}$ are all zero. Therefore, the orders of vanishing of $h_{p}(\tau)$ at points of $X(\Gamma_{1}(p))$ except for the cusps $r/p$ for $2\leq r\leq\frac{p-1}{2}$ are all non-negative. For any $2\leq r\leq \frac{p-1}{2}$, suppose $r=\tilde{\alpha}^{j} \in (\mathbb{Z}/p\mathbb{Z})^{\times}/\{\pm1\} = \langle \alpha \rangle$, where $\tilde{\alpha} \equiv \alpha^{-1} \pmod{p}$. Then by Lemma~\ref{fgh},
$$
{\rm ord}_{r/p}(U_{p}(f_{\vec{a}_{p}}))={\rm ord}_{i\infty}(U_{p}(f_{\vec{a}_{p}})|_{a_{0}/2}^{j}\gamma_{p})={\rm ord}_{i\infty}(U_{p}(f_{\sigma_{p}^{j}(\vec{a}_{p})}))\geq \ell_{j}={\rm ord}_{r/p}(f_{\vec{a}_{p}}).
$$
Therefore, $h_{p}(\tau)$ is a holomorphic function on $X(\Gamma_{1}(p))$, which must be the constant~1. Hence, $U_{p}(f_{\vec{a}_{p}})=p^{2}f_{\vec{a}_{p}}$.
\end{proof}

\begin{Corollary} \label{hl} For all positive integers $n$ and $j$ and primes $p$ in the indicated range, \label{square}
\begin{align*}
      P_{6,1,0\ldots,0,-4}(p^{j}n-\ell) &\equiv 0  \pmod{p^{2j}}, \quad 5 \le p \le 101,
\end{align*}
and 
\begin{align*}
    P_{4,-2,0,\ldots,0,1,-2,1}(p^{j}n-\ell) &\equiv 0 \pmod{p^{j}}, \quad 11\le p\le 101, \\
    P_{4,-2,0,\ldots,0,1,-1}(p^{j}n-\ell) &\equiv 0 \pmod{p^{j}}, \quad 11\le p\le 101,   \\
    P_{4,1,1,0,\ldots,0,-2,-2}(p^{j}n-\ell) &\equiv 0 \pmod{p^{j}}, \quad 11\le p\le 101.
\end{align*}
\end{Corollary}

\begin{proof}[Proof of Corollary~\ref{hl}]
This follows from showing that the necessary and sufficient conditions given in Theorem~\ref{61004} are all met for the cases considered
\end{proof}

According to Theorem~\ref{z1}, products satisfying Ramanujan type congruences modulo $p$, fall into equivalence classes determined by the permutation $\sigma_p$. For the primes $p=5,7,11$, each product satisfying a Ramanujan type congruence is mapped to an
additional $(p-3)/2$ distinct products satisfying
Ramanujan type congruences. Thus, for $5\le p \le 11$, the total number of products $f_{\vec{a}} \in M_{k}(\Gamma_{1}(p))$ satisfying Ramanujan type
congruences is a multiple of $(p-1)/2$. This follows from the fact that the parameter $(p-1)/2$ is prime for $p=5,7,11$. For these primes, the
  Stabilizer Lemma and Lagrange's Theorem restrict the orbit sizes to $(p-1)/2$. This is not the
  case at level $p=13$, since this
  is the smallest prime for which $(p-1)/2=6$ is composite. For
  $p=13$, the possible sizes of
  the orbits under the corresponding group action are $2$, $3$ and
  $6$. 

In general, we may deduce the sizes of equivalence
  classes of products satisfying Ramanujan type congruences
  according to the form of the exponent list. Specifically, since subgroups of $\langle\sigma_{p}\rangle$ are all cyclic,  for each divisor $m$ of $(p-1)/2$, if $\langle\sigma_{p}^{m}\rangle$ is the stabilizer subgroup of $x=(a_{1},\ldots,a_{(p-1)/2})$ in $\langle\sigma_{p}\rangle$, then the orbit of $x$ under the action of $\langle\sigma_{p}\rangle$ is given by $\{\sigma_{p}^{\ell}(x)\}$ for $\ell\in \mathbb{Z}/m\mathbb{Z}$. More specifically, for a nontrivial subgroup $\langle\sigma_{p}^{m}\rangle$, writing $d=(p-1)/(2m)$, one has the following expansion as a product of disjoint cycles:
$$
\sigma_{p}^{m}=(1,\alpha^{m},\ldots,\alpha^{(d-1)m})(\alpha,\alpha^{m+1},\ldots,\alpha^{(d-1)m+1})\cdots(\alpha^{m-1},\alpha^{2m-1},\ldots,\alpha^{(d-1)m+m-1}).
$$
\begin{Corollary}
 Suppose $m \mid (p-1)/2$. If $\langle\sigma_{p}^{m}\rangle$ is the stabilizer of $x=(a_{1},\ldots,a_{(p-1)/2})$, then $a_{i}=a_{j}$, if and only if $i\equiv j\ \left({\rm mod}^{\times}\,\,{\alpha^{m}}\right)$; that is, $i=\alpha^{km}j$, for some $k\in \mathbb{N}$.
\end{Corollary}

This allows us to enumerate the number of products in an orbit according to the form of their vector of exponents. For instance, the fact that there are exactly $15$ products for the case $p=13$ and $a_{0}=4$ in Theorem~\ref{mainth1} follows from the fact that the orbits of each of the exponents $(4,1, 1, 0, 0,
-2, -2)$ and $(4,1,0,0,-1,-2,0)$
under action by $\langle (1,2,4,5,3,6) \rangle$ are of size $6$, while
the order of the orbit of $(4,1, -2, -2, 0, 1, 0)$ is $3$, since clearly, $\langle\sigma_{13}\rangle$ has two nontrivial subgroups,  $\langle\sigma_{13}^{2}\rangle=\langle(1,4,3)(2,5,6)\rangle$ and $\langle\sigma_{13}^{3}\rangle=\langle(1,5)(2,3)(4,6)\rangle$, for which the former is the stabilizer subgroup of $x=(a,b,a,a,b,b)$, and the latter is the stabilizer subgroup of $x=(a,b,b,c,a,c)$.


Recall from Definition \ref{defcong} that Ramanujan type congruences that are satisfied for only finitely many powers of $p$ are called chimeral. We have verified that there are no chimeral Ramanujan type congruences for quotients of Klein forms with exponents satisfying \eqref{ineq} for the primes $5, 7$ up to weight $3$. However, there are instances of products satisfying chimeral congruences of higher weights for these primes. 
For primes $p\ge 13$, there are many examples of products satisfying chimeral congruences up to weight $3$. As is evident in the next theorems, products satisfying chimeral congruences modulo $p$ for low weights exhibit patterns paralleling those for non-chimeral Ramanujan type congruences. In particular, common forms at lower levels appear to determine corresponding forms at higher prime levels. The recurring patterns are exhibited in the next theorems for products satisfying chimeral Ramanujan type congruences of weight $3$ and level $11\le p \le 19$. These lists share elements of the same form, with the defining exponents differing only by intermediate zeros.

\begin{Theorem}[Level 11, Weight 3, Chimeral] \label{chi11}
The following and their images under $\sigma^{i}$, $0\le i \le 4$ constitute a complete list of exponents for level $11$ weight $3$
forms satisfying  chimeral Ramanujan type congruences. All products
satisfy chimeral congruences modulo $11$:
 \begin{align*}
&  (6,-5, 3, 0, -2, 1), (6,-4, 1, 1, -3, 2), (6,-1, -1, -3, 2, 0), (6,1, 0, 1, -3, -2),\\
& (6,-3, 0, 3, -2, -1), (6,-3, 2, -1, -2, 1), (6, -1, 5, 1, 2, 2), (6,-5, 4, -3, 1, 0).         
  \end{align*}
\end{Theorem}

\begin{Theorem}[Level 13, Weight 3, Chimeral]
The following and their images under $\sigma^{i}$, $0\le i \le 5$ constitute a complete list of exponents for level $13$ weight $3$
forms satisfying  chimeral Ramanujan type congruences. All products
satisfy chimeral congruences modulo $13$:
\begin{align*}
    &(6,-5,3,0,-2,1,0),(6,-4,1,0,1,-3,2),(6,-1, -1, -3, 2, 0, 0),
    (6,3, -2, 1, 0, -3, -2), \\ &(6,1, 0, 0, 1, -3, -2),(6,-3,2,-1,-2,0,1),  (6,5, 2, 2, 1, 1, -2), (6,-5,4,-3,1,0,0), \\  &(6,-4,2,-2,1,-1,1),  (6,-3,0,-1,-1,4,-2),
    (6,-2,-2,1,3,-2,-1), (6,-2,-2,2,1,-2,0), \\ & (6,-2,1,-1,1,-1,-1), (6,-7,7,-3,2,-3,1).
\end{align*}
\end{Theorem}

 \begin{Theorem}[Level 17, Weight 3, Chimeral]
The following and their images under $\sigma^{i}$, $0\le i \le 7$ constitute a complete list of exponents for level $17$ weight $3$
forms satisfying  chimeral Ramanujan type congruences. All products
satisfy chimeral congruences modulo $17$:
\begin{align*}
 & (6,-5, 3, 0, -2, 1, 0, 0, 0), (6,-4, 1, 0,
  0, 0, 1, -3, 2), (6,-1, -1, -3, 2, 0, 0, 0, 0), \\ 
& (6,3,-2,0,0,1,0,-3,-2), (6,1, 0, 0, 0, 0, 1, -3, -2), (6,-3, 2, -1, -2, 0, 1, 0, 0), \\ & (6,5, 2, 2, 1, 1, -2, 0, 0), (6,-5, 4, -3, 1, 0, 0, 0, 0).
\end{align*}
 \end{Theorem}
 \begin{Theorem}[Level 19, Weight 3, Chimeral]
 The following and their images under $\sigma^{i}$, $0\le i \le 8$ constitute a complete list of exponents for level $19$ weight $3$
forms satisfying  chimeral Ramanujan type congruences. All products
satisfy chimeral congruences modulo $19$:
 \begin{align*}
   &  (6,-5,3,0,-2,1,0,0,0,0), (6,-4,1,0,0,0,0,1,-3,2),(6,-1, -1, -3, 2, 0, 0, 0, 0, 0), \\ &(6,3,-2,0,0,0,1,0,-3,-2),(6,1, 0, 0, 0, 0, 0, 1, -3, -2), (6,-3,2,-1,-2,0,1,0,0,0) \\ & (6,5, 2, 2, 1, 1, -2, 0, 0, 0),(6,-5,4,-3,1,0,0,0,0,0).
 \end{align*}
 \end{Theorem}


The products satisfying chimeral Ramanujan type congruences in the last theorems lead to the products of higher level satisfying congruences. Computations indicate that chimeral congruences for these exponent forms hold for all larger primes.
\begin{Conjecture}
\label{hl1} Denote by $\vec{\bf{0}}_{k}$ the zero vector of length $k$. For primes $p$ in the indicated range, $U_{p}(f_{a_0, a_1, \ldots, a_{(p-1)/2}}) \equiv 0 \pmod{p}$ if  $( a_{0}, \ldots, a_{(p-1)/2})$ lies in 
  \begin{align*}
 \left \{ \begin{array}{c}  (6,-5, 3, 0, -2, 1,\vec{\bf{0}}_{(p-11)/2}), (6, 2, -4, -3, 1, 1, \vec{\bf{0}}_{(p-11)/2}),  \\ (6,-1,-1,-3,2, \vec{\bf{0}}_{(p-9)/2}), (6,1,  \vec{\bf{0}}_{(p-9)/2}, 1, -3, -2), \\ (6,-5,4,-3,1, \vec{\bf{0}}_{(p-9)/2})
                    \end{array}  \right \},\quad  p\ge 11,
        \end{align*} 
  \begin{align*}
 \left \{ \begin{array}{c}  (6,3,-2,\vec{\bf{0}}_{(p-13)/2},1,0,-3,-2), (6, -3, 2, -1, -2, 0,1, \vec{\bf{0}}_{(p-13)/2}), \\ (6,5,2,2,1,1,-2,\vec{\bf{0}}_{(p-13)/2})
                    \end{array}  \right \},\quad  p\ge 13.
        \end{align*} 
\end{Conjecture}

 The search for complete lists of chimeral and non-chimeral Ramanujan type congruences for $P_{a_{0},\ldots,a_{(p-1)/2}}(n-\ell)$ can be extended to larger $a_{0}$ by relying on Theorem~\ref{z1} and Algorithm~\ref{a1} discussed in the next section. 
Our computations indicate there is no non-chimeral Ramanujan type congruence modulo $13$ for $P_{a_{0},\ldots,a_{6}}(n-\ell)$ for $8\leq a_{0}~\leq~16$. 
However, the search becomes less computationally feasible since the size of the associated bounded polyhedron increases as $p$ and $a_{0}$ get larger \cite[Theorem~5.3]{BP}. Products of higher weight satisfying chimeral congruences abound. For example, the product $f_{12,-13,17,-9,4,-7,2}\in M_{6}(\Gamma_{1}(13))$ and its orbit under $\sigma_{13}$ satisfy a chimeral Ramanujan type congruence modulo~$13^{4}$.

\section{Discussion of the Algorithm and computations}

This section includes a brief discussion of the computations and algorithms used.
To generate all relevant modular forms $f_{\vec{a}}$ for each level and weight, a full set of lattice points satisfying Lemma \ref{m1r} must be derived. Mining these lattice points for potential congruences is one of the most time-consuming and memory-intensive parts of the process. For larger primes and weights, the analysis required a distributed computing approach. The number of lattice points satisfying the inequality and congruence conditions of Lemma \ref{m1r} at each level and weight is given in Table 2. 
\begin{table}[h] 
    \centering 
    \begin{tabular}{c | ccccc} 
\toprule
   $p / a_{0}$ & $a_0 = 2$ & $ a_0 = 4$ &  $a_0 = 6$ &  $a_0 = 8$ &  $a_0 = 10$ \\
 \midrule
  5 & 2 &  4 & 6 & 9 & 12 \\
  7 &   6 &  18 & 39 & 72 & 120 \\ 
  11 &   25 &  226 & 1,000 & 3,126 & 7,877 \\ 
   13 &   42 & 684 & 4,388 & 17,976 & 56,076 \\ 
   17 &   96 & 4944 & 68,288 & 486,469 & 2,339,800 \\ 
   19 &  135 & 12,195 & 248,489 & 2,339,379 & 13,997,547 \\
\bottomrule
\end{tabular}
\caption{Number of lattice points satisfying \eqref{ineq}}
    \label{tab:my_label}
\end{table}

\begin{algorithm}[H]\label{a1}
\SetAlgoLined
\KwResult{
For a given $a_{0} \in 2\Bbb Z^{+}$, return a complete set $R_{p}$ of $(a_{1}, a_{2}, \ldots, a_{(p-1)/2})$ such that $f_{a_{0}, a_{1}, \ldots, a_{(p-1)/2}}(\tau)\in  
   {M}_{a_{0}/2}(\Gamma_{1}(p))$ satisfies a Ramanujan type
   congruence modulo $p$.  Return sets $C_{p}, C_{p^{2}},
   \ldots C_{p^{r}}$ corresponding to chimeral congruences modulo
   $p^{i}$.} \ \\ 
\textbf{Initialize:}

\Indp \Indp $L(a_{0}) =\{(a_{0}, a_{1}, \ldots, a_{(p-1)/2}) \mid$ System
\eqref{ineq} is satisfied$\}$\; 
$d =  \dim  {M}_{a_{0}/2}(\Gamma_{1}(p))$\;
$B_{p} =(b_{i})_{i=1}^{d},$  an ordered
basis for $ {M}_{a_{0}/2}(\Gamma_{1}(p))$\; $A_{p}$
the matrix of $U_{p}$ with respect to $B_{p}$\;
$V_{p} = \{ v \mid A_{p}v = \lambda v$, $\lambda = p^{r}$
 or $p^{r}\delta$, where $r\ge 1$, $\delta$ is algebraic integer\}
\\
\Indm \Indm
\vspace{0.1in}
\textit{\textbf{Step 1.} Form set $S_{p}$ of candidate lattice points satisfying Ramanujan type congruences.} \\
\For{$\ell \in L(a_{0})$}{
  set $F_{\ell,p}=(c_{np})_{0\le n\le 2}$, where $f_{\ell} = \sum_{n=0}^{\infty} c_{n}q^{n}$\;
  \If{$F_{\ell,p} \cdot F_{\ell,p} \equiv 0 \pmod{p}$}{
   $S_{p} \gets S_{p} \cup \{\ell\}$ $\; 
  $\; 
  }}
  \vspace{0.1in}
 \textit{\textbf{Step 2.} Form subset $S_{p}^{(1)}\subseteq S_{p}$ such that the orbit of the candidate lattice points are in $S_{p}$.} \\
  \For{$\ell \in S_{p}$}{
  \If{$\sigma_{p}^{j}(\ell) \in S_{p},\ \text{for}\ 1\le j \le (p-3)/2$}{$S_{p}^{(1)} \gets S_{p}^{(1)} \cup \{\ell\}$ }
  }
\vspace{0.1in}
\textit{\textbf{Step 3.} Compute eigendecompositions \& chimeral sets.} \\
\For{$\ell \in S_{p}^{(1)}$}{set $v_{f_{\ell}} = (v_{k})_{1\le k \le d}$, where $f_{\ell}(\tau) =
  \sum_{k=1}^{d} v_{k} b_{k} $\; 
    \eIf{$v_{f_{\ell}}$ is in the span of $V_{p}$ and satisfies the hypotheses of Lemma~\ref{upeee}}{
  $R_{p} \gets R_{p} \cup \{\ell \}$\;}{$i = 1$, set $F_{\ell,p}=(c_{np})_{0\le n\le 5}$, where $f_{\ell} = \sum_{n=0}^{\infty} c_{n}q^{n}$ \;
\While{$F_{\ell,p^{i}}\cdot F_{\ell,p^{i}} \equiv 0 \pmod{p^{i}}$}{
$i \gets i+1$}
$C_{p^{i}} \gets C_{p^{i}} \cup \{\ell\}$
}
  }
 \caption{Find all $f_{a_{0}, \ldots,
       a_{\frac{p-1}{2}}}(\tau)\in  {M}_{{a_{0}}/{2}}(\Gamma_{1}(p))$
   satisfying Ramanujan type congruences}
\end{algorithm}


An optimization of the computation of the initial series coefficients in Step 1 of the algorithm is important. 
When the data set of lattice points satisfying Lemma \ref{m1r} is very large, computing even the first few coefficients of each series required a lot of time. We significantly decreased the run-time by using symbolic expressions for the coefficients of the general product $f_{\vec{a}} = \sum c_{n}(\vec{a})q^{n}$ in terms of undetermined exponents $\vec{a}$: $c_{0}(\vec{a}), c_{p}(\vec{a}), c_{2p}(\vec{a})$. Testing each lattice point $\ell = \vec{a}$
\begin{align} \label{to} c_{0}(\vec{a}) \equiv c_{p}(\vec{a}) \equiv c_{2p}(\vec{a}) \equiv 0 \pmod{p} \end{align} 
with symbolic representations was much more efficient than computing the series expansion up to order $q^{2p}$ for each product. 
This refinement and the elimination of orbits of lattice points not satisfying Ramanujan type congruences significantly decreased the number of points we need to further sift. These preliminary steps allowed us to eliminate most of the lattice points before the comparatively computationally expensive basis representations are constructed and used to find eigenrepresentations or corresponding expansions to show non-congruences for chimeral cases.

\end{document}